\documentclass[12pt,a4paper]{article}

\usepackage{amsmath,amssymb,amsthm,color}
\usepackage{verbatim}
\usepackage{graphicx,epstopdf}
\usepackage{subfigure}
\usepackage{amscd}

\definecolor {green}{rgb}{0,0.4,0}
\usepackage{bbm,ascii,hyperref}
\usepackage{mysec,times1}

\usepackage{cancel}
\usepackage[normalem]{ulem}

\newcommand{\BP}{\mathbf{P}}
\newcommand{\BE}{\mathbf{E}}
\newcommand{\EE}{\mathbb{E}}
\newcommand{\PP}{\mathbb{P}}

\newcommand{\CP}{\varLambda}

\newcommand{\CalD}{\mathcal{D}}

\newcommand{\bq}{\mathfrak{q}}
\newcommand{\ii}{\mathrm{i}}
\newcommand{\rd}{\mathrm{d}}
\newcommand{\re}{\mathrm{e}}
\newcommand{\RR}{\mathbb{R}}
\newcommand{\NN}{\mathbb{N}}
\newcommand{\ZZ}{\mathbb{Z}}
\newcommand{\MDP}{\text{MDP}}
\newcommand{\bfOne}{\mathbbm{1}}
\newcommand{\DOT}{\raisebox{-.1pc}{\mbox{\Large\boldmath{$\cdot$}}}}
\newcommand{\subDOT}{\raisebox{-.08pc}{\mbox{\boldmath{$\cdot$}}}}
\newcommand{\const}{\mathrm{const}}

\newcommand{\Li}{\operatorname{Li}\nolimits}
\newcommand{\floor}[1]{\lfloor#1\rfloor}
\newcommand{\ceiling}[1]{\lceil#1\rceil}

\newcommand{\myp}{\mbox{$\:\!$}}
\newcommand{\mypp}{\mbox{$\;\!$}}

\newcommand{\myn}{\mbox{$\;\!\!$}}
\newcommand{\mynn}{\mbox{$\:\!\!$}}

\def\MR#1{\href{http://www.ams.org/mathscinet-getitem?mr=#1}{MR#1}}

\numberwithin{equation}{section}

\newtheorem{theorem}{Theorem}[section]
\newtheorem{proposition}[theorem]{Proposition}
\newtheorem{lemma}[theorem]{Lemma}
\newtheorem{corollary}[theorem]{Corollary}
\theoremstyle{definition}
\newtheorem{assumption}{Assumption}[section]
\newtheorem{definition}{Definition}[section]
\theoremstyle{remark}
\newtheorem{remark}{Remark}[section]

\newtheorem{example}{Example}[section]

\begin{document}

\title{Limit shape of minimal difference partitions\\
and fractional statistics}

\author{Leonid V.\ Bogachev\myp$^{\rm
a}$
%\thanks{Department of Statistics, School of Mathematics,
%University of Leeds, Leeds LS2 9JT, UK. Email:
%\href{mailto:L.V.Bogachev@leeds.ac.uk}{L.V.Bogachev@leeds.ac.uk}}\
and Yuri V.\ Yakubovich\myp$^{\rm b}$
%\thanks{Department of
%Probability Theory and Mathematical Statistics, Saint Petersburg
%State University, Universitetsky Prospekt 28, Peterhof, Saint
%Petersburg, 198504, Russia. Email:
%\href{mailto:Y.Yakubovich@spbu.ru}{Y.Yakubovich@spbu.ru}}
}
\date{\small $^{\,\rm a}$\myp Department of Statistics, School of Mathematics, University of Leeds, Leeds LS2 9JT,
UK.\\
Email:
\href{mailto:L.V.Bogachev@leeds.ac.uk}{L.V.Bogachev@leeds.ac.uk}\\[.2pc]
$^{\,\rm b}$\myp Department of Probability Theory and Mathematical
Statistics, St.~Petersburg State University,
%Universitetsky Prosp.\ 28, Peterhof, Saint Petersburg, 198504,
7/9~Universitetskaya nab., St.~Petersburg, 199034,
Russia.\\ Email:
\href{mailto:Y.Yakubovich@spbu.ru}{Y.Yakubovich@spbu.ru}}

%\date{17 September 2018}

\maketitle
\begin{abstract}
The class of \emph{minimal difference partitions} $\MDP(q)$ (with
gap $q$) is defined by the condition that successive parts in an
integer partition differ from one another by at least $q\ge0$. In a
recent series of papers by A.~Comtet and collaborators, the
$\MDP(q)$ ensemble with uniform measure was interpreted as a
combinatorial model for quantum systems with \emph{fractional
statistics}, that is, interpolating between the classical
Bose--Einstein ($q=0$) and Fermi--Dirac ($q=1$) cases. This was done
by formally allowing values $q\in(0,1)$ using an analytic
continuation of the limit shape of the corresponding Young diagrams
calculated for integer~$q$. To justify this ``replica-trick'', we
introduce a more general model based on a variable MDP-type
condition encoded by an integer sequence $\bq=(q_i)$, whereby the
(limiting) gap $q$ is naturally interpreted as the Ces\`aro mean of
$\bq$. In this model, we find the family of limit shapes
parameterized by $q\in[0,\infty)$ confirming the earlier answer, and
also obtain the asymptotics of the number of parts.


%\\
%05A17 Partitions of integers\\
%60C05 Combinatorial probability\\
%82B10 Quantum equilibrium statistical mechanics (general)

\end{abstract}

\section{Introduction}
\subsection{Integer partitions and the limit shape}\label{sec:1.1}
An \emph{integer partition} is a decomposition of a given natural
number into an \emph{unordered sum} of integers; for example,
$35=8+6+6+5+4+2+2+1+1$. That is to say, a non-increasing sequence of
integers $\lambda_1\ge\lambda_2\ge\dots\ge 0$,
$\lambda_i\in\NN_0:=\NN\cup \{0\}=\{0,1,2,\dots\}$ is a partition of
$n\in\NN_0$ if $n=\lambda_1+\lambda_2+\cdots$, which is expressed as
$\lambda\vdash n$. Zero terms are added as a matter of convenience,
without causing any confusion. The non-zero terms
$\lambda_i\in\lambda$ are called the \emph{parts} of the partition
$\lambda$. We formally allow the case $n=0$ represented by the
``empty'' partition $\varnothing=(0,0,\dots)$, with no parts. The
set of all partitions $\lambda\vdash n$ is denoted by $\CP(n)$, and
$\CP:=\cup_{n\in\NN_0}\CP(n)$ is the collection of \emph{all}
integer partitions. For a partition $\lambda=(\lambda_i)\in\CP$, the
sum $N(\lambda):=\lambda_1+\lambda_2+\cdots$ is referred to as its
\textit{weight} (i.e., $\lambda\vdash N(\lambda)$), and the number
of its parts
$K(\lambda):=\#\{\lambda_i\in\lambda\colon\lambda_i>0\}$ is called
the \textit{length} of~$\lambda$. Thus, for $\lambda\in\CP(n)$, we
have $N(\lambda)=n$ but $K(\lambda)\le n$.

\begin{figure}[h]
\begin{picture}(60,180) \put(0,8){
%\put(127,0){\includegraphics[width=8.27cm]{YoungDiagram2.pdf}}
\put(127,0){\includegraphics[width=8.27cm]{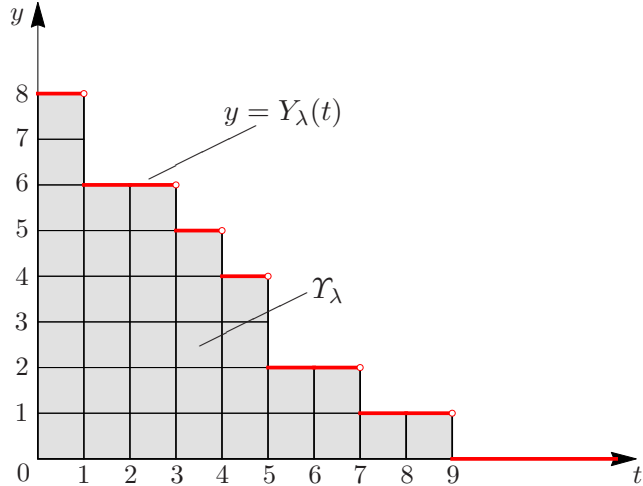}}
\put(123,172){\mbox{\footnotesize $y$}}
%\put(139,170){$u$}
\put(356,-3){\mbox{\footnotesize$t$}}
%\thicklines
\put(234,69){\line(-2,-1){40}}%
%\put(234,69){\vector(-2,-1){40}}%
\put(236,67){\mbox{\small $\varUpsilon_\lambda$}}%
%\put(215,132){\vector(-2,-1){40}}%
\put(215,132){\line(-2,-1){40}}%
\put(203,134){\mbox{\small $y=Y_\lambda(t)$}}
%\put(171,64){\mbox{\small $\varUpsilon_\lambda$}}
%\put(291,13){\tiny$D_9(\la){=}1\ge q_0{=}1$}
%\put(256,30.5){\tiny$D_7(\la){=}1\ge q_2{=}1$}
%\put(222,56){\tiny$D_5(\la){=}2\ge q_4{=}1$}
%\put(204.5,82.5){\tiny$D_4(\la){=}1\ge q_5{=}0$}
%\put(187,99.5){\tiny$D_3(\la){=}1\ge q_6{=}1$}
%\put(153,124.5){\tiny$D_1(\la){=}2{\ge} q_8{=}1$}
\put(125.5,-2.2){\mbox{\footnotesize$0$}}
%\put(217,-7){\mbox{\footnotesize$5$}}
%\put(300,-7){$10$}
\put(124.8,20.6){\mbox{\footnotesize$1$}}
%\put(126,37){\mbox{\footnotesize$2$}}
\put(125,37.8){\mbox{\footnotesize$2$}}
\put(124.8,55){\mbox{\footnotesize$3$}}
\put(125,72){\mbox{\footnotesize$4$}}
\put(125,90){\mbox{\footnotesize$5$}}
\put(125,107.3){\mbox{\footnotesize$6$}}
\put(125,124.3){\mbox{\footnotesize$7$}}
\put(125,141){\mbox{\footnotesize$8$}}
\put(148,-3.0){\mbox{\footnotesize$1$}}
\put(165.6,-3.0){\mbox{\footnotesize$2$}}
\put(183,-3.0){\mbox{\footnotesize$3$}}
\put(200,-3.0){\mbox{\footnotesize$4$}}
\put(217.5,-3.0){\mbox{\footnotesize$5$}}
\put(234.6,-3.0){\mbox{\footnotesize$6$}}
\put(251.5,-3.0){\mbox{\footnotesize$7$}}
\put(269.1,-3.0){\mbox{\footnotesize$8$}}
\put(286.3,-3.0){\mbox{\footnotesize$9$}} }
\end{picture}
\caption{The Young diagram $\varUpsilon_\lambda$ (shaded) of a
partition $\lambda=(8,6,6,5,4,\allowbreak 2,2,\allowbreak
1,1,\allowbreak 0,\dots)$, with weight $N(\lambda)=35$ and length
$K(\lambda)=9$. Note that the parts $\lambda_i>0$ are represented by
the successive columns of the diagram. The graph of the step
function $t\mapsto Y_\lambda(t)$ (shown in red in the online
version) gives the upper boundary of
$\varUpsilon_\lambda$.}\label{fig:yd}
\end{figure}

A partition $\lambda=(\lambda_1,\lambda_2,\dots)\in\CP$ is
succinctly visualized by its \emph{Young diagram}
$\varUpsilon_{\myn\lambda}$ formed by left- and bottom-aligned
column blocks with $\lambda_1,\lambda_2,\dots$ unit square cells,
respectively. In particular, the area of the Young diagram
$\varUpsilon_{\myn\lambda}$ equals the partition weight
$N(\lambda)$. The upper boundary of $\varUpsilon_{\myn\lambda}$ is a
non-increasing step function $Y_\lambda\colon[\myp0,\infty)\to\NN_0$
(see Fig.~\ref{fig:yd} for illustration). Note that
$\inf\{t\ge0\colon Y_\lambda(t)=0\}$ coincides with the length
$K(\lambda)$.

Theory of integer partitions is a classical branch of discrete
mathematics and combinatorics dating back to Euler, with further
fundamental contributions due to Hardy, Ramanujan, Rademacher and
many more (see \cite{Andrews} for a general background). The study
of asymptotic properties of \emph{random} integer partitions (under
the uniform distribution) was pioneered by Erd\H{o}s \& Lehner
\cite{ErdosLehner}, followed by a host of research which in
particular discovered a remarkable result that, under a suitable
rescaling, the Young diagrams
$\varUpsilon_\lambda$ of typical partitions $\lambda$ of a
large integer $n$ are close to a certain deterministic \textit{limit
shape}. For \emph{strict} partitions (i.e.,\ with distinct parts)
this result was (implicitly) contained already
in~\cite{ErdosLehner}; for
\emph{plain} partitions (i.e., without any
restrictions), the limit shape was first identified by Temperley
\cite{Temperley} in relation to the equilibrium shape of a growing
crystal, and obtained more rigorously much later by Vershik (as
pointed out at the end of \cite{VK1985}) using some asymptotic
estimates by Szalay \& Tur\'an \cite{ST}. An alternative proof in
its modern form was outlined by Vershik \cite{Vershik96} and
elaborated by Pittel~\cite{Pittel}, both using the conditioning
device\footnote{The randomization trick, often collectively called
``Poissonization'', is well known in the general enumerative
combinatorics (see, e.g., Kolchin et al.~\cite{Kolchin}). In the
context of integer partitions, it was introduced  by
Fristedt~\cite{Fristedt}.} based on a suitable randomization of the
integer $n$ being partitioned.

\begin{center}
\begin{figure}[h]
%\linethickness{0.3mm}
\thinlines
%\mbox{}\hspace{2.5pc}\includegraphics[width=2.50in]{fig2a.eps}\hspace{.5in}
%\includegraphics[width=2.50in]{fig2b.eps}
%
\mbox{}\hspace{2.5pc}\includegraphics[height=2.8in]{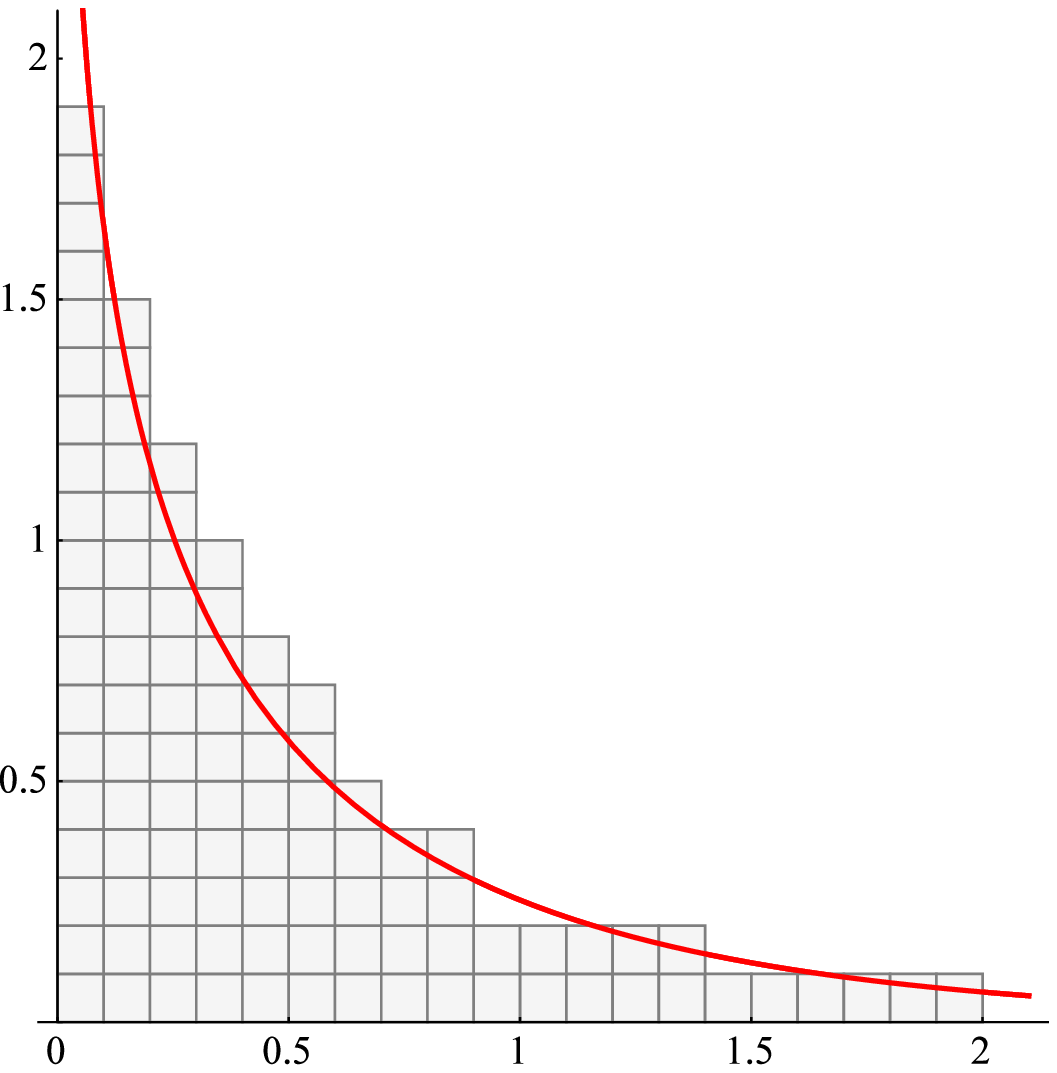}\hspace{2pc}
\includegraphics[height=2.8in%width=1.30in
]{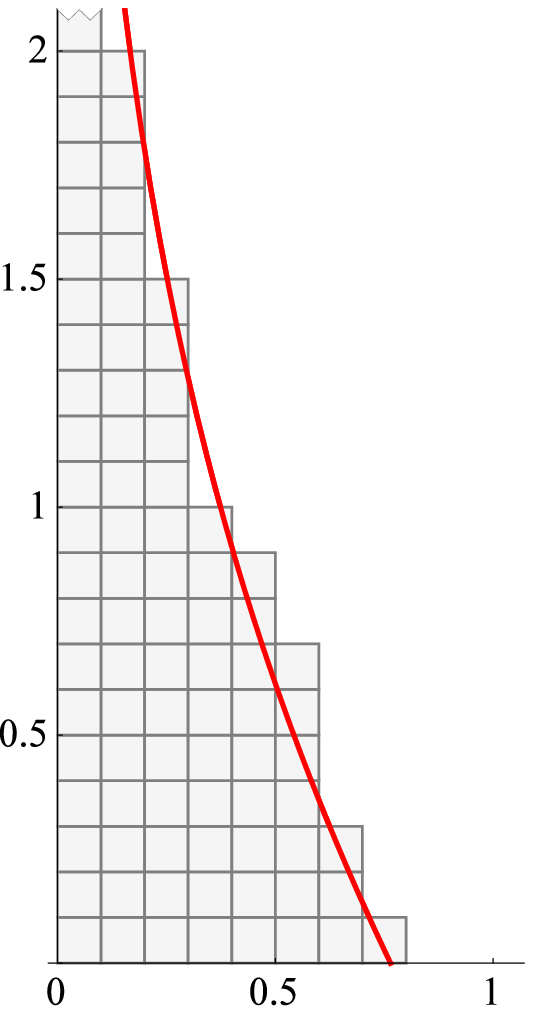}
\put(-270,100){\mbox{\small $\re^{-x\pi/\sqrt{6}}+\re^{-y\pi/\sqrt{6}}=1$}} %
\put(-40,100){\mbox{\small $\re^{\myp x\pi/\sqrt{12}}=\re^{-y\pi/\sqrt{12}}+1$}} %
%\put(-315,-6){\mbox{\footnotesize  (a)}} %
%\put(-97,-6){\mbox{\footnotesize  (b)}}  %
%H1
\put(-132,2.5){\mbox{\footnotesize$x$}}
%H2
\put(2,2.5){\mbox{\footnotesize$x$}}
%V1
\put(-328,205){\mbox{\footnotesize$y$}}
%V2
\put(-101,205){\mbox{\footnotesize$y$}} %
\caption{The limit shape (shown in red in the online version) for
plain partitions (left) and strict partitions (right), both under
the scaling $Y_\lambda(t)\mapsto n^{-1/2}\,Y(n^{1/2}\mypp t)$
\,($\lambda\vdash n$) as $n\to\infty$. The scaled Young diagrams
(shaded in grey) represent integer partitions uniformly sampled with
$n=100$. On the right picture, the largest part (depicted as the
leftmost column) is only partially shown; in fact, here
$\lambda_1=35$.} \label{fig1+}
\end{figure}
\end{center}

Under the natural rescaling of Young diagrams
$\varUpsilon_{\myn\lambda}$ of partitions $\lambda\vdash n$ by
$\sqrt{n}\myp$ in each coordinate,\footnote{See, however,
Remark~\ref{rm:1.6} below.} the limit shape for these two classical
ensembles is determined, \strut{}respectively, by the equations
$\re^{-x\pi/\sqrt{6}}+\re^{-y\pi/\sqrt{6}}=1$ (plain partitions) and
$\re^{\myp x\pi/\sqrt{12}}=\re^{-y\pi/\sqrt{12}}+1$ (strict
partitions); see Fig.~\ref{fig1+}. Note that in the latter case, the
limit shape hits zero at $x=c_1= \pi^{-1}\sqrt{12}\,\log
2\doteq0.764304$; this implies that the number of parts $K(\lambda)$
in a typical strict partition $\lambda\vdash n$ grows like
$c_1\sqrt{n}$ as $n\to\infty$. In contrast, for plain partitions the
number of parts grows faster than $\sqrt{n}$\myp; more precisely,
$K(\lambda)\sim c_0\myp\sqrt{n}\,\log n$, where
$c_0=\sqrt{6}/(2\pi)$~\cite{ErdosLehner}.

To date, many limit shape results are known for
integer partitions subject to various restrictions
(see, e.g.,
Bogachev~\cite{Bogachev}, Yakubovich~\cite{Yakubovich}, and also a
review in
 DeSalvo \& Pak \cite{Pak2}).
Deep connections between statistical properties of quantum systems
(where discrete random structures naturally arise due to
quantization) and asymptotic theory of random integer partitions are
discussed in a series of papers by Vershik
\cite{Vershik96,Vershik-UMN}. Note that the idea of conditioning in
problems of quantum statistical mechanics was earlier promoted by
Khinchin~\cite{Khinchin} who advocated systematic use of local limit
theorems of probability theory as a tool to prove the equivalence of
various statistical ensembles in the thermodynamic limit.

From the point of view of statistical mechanics, it is
conventional\footnote{For a historic background, see older papers by
Auluck \& Kothari \cite{Auluck} and Temperley \cite{Temperley}, and
Vershik~\cite{Vershik-UMN} for a modern exposition.} to interpret
the integer partition $\lambda=(\lambda_i)\in\CP$ as the energy
spectrum in a sample configuration (state) of quantum gas, with
$K(\lambda)=\#(\lambda_i>0)$ particles and the total energy
$\sum_i\lambda_i=N(\lambda)$. Note that decomposition into a sum of
integers is due to the quantization of energy in quantum mechanics,
while using \emph{unordered} partitions corresponds to the fact that
quantum particles are indistinguishable. In this context, the limit
shape of Young diagrams associated with random partitions (for
instance, under the uniform measure) is of physical interest as it
describes the asymptotic distribution of particles in such ensembles
over the energy domain.

\subsection{Minimal difference partitions}
For a given $q\in\NN_0$, the class of \textit{minimal difference
partitions with gap $q$}, denoted by $\MDP(q)$, is the set of
integer partitions $\lambda=(\lambda_1,\lambda_2,\dots)$ subject to
the restriction $\lambda_i-\lambda_{i+1}\ge q$ whenever
$\lambda_i>0$. Two important special cases of the $\MDP(q)$ are
furnished by the values $q=0$ corresponding to \emph{plain
partitions} (i.e., with no restrictions), and $q=1$ leading to
\emph{strict partitions} (i.e., with different parts).

In this paper, we propose a natural generalization of the $\MDP$
property as follows.

\begin{definition}
For a given sequence $\bq=(q_i)_{i\in\NN_0}$ of non-negative
integers (with the convention that $q_0\ge 1$), we define
$\CP_{\bq}\equiv\MDP(\bq)$ to be the set of all integer partitions
$\lambda=(\lambda_i)$ subject to the variable MDP-type condition
\begin{equation}\label{eq:restr}
\lambda_i-\lambda_{i+1}\ge
q_{k-i},\qquad
i=1,\dots,k,
\end{equation}
where $k$ is the number of (non-zero) parts in the partition
$\lambda$. By convention, the empty partition $\varnothing$
satisfies~\eqref{eq:restr}. The sequence $\bq$ is referred to as the
\emph{gap sequence}.
\end{definition}

\begin{remark}
For $i=k$, the inequality \eqref{eq:restr} specializes to
$\lambda_{k}-\lambda_{k+1}\equiv \lambda_{k}\ge q_0$. That is to
say, the smallest part of the partition $\lambda=(\lambda_i)$ is
required to be not less than $q_0\ge1$ (which really poses a
restriction only if $q_0>1$).
\end{remark}

\begin{remark}\label{rm:1.2}
The partition model \eqref{eq:restr} appeared earlier (without any
name) in a paper by Bessenrodt and Pak \cite[\S\mypp4]{Bessenrodt}
devoted to partition bijections, in connection with
\emph{generalized Sylvester's transformation} $\lambda_{k-i}\mapsto
\lambda_{k-i}+\sum_{j=0}^i q_i$ \,($i=0,\dots,k-1$), extending the
classical case $q_i\equiv 1$.
\end{remark}

\begin{remark}
Alternatively, one could consider partitions subject to similar
restrictions as \eqref{eq:restr} but in the reverse order relative
to the sequence $\bq$,
\begin{equation*}
\lambda_i-\lambda_{i+1}\ge q_{i},\qquad
i=1,\dots,k.
\end{equation*}
However, the model \eqref{eq:restr} is preferable in view of the
physical interpretation of parts $\lambda_i$ as successive energy
levels in a configuration (state) of a quantum system
\cite{Vershik-UMN}, which makes it more natural to enumerate the
energy gaps starting from the minimal level
$\lambda_{k}=\min\{\lambda_i\colon \lambda_i>0\}$.
%One more
%motivation for choosing the model \eqref{eq:restr} is its earlier
%appearance (without any name) in a paper by Bessenrodt and Pak
%\cite[\S\,4]{Bessenrodt} (devoted to partition bijections) in
%connection with generalized Sylvester's transformation (where the
%classical.
\end{remark}

Throughout the paper, we impose the following
\begin{assumption}\label{as:main}
The gap sequence $\bq=(q_i)$ satisfies the asymptotic regularity
condition
\begin{equation}\label{eq:assq*}
Q_k:=\sum_{i=0}^{k-1} q_i = q\myp k+O(k^{\beta}) \qquad
(k\to\infty),
\end{equation}
with some $q\ge 0$ and $0\le\beta<1$.
\end{assumption}

Note that under Assumption \ref{as:main} the sequence $\bq=(q_i)$
has a well-defined Ces\`aro mean, referred to as the \emph{limiting
gap},
\begin{equation}\label{eq:q*}
\lim_{k\to\infty}k^{-1} Q_k= q\ge 0.
\end{equation}

\begin{remark}
In the case $q=0$, the asymptotic relation \eqref{eq:assq*}
accommodates sequences $(Q_k)$ that are irregularly growing
(provided the growth is sublinear) or even bounded ($\beta=0$),
including the case $Q_k\equiv1$ corresponding to plain
(unrestricted) integer partitions.
\end{remark}

For $q=0$ (when the leading term in \eqref{eq:assq*} vanishes), it
is still possible to derive the limit shape results under our
standard Assumption~\ref{as:main}. However, to obtain the
asymptotics of the typical MDP length $K(\lambda)$,
more regularity should be assumed by
specifying the behaviour of the remainder term
$O(k^\beta)$.

\begin{assumption}[$q=0$]\label{as:main+}
The gap sequence $\bq=(q_i)$ satisfies the asymptotic regularity
condition
\begin{equation}\label{eq:assq*+}
Q_k:=\sum_{i=0}^{k-1} q_i = \tilde{q}\myp
k^{\beta}+O(k^{\tilde{\beta}}\myp) \qquad (k\to\infty),
\end{equation}
with some $\tilde{q}\ge 0$ and $0\le\tilde{\beta}<\beta<1$.
\end{assumption}

\begin{remark}\label{rem:1.4}
The utterly degenerate case $\tilde{q}=0$ and $\tilde{\beta}=0$ in
Assumption \ref{as:main+} is equivalent to Assumption \ref{as:main}
with $q=0$ and $\beta=0$. In this case, we have $Q_k=O(1)$ as
$k\to\infty$, and since $q_i\in\NN_0$, this implies that $q_i=0$ for
all sufficiently large $i$. Clearly, the first few non-zero terms in
the sequence $\bq=(q_i)$ (i.e., in the MDP conditions
\eqref{eq:restr}) do not affect any limiting results, and so
effectively such a model is identical with the classical case of
plain partitions ($q_0=1$ and
$q_i\equiv0$ for $i\in\NN$).
\end{remark}

\subsection{Main result}\label{sec:1.3}
For $n\in\NN_{0}$, consider the subset $\CP_{\bq}(n)=\CP_{\bq}\cap
\CP(n)$ comprising $\MDP(\bq)$ partitions of weight $N(\lambda)=n$.
For example, the partition $\lambda=(8,6,6,5,4,2,2,1,1,0,0,\dots)$
used in Fig.~\ref{fig:yd} fits into the \MDP-space $\CP_{\bq}(35)$
with the alternating sequence $\bq=(1,0,1,0,1,0,\dots)$.
%; also, $\lambda\in\CP_{\tilde{\bq}}$ with
%$\tilde{\bq}=(1,0,1,0,2,1,1,0,2,0,\dots)$.
Suppose that each (non-empty) space $\CP_{\bq}(n)$ is endowed with
uniform probability measure denoted
by $\nu_n^{\bq}$. We are interested in asymptotic
properties (as $n\to\infty$) of this and similar measures on MDP
spaces; in particular, we find the limit shape of properly scaled
Young diagrams associated with partitions $\lambda\in\CP_{\bq}(n)$
and prove exponential bounds for deviations from the limit shape.

Let us state one of our main results, slightly simplifying the
notation as compared to the more general case treated in
Section~\ref{sec:pdls}. For every $q\ge0$, define the function
\begin{equation}\label{eq:LS_intro}
\varphi(t;q):=\max\{0,-q\myp t-\log\myn(1-\re^{-t})\},\qquad t>0,
\end{equation}
and let $T_q:=\inf\{t>0\colon \varphi(t;q)=0\}$; that is, $T_q$ is
the unique root of the equation
\begin{equation}\label{eq:Tq-intro}
q=-T_q^{-1}\log (1-\re^{-T_q})
\end{equation}
(with the convention $T_0:=+\infty$). The area under the graph of
$\varphi(t;q)$ is computed as
\begin{equation}\label{eq:theta-intro}
\vartheta_q^2:=\int_0^{T_q} \!\varphi(t;q)\,\rd{t}=-\frac{q\myp
T_q^2}{2}+\Li_2(1)-\Li_2(\re^{-T_q}),
\end{equation}
where $\Li_2(\cdot)$ denotes the \emph{dilogarithm} function (see,
e.g., \cite[p.\mypp1]{Lewin}),
\begin{equation}\label{eq:Dilog}
\Li_2(x):=-\int_0^x \frac{\log\myn(1-u)}{u}\,\rd{u}\equiv
\sum_{k=1}^\infty\frac{x^k}{k^2},\qquad 0\le x\le 1.
\end{equation}
Note that $\Li_2(1)=\zeta(2)=\pi^2/6$. It is easy to check from
\eqref{eq:Tq-intro} that $\lim_{q\downarrow 0}q\myp T_q^2=0$, so
using \eqref{eq:theta-intro} we obtain
\begin{equation}\label{eq:q=0}
\vartheta_0=\lim_{q\downarrow0}\vartheta_q=\sqrt{\Li_2(1)}=\frac{\;\pi}{\sqrt{6}}\,.
\end{equation}

Finally, observe that, setting $x=\re^{-T_q}$ in the well-known
identity\footnote{This identity can be obtained from the definition
\eqref{eq:Dilog} by integration by parts.} \cite[Eq.\:(1.11),
p.\,5]{Lewin}
\begin{equation}\label{eq:Li2=}
\Li_2(x)+\Li_2(1-x)=\Li_2(1)-\log x\cdot \log\myn(1-x),
\end{equation}
and using the equation \eqref{eq:Tq-intro}, the expression
\eqref{eq:theta-intro} is rewritten in a more appealing form,
\begin{equation}\label{eq:theta-intro-alt}
\vartheta_q^2=\frac{q\myp T_q^2}{2}+\Li_2(1-\re^{-T_q}),
\end{equation}
where the terms on the right-hand side can be given a meaningful
geometric interpretation (see details in Section~\ref{sec:4.4}).

\begin{theorem}[Limit shape in $\CP_\bq(n)$]\label{th:main*}
Let the sequence $\bq=(q_i)$ satisfy Assumption
\textup{\ref{as:main}}, with $q\ge0$. Then, for every $t_0>0$ and
any $\varepsilon>0$, we have
\begin{equation}\label{eq:ls-intro}
\lim_{n\to\infty}\,\nu^\bq_n\Bigl\{ \lambda\in\CP_\bq(n)\colon
\sup_{t\ge t_0}\myp\bigl|n^{-1/2}\myp Y_\lambda(t\myp n^{1/2})
-\vartheta_q^{-1}\varphi(t\myp\vartheta_q;q)\bigr|>\varepsilon
\Bigr\} =0,
\end{equation} where $Y_\lambda(\cdot)$ denotes the upper boundary of the
Young diagram $\varUpsilon_\lambda$
and $\vartheta_q$ is given
by~\eqref{eq:theta-intro-alt}.
\end{theorem}

In view of formula \eqref{eq:LS_intro}, in the Cartesian coordinates
\begin{equation}\label{eq:cartesian}
x=t\myp\vartheta_q,\qquad y=\varphi(t\myp\vartheta_q;q)
\end{equation} the limit shape \eqref{eq:ls-intro} is given by the
equation
\begin{equation}\label{eq:ls-intro-xy}
\re^{-y}=\re^{qx}\myp(1-\re^{-x}).
\end{equation}
Clearly, $y=y(x)$ is a continuous decreasing function (as long
as $y(x)>0$), hitting zero at $x=T_q$ for $q>0$ (see
equation~\eqref{eq:Tq-intro}) and with $\lim_{x\to\infty}y(x)=0$ for
$q=0$.
\begin{remark}\label{rm:1.6}
It is common to scale Young diagrams via reducing their area $n$
to~$1$ \cite{Vershik-UMN}. In our case, this leads to the additional
rescaling in the expression of the limit shape
(see~\eqref{eq:ls-intro}). Instead, it is more natural to work with
the intrinsic scaling \eqref{eq:cartesian} to produce a simpler
equation for the limit shape \eqref{eq:ls-intro-xy} but where the
limiting area $\vartheta_q^2$ varies with $q$
(see~\eqref{eq:theta-intro-alt}). See the precise corresponding
assertions in Section~\ref{sec:pdls}.
\end{remark}

\begin{example}\label{ex:ls} Let us specialize the notation introduced before Theorem
\ref{th:main*} for a few simple values of $q\ge0$, including all
cases where closed expressions for $T_q$ and $\vartheta_q$ in
elementary functions are available.
\begin{itemize}
\item
$q=0$: here $T_0=\infty$,
$\vartheta_0=\sqrt{\Li_2(1)}=\pi/\sqrt{6}\doteq1.282550$, and the
limit shape \eqref{eq:ls-intro-xy} specializes to (cf.\
Vershik~\cite[p.\,99]{Vershik96})
$$
\re^{-x}+\re^{-y}=1.
$$

\item
$q=1$: from the equation \eqref{eq:Tq-intro} we get $T_1=\log 2
\doteq 0.693147$. By virtue of Euler's result
(see~\cite[Eq.\:(1.16), p.\,6]{Lewin})
$$
\Li_2\bigl(\tfrac12\bigr)=\frac{\pi^2}{12}-\frac{(\log 2)^2}2,
$$
we obtain from~\eqref{eq:theta-intro-alt}
$$
\vartheta_1^2=\frac{T_1^2}2
+\Li_2\bigl(\tfrac12\bigr)=\frac{\pi^2}{12}.
$$
Hence, $\vartheta_1=\pi/\sqrt{12}\doteq 0.906900$ and the limit
shape \eqref{eq:ls-intro-xy} is reduced to (cf.\
Vershik~\cite[p.\mypp100]{Vershik96})
$$
\re^{x}-\re^{-y}=1.
$$

\item
$q=2$: the equation \eqref{eq:Tq-intro} (quadratic in
$z=\re^{-T_2}$) solves to give
$T_{2}=\log\left(\frac{1+\sqrt{5}}{2}\right)\doteq 0.481212$. Hence,
we find $1-\re^{-T_2}=\tfrac{3-\sqrt{5}}{2}$. Using a known
expression for the dilogarithm at this point (see \cite[Eq.\:(1.20),
p.\mypp7]{Lewin}), we obtain from~\eqref{eq:theta-intro-alt}
\begin{align*}
\vartheta^2_{2}=T_2^2+\Li_2\biggl(\frac{3-\sqrt{5}}{2}\biggr)&=
\log^2\biggl(\frac{1+\sqrt{5}}{2}\biggr)+
\frac{\pi^2}{15}-\frac14\log^2\biggl(\frac{3-\sqrt{5}}{2}\biggr)=\frac{\pi^2}{15},
\end{align*}
which gives $\vartheta_2=\pi/\sqrt{15}\doteq 0.811156$ (cf.\
Romik~\cite{Romik}).

\item
$q=\frac12$: solving the equation \eqref{eq:Tq-intro} we get
$T_{1/2}=\log\left(\frac{3+\sqrt{5}}{2}\right)\doteq 0.962424$.
Hence, $1-\re^{-T_{1/2}}=\tfrac{\sqrt{5}-1}{2}$. Using another exact
value of dilogarithm \cite[Eq.\:(1.20), p.\mypp7]{Lewin}, formula
\eqref{eq:theta-intro-alt} yields
$$
\vartheta^2_{1/2}=\frac{T_{1/2}^2}{4}+\Li_2\biggl(\frac{\sqrt{5}-1}{2}\biggr)
=\frac14 \log^2\biggl(\frac{3+\sqrt{5}}{2}\biggr)+ \frac{\pi^2}{10}
-\log^2\biggl(\frac{\sqrt{5}-1}{2}\biggr)=\frac{\pi^2}{10},
$$
so that $\vartheta_{1/2}=\pi/\sqrt{10}\doteq 0.993459$.

\item
$q=3$: an exact value of $T_{3}$ can be found by solving the
equation \eqref{eq:Tq-intro} (cubic in $z=\re^{-T_3}$), but no
elementary expression is available for $\Li_2(1-\re^{-T_3})$
(cf.~\cite{Lewin}). It is easy to find numerically $T_{3}\doteq
0.382245$ and $\vartheta_{3}\doteq 0.752618$ (cf.\ \cite[Fig.\,3,
p.\,8]{Comtet2}).

\item
$q=\frac13$: numerical values are given by $T_{1/3}\doteq 1.146735$
and $\vartheta_{1/3}\doteq 1.038508$.

\end{itemize}
\end{example}

\subsection{MDP and fractional statistics}

The special case of the $\MDP(\bq)$ model with a \emph{constant} gap
sequence $q_i\equiv q\in\NN_0$ in \eqref{eq:restr} was considered in
a series of papers by Comtet et al.\ \cite{Comtet1,Comtet2,Comtet3}
in connection with fractional exclusion statistics of quantum
particle systems (see \cite{Khare}, \cite{Lerda} or \cite{Murthy}
for a ``physical'' introduction to this area). These authors
obtained the limit shape of $\MDP(q)$ using a physical
argumentation. In particular,  it was observed that the analytic
continuation of the limit shape, as a function of $q\in\NN_0$, into
the range $q\in(0,1)$ (the so-called \textit{replica trick}) may be
interpreted as a quantum gas obeying fractional exclusion
statistics, thus furnishing a family of probability measures
``interpolating'' between the Bose--Einstein statistics ($q=0$) and
the Fermi--Dirac statistics ($q=1$).

In the present work,\footnote{A short announcement of our approach
(in the case $q>0$) appeared in~\cite{BoYu}.} we provide a
combinatorial justification of this physical construction by working
with a more general $\MDP(\bq)$ model satisfying
Assumption~\ref{as:main}. In addition to many deterministic examples
with such a property, the assumption \eqref{eq:assq*} (and hence
\eqref{eq:q*}) holds almost surely for sequences of independent
random variables $\bq=(q_i)$ satisfying mild conditions, thus
providing a stochastic version of the $\MDP(\bq)$ model (see
Section~\ref{sec:MDPRE} below).

As was observed by Comtet et al.\ \cite{Comtet1}, another model of
statistical physics leading to the MDP-type constraint is the
one-dimensional \emph{quantum Calogero model} with harmonic
confinement (see \cite{Polychronakos} for a review and further
references therein), defined by the Hamiltonian of a $k$-particle
system with spatial positions $(x_i)_{i=1}^{k}$ on a line,
\[
H_q:=-\frac{1}{2}\sum_{i=1}^{k} \frac{\partial^2}{\partial x_i^2}
+\sum_{1\le i<j\le k}\frac{q\myp(q-1)}{(x_i-x_j)^2}
+\frac{1}{2}\sum_{i=1}^{k} x_i^2.
\]
This model is exactly solvable, and the solution can be expressed in
terms of the pseudo-excitation numbers $\lambda_i$ satisfying the
condition $\lambda_i-\lambda_{i+1}\ge q$, with a positive real~$q$.

%\begin{remark} At first glance, it might seem that $q=1$ is a
%special point since the singular part of the potential changes its
%sign there;
%however, this is not the case due to quantum effects.
%Likewise, we will see that the
%$\MDP(\bq)$ model has no singularity at $q=1$.
%%\begin{color}{blue}\sout{Similarly, our family of
%%combinatorial models with the restrictions \eqref{eq:restr} and the
%%regularity condition \eqref{eq:assq*} obviously has no singularity
%%at $q=1$, hence suggesting that nothing special should
%%happen in the Calogero model.}\end{color}
%\end{remark}

As is common in such models (cf.~\cite{Haldane}), an analogue of
Pauli's exclusion principle is not strictly local for models
$\MDP(\bq)$ with sequences $\bq=(q_i)$ not degenerating to the
trivial sequences $q_i\equiv 0$ or $q_i\equiv1$ ($i\in\NN$). Indeed,
the occurrence of part $\lambda_{i}=j$ rules out a few adjacent
values, that is, $\lambda_{i-1}\notin \{j,j+1,\dots,j+q_{k-i+1}-1\}$
if $q_{k-i+1}>0$ or $\lambda_{i+1}\notin
\{j,j-1,\dots,j-q_{k-i}+1\}$ if $q_{k-i}>0$, but the actual index
$k-i$ is determined by the entire partition $\lambda=(\lambda_i)$
through the rank of the part $\lambda_{i}=j$ among all (ordered)
parts $\lambda_i$, together with the total number $k$ of non-zero
parts in $\lambda$.

\begin{remark}
Heuristically, the requirement $\lambda_i-\lambda_{i+1}\ge q$ with
$q\in(0,1]$ may be interpreted, at least for integer $m:=q^{-1}$, as
saying that $\lambda_i-\lambda_{i+m}\ge 1$ as long as $\lambda_i>0$,
that is, to prohibit more than $m=q^{-1}$ equal parts; in other
words, no part counts bigger than $q^{-1}$ are allowed. For $q=1$
this indeed translates as only strict partitions being permissible.
In the general case, this interpretation turns out to be true for
the \emph{expected} part counts (see \cite[\S\mypp5.2]{Khare});
however, literal restriction that the part counts do not exceed
$q^{-1}$ leads to a different model called \textit{Gentile's
statistics} \cite[\S\mypp5.5]{Khare}. The limit shape of partitions
under Gentile's statistics was found in \cite[\S\mypp9]{Pak} (see
also~\cite{Yakubovich} where a rigorous proof is given).
\end{remark}

The rest of the paper is organized as follows. In Section
\ref{sec:defs}, several measures on minimal difference partitions
are introduced, and certain relations between them are stated.
Section~\ref{sec:m} is devoted to finding the typical length of
MDPs. In Section~\ref{sec:pdls} the main results concerning the
limit shape of MDPs, both with a restricted and unrestricted length
growth, are proved. In fact, we obtain
sharp exponential bounds for
deviations from the limit shape. Section~\ref{sec:lsalt} describes
an alternative approach to the limit shape based on a partition
bijection that effectively removes the MDP-constraint. In
Section~\ref{sec:MDPRE}, we extend our results to the case of random
sequences~$\bq$. Finally, the Appendix contains proof of the two
technical propositions stated in Section~\ref{sec:defs}, which
establish the equivalence of ensembles.

\section{Probability measures on the MDP spaces}
\label{sec:defs}

\subsection{Basic definitions and notation}
In this paper, we shall use several probability measures on MDPs and
other partition spaces. In the present section we describe them and
establish some properties. First we introduce notation for some
functionals on partitions we shall need. If one fixes a probability
measure on partitions, these functionals become random variables.

Let $\lambda=(\lambda_1,\lambda_2,\dots)$
($\lambda_1\ge\lambda_2\ge\dots\ge0$) be an integer partition,
$\lambda\in\CP$. Recall that $N(\lambda)=\lambda_1+\lambda_2+\cdots$
and $K(\lambda)=\#\{\lambda_i\in\lambda\colon\lambda_i>0\}$. An
equivalent description of a partition~$\lambda$ can be given in
terms of the consecutive differences
$D_j(\lambda)=\lambda_j-\lambda_{j+1}$; obviously,
\begin{equation}\label{eq:NK}
\lambda_i=\sum_{j\ge i} D_j(\lambda),\qquad N(\lambda)=\sum_{j\ge 1}
jD_j(\lambda),\qquad K(\lambda)=\max\{j\colon D_j(\lambda)>0\}.
\end{equation}
Consider the function
\begin{equation}\label{eq:Phi}
Y_\lambda(t):=\sum_{j>t} D_j(\lambda),\qquad t\ge 0,
\end{equation}
Clearly, the map $t\mapsto Y_\lambda(t)$ is non-increasing,
piecewise constant, and right-continuous. From \eqref{eq:NK}, it is
also easy to see that $Y_\lambda(t)=\lambda_{\lfloor t\rfloor+1}$
($t\ge0$), with $\lfloor\cdot\rfloor$ denoting the floor function
(i.e., integer part). The \textit{Young diagram}
$\varUpsilon_\lambda$ of a partition $\lambda$ is defined as the
closure of the planar set
\begin{equation*}
\{(t,u)\in\RR^2\colon \,t\ge0,\ 0\le u\le Y_\lambda(t)\}.
\end{equation*}
That is to say, the Young diagram $\varUpsilon_\lambda$ is the union
of (left- and bottom-aligned) column blocks with $\lambda_1$,
$\lambda_2$, \dots\ unit squares, respectively; in particular, the
function $t\mapsto Y_\lambda(t)$ defines its upper boundary
(cf.\ Section~\ref{sec:1.1}). We
shall often identify the Young diagram $\varUpsilon_\lambda$ with
the (graph of the) function $Y_\lambda(t)$ (see Fig.~\ref{fig:yd}).

The measure most important for us is the aforementioned uniform
measure $\nu^\bq_n$ on the set $\CP_\bq(n)$:
$$
\nu^\bq_n(\lambda):=\frac{1}{p_\bq(n)}\quad
(\lambda\in\CP_\bq(n)),\qquad p_\bq(n):=\#\CP_\bq(n).
$$
The space $\CP_\bq(n)$ can be further decomposed as a disjoint union
of the sets $\CP_\bq(n,k):=\{\lambda\in\CP_\bq(n)\colon
K(\lambda)=k\}$, and one can introduce the uniform measures on these
spaces,
$$
\nu^\bq_{n,\myp k}(\lambda):=\frac{1}{p_\bq(n,k)}\quad
(\lambda\in\CP_\bq(n,k)),\qquad p_\bq(n,k):=\#\CP_\bq(n,k).
$$
Note that $\nu^\bq_{n,\myp k}$ can be viewed as the measure
$\nu^\bq_n$ conditioned on the event $\{K(\lambda)=k\}$; indeed, for
any $\lambda\in \CP_\bq(n,k)$,
\begin{align*}
\nu^\bq_{n}(\lambda\mypp|\mypp
K(\lambda)=k)&=\frac{\nu^\bq_{n}(\lambda)}{\nu^\bq_{n}(K(\lambda)=k)}\\
&=\frac{1/p_\bq(n)}{\sum_{\lambda\in\CP_\bq(n,k)}1/p_\bq(n)}
=\frac{1}{p_\bq(n,k)}=\nu^\bq_{n,\myp k}(\lambda).
\end{align*}
This conditional measure is somewhat simpler than $\nu^\bq_n$
itself, since there exists a product expression for the Laplace
generating function of $p_\bq(n,k)$ with respect to $n$ (for any
fixed~$k$).

To establish such an expression, the following simple observation is
useful. Define
\[ \CalD_\bq(k):=\{(d_1,\dots,d_k)\in\NN_0^k\colon \,d_j\ge q_{k-j},\
j=1,\dots, k\},\qquad k\in\NN.
\]
Then the $\MDP(\bq)$ condition \eqref{eq:restr} implies that
$\lambda\in\CP_\bq(\DOT,k):=\bigcup_{n\ge 0}\CP_\bq(n,k)$ if and
only if $(D_1(\lambda),\dots,D_k(\lambda))\in\CalD_\bq(k)$ and
$D_j(\lambda)=0$ for all $j>k$. Hence, the space $\CP_\bq(\DOT,k)$
is in one-to-one correspondence with the set $\mathcal{D}_\bq(k)$.
Moreover, using the second of the formulas \eqref{eq:NK}, the
Laplace generating function $F_\bq(z,k)$ ($z\ge0$) of the sequence
$(p_\bq(n,k))_{n\ge0}$ (with $k\ge 0$ fixed) is evaluated as
$F_\bq(z,0)=1$ and for $k\ge 1$
\begin{align}
\notag F_\bq(z,k):=\sum_{n=0}^\infty
p_\bq(n,k)\,\re^{-zn} &=\sum_{n=0}^\infty
\sum_{\lambda\in\CP_\bq(\subDOT,\myp
k)}\!\bfOne_{\{N(\lambda)=n\}} \,\re^{-zN(\lambda)}\\
\notag
&=\sum_{\lambda\in\CP_\bq(\subDOT,\myp k)}\re^{-zN(\lambda)}
=\prod_{j=1}^k\sum_{d_j=\myp q_{k-j}}^\infty \re^{-z j\myp d_j}\\
&=\prod_{j=1}^k \frac{\re^{-z j\myp
q_{k-j}}}{1-\re^{-zj}}=\frac{\re^{-zs_k}}{(1-\re^{-z})
\cdots(1-\re^{-zk})}\,,\label{eq:gfpnm}
\end{align}
where we set
\begin{equation}\label{eq:s}
s_k:=\sum_{j=1}^{k}j\myp q_{k-j}\equiv \sum_{i=1}^{k} Q_i,\qquad
k\in\NN,
\end{equation}
with $Q_i$ defined in \eqref{eq:assq*}. In particular, $s_k\ge k$
for all $k\ge 1$ (because $Q_i\ge q_0\ge1$, see \eqref{eq:assq*});
moreover, the asymptotic condition \eqref{eq:assq*} implies that,
for $q\ge 0$,
\begin{equation}\label{eq:sksim}
s_k=\frac{q\myp k^{2}}{2}+O(k^{\beta+1})\qquad (k\to\infty).
\end{equation}

\begin{remark}
The product structure of $F_\bq(z,k)$ revealed in \eqref{eq:gfpnm}
is similar to that of multiplicative measures introduced by
Vershik~\cite{Vershik96}. However, there are some distinctions from
multiplicative measures. Firstly, the partition length $K(\lambda)$
must be fixed to obtain independence. Secondly, the role of the part
counts which become independent after randomization of
$N(\lambda)=n$ is played here by the differences $D_j(\lambda)$.
\end{remark}
Let us define an auxiliary probability measure $\mu^\bq_{z,\myp k}$
on the space $\CP_\bq(\DOT,k)$ (parameterized by $z>0$) by setting
\begin{equation}\label{eq:mu_zk}
\mu^\bq_{z,\myp
k}(\lambda):=\frac{\re^{-zN(\lambda)}}{F_\bq(z,k)},\qquad \lambda\in
\CP_\bq(\DOT,k).
\end{equation}
Note that, for every $z>0$, the measure $\mu^\bq_{z,\myp k}$
conditioned on the event $\{N(\lambda)=n\}$ coincides with the
uniform measure $\nu^\bq_{n,\myp k}$ on the space $\CP_\bq(n,k)$;
indeed, according to \eqref{eq:mu_zk} we have, for any
$\lambda\in\CP_\bq(n,k)$,
\begin{align}\notag
\mu^\bq_{z,\myp k}(\lambda\mypp|\mypp
N(\lambda)=n)&=\frac{\mu^\bq_{z,\myp k}(\lambda)}{\mu^\bq_{z,\myp
k}\left\{N(\lambda)=n\right\}}\\
\label{eq:mu/nu}
&=\frac{\re^{-zn}/F_\bq(z,k)}{\sum_{\lambda\in\CP_\bq(n,\myp
k)}\re^{-zn}/F_\bq(z,k)}=\frac{1}{\#\CP_\bq(n,\myp
k)}=\nu^\bq_{n,\myp k}(\lambda).
\end{align}

The following fact will be instrumental below.
\begin{lemma}\label{lm:D_j}
Under the measure $\mu^\bq_{z,\myp k}$, the differences
$(D_j(\lambda))_{j=1}^k$ are independent random variables such that
the marginal distribution of $D_j(\lambda)-q_{k-j}\in\NN_0$ is
geometric with parameter $1-\re^{-zj}$
\textup{(}$j=1,\dots,k$\textup{)}\textup{;} that is, for any
$(d_1,\dots,d_k)\in \CalD_\bq(k)$,
$$
\mu^\bq_{z,\myp k}\{\lambda\in\CP_\bq(\DOT,k)\colon
D_j(\lambda)=d_j,\ j=1,\dots,k\}=\prod_{j=1}^k
(1-\re^{-zj})\mypp\re^{-zj(d_j-q_{k-j})}.
$$
In particular, the expected values are given by
\begin{equation}\label{eq:E(D)}
\BE^\bq_{z,\myp k} [D_j(\lambda)]=
q_{k-j}+\frac{\re^{-zj}}{1-\re^{-zj}}\qquad (j=1,\dots,k).
\end{equation}
\end{lemma}
\proof The claim easily follows from the representation of
$N(\lambda)$ through $(D_j(\lambda))$ (see~\eqref{eq:NK}) and the
product structure of the Laplace generating
function~\eqref{eq:gfpnm}.
\endproof

Similarly, we can assign the weight $\re^{-z N(\lambda)}$ to each
partition $\lambda\in\CP_\bq=\bigcup_{k=0}^\infty\! \CP_\bq(\DOT,k)$
normalized by
\begin{align}\label{eq:Fk}
F_\bq(z):=\sum_{\lambda\in\CP_\bq}\re^{-zN(\lambda)}
&=1+\sum_{k=1}^\infty
F_\bq(z,k)\\
\label{eq:F}
&=1+\sum_{k=1}^\infty\frac{\re^{-zs_k}}{(1-\re^{-z})
\cdots(1-\re^{-zk})}\,.
\end{align}
Note that the series (\ref{eq:F}) converges for all $z>0$, since it
is bounded by the convergent series $\sum_k
\re^{-zk}(1-\re^{-z})^{-1}\cdots(1-\re^{-zk})^{-1}=
\prod_j\myn(1-\re^{-zj})^{-1}$. This way, we get the probability
measure
\begin{equation}\label{eq:mu_q}
\mu^\bq_z(\lambda):=\frac{\re^{-zN(\lambda)}}{F_\bq(z)},\qquad
\lambda\in\CP_\bq\myp.
\end{equation}
Similarly to \eqref{eq:mu/nu}, it is easy to check that the measure
$\mu^\bq_z$ conditioned on $\{N(\lambda)=n\}$ coincides with the
uniform measure $\nu^\bq_n$ on $\CP_\bq(n)$,
\[
\mu^\bq_z(\lambda\mypp|\mypp
N(\lambda)=n)=\frac{1}{p_\bq(n)}=\nu^\bq_n(\lambda),\qquad
\lambda\in \CP_\bq(n).
\]
Furthermore, the definition (\ref{eq:mu_q}) implies
\begin{equation}\label{eq:muxk}
\mu^\bq_z\{\lambda\in\CP_\bq\colon K(\lambda)=k\}
=\frac{F_\bq(z,k)}{F_\bq(z)},\qquad k\in\NN.
\end{equation}

\begin{figure}[ht]
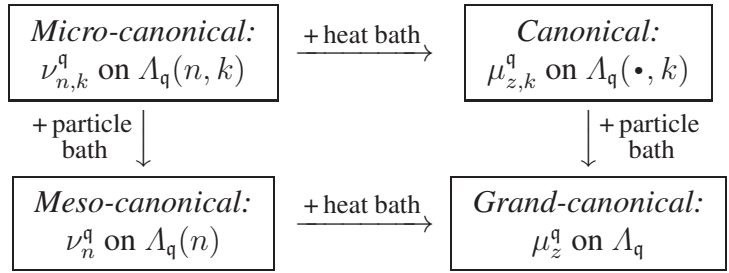

$$
\begin{CD}
\fbox{$\begin{array}{c}
\textit{Micro-canonical:}\\
\nu^\bq_{n,k}\text{ on }\CP_\bq(n,k)\end{array}$}
@>\mbox{\footnotesize +\,heat bath\,}>> \fbox{$\begin{array}{c}
\textit{Canonical:}\\
\mu^\bq_{z,k}\text{ on }\CP_\bq(\DOT,k)\end{array}$}\\
@V\mbox{\footnotesize$\begin{array}{c}\text{+\,particle}\\[-.2pc]\,\text{bath}\end{array}$}\hspace{-.4pc}VV
@VV\mbox{\hspace{-.4pc}\footnotesize$\begin{array}{c}\text{+\,particle}\\[-.2pc]\,\text{bath}\end{array}$}V\\
\fbox{$\begin{array}{c}
\textit{Meso-canonical:}\\
\nu^\bq_{n}\text{ on }\CP_\bq(n)\end{array}$} @>\mbox{\footnotesize
+\,heat bath\,}>> \fbox{$\begin{array}{c}
\textit{Grand-canonical:}\\
\mu^\bq_{z}\text{ on }\CP_\bq\end{array}$}
\end{CD}
$$
\caption{Schematic diagram illustrating the relation between
different MDP-ensembles. The integer parameters $n$ and $k$ are
interpreted as the total energy of the (quantum) system and the
number of particles, respectively. The arrows ``heat bath'' and
``particle bath'' indicate that fixation of energy or the number of
particles is lifted.}\label{fig:ensembles}\end{figure} We finish
this subsection by a comment linking the above MDP spaces and
probability measures on them with the general nomenclature of
ensembles in statistical mechanics (see, e.g., the monographs by
Huang \cite{Huang} or Greiner et al.~\cite{Greiner}). Under the
quantum interpretation of integer partitions
$\lambda=(\lambda_i)\in\CP$ briefly mentioned in Section
\ref{sec:1.1}, the $\MDP(\bq)$ restriction determines the exclusion
rules for permissible energy levels $(\lambda_i)$. In general, the
weight $N(\lambda)$ (total energy) and length $K(\lambda)$ (number
of particles) are random. Fixing one or both of these parameters
leads to different measures on the corresponding spaces, and
therefore determines different ensembles. In particular, a
completely isolated system, with fixed $N(\lambda)=n$ and
$K(\lambda)=k$ and under uniform measure $\nu^\bq_{n,k}$ on the
corresponding space $\CP_\bq(n,k)$, has the meaning of
\emph{micro-canonical} MDP ensemble. When, say, the fixation
$N(\lambda)=n$ is lifted (which may be thought of as connecting the
system to a \emph{heat bath}, whereby thermal equilibrium is settled
through exchange of energy with the bath), we get an enlarged space
$\CP(\DOT,k)$ with the measure $\mu^{\bq}_{z,k}$, which is
interpreted as the \emph{canonical} ensemble, with a fixed number of
particles~$k$. Furthermore, removing the latter constraint (which,
similarly, is achieved by putting the canonical ensemble into a
\emph{particle bath} allowing free exchange of particles) leads to
the space $\CP_\bq$ with the measure $\mu^\bq_{z}$, which is
referred to as the \emph{grand canonical} ensemble (see the
schematic diagram in Fig.\:\ref{fig:ensembles}).

Note however that the space $\CP_\bq(n)$ (i.e.,\ with a fixed energy
$N(\lambda)=n$ and endowed with uniform measure $\nu^\bq_n$), which
is most natural from the combinatorial point of view, is missing in
this picture; indeed, it may not be physically meaningful to talk
about systems with fixed energy and free number of particles. But
logically, it is perfectly possible to interchange the order of
relaxations described above and first lift the condition
$K(\lambda)=k$ by connecting the micro-canonical system to a
particle bath; we take the liberty to call the resulting ensemble
\emph{meso-canonical},\footnote{This is just a placeholder in lieu
of an established physical term.} indicating an intermediately
coarse partitioning of the phase space (cf.~\cite{3levels}).
Finally, removing the remaining constraint $N(\lambda)=n$ (by
connecting the system further to a heat bath) we again obtain the
grand canonical ensemble.

\subsection{Asymptotic equivalence of ensembles}
For $q\ge0$, define the function
\begin{equation}\label{eq:deftheta}
\vartheta_q(t):=\sqrt{\tfrac12\mypp q\myp t^2+\Li_2(1-\re^{-t})},
\qquad \quad t>0,
\end{equation}
where $\Li_2(\cdot)$ is the dilogarithm (see~\eqref{eq:Dilog}).
Recall that $T_q>0$ is the unique solution of the equation
(cf.~\eqref{eq:Tq-intro})
\begin{equation}\label{eq:deftd}
\re^{-q T_q}=1-\re^{-T_q}.
\end{equation}
Note that the value $\vartheta_q(T_q)$ coincides with the notation
$\vartheta_q$ introduced in~\eqref{eq:theta-intro-alt}.

The following curious identity will be explained in
Section~\ref{sec:4.4}.
\begin{lemma}\label{lm:duality_identity}
For all $q>0$, we have
\begin{equation}\label{eq:Tq-identity}
T_{q^{-1}}=q T_q.
\end{equation}
\end{lemma}

\proof Rewriting the equation \eqref{eq:deftd} in the form
$\re^{-q^{-1}(qT_q)}=1-\re^{-q T_q}$, we see that $\tau=qT_q$
satisfies \eqref{eq:deftd} with $q$ replaced by $q^{-1}$. By
uniqueness, this implies~\eqref{eq:Tq-identity}.
\endproof

The next proposition establishes an asymptotic link between the
measures $\mu^\bq_z$ and $\nu^\bq_n$.
\begin{proposition}\label{pr:munu}
Suppose that the sequence\/ $\bq$ satisfies the
condition~\eqref{eq:assq*}. Let $\{A_z\}_{z>0}$ be a family of
subsets of the space $\CP_\bq$ such that, for some positive constant
$\kappa$,
\begin{equation}\label{eq:limsupmuz}
\limsup_{z\downarrow0} z^{\kappa} \log \mu^\bq_z(A_z)<0.
\end{equation}
Then there exists a sequence $(z_n)$ such that
\begin{equation}\label{eq:xn}
\lim_{n\to\infty}z_n \sqrt{n}=\vartheta_q \equiv \vartheta_q(T_{q})
\end{equation}
and \begin{equation}\label{eq:log(nu)<0a}
\limsup_{n\to\infty}\,n^{-\kappa/2}\log \nu^\bq_n(A_{z_n})<0.
\end{equation}
\end{proposition}

There is a similar connection between the measures $\mu^\bq_{z,\myp
k}$ and $\nu^\bq_{n,\myp k}$, provided that $z\downarrow0$,
$k\to\infty$ and $n\to\infty$ in a coordinated manner.

\begin{proposition}\label{pr:mumnum}
Let the sequence $\bq=(q_i)$ satisfy Assumption
\textup{\ref{as:main}}. Let a family of sets $A_{z,\myp
k}\subset\CP_\bq(\DOT,k)$ \textup{($z>0$, $k\in\NN$)} be such that,
for some constant $\kappa>0$,
\begin{equation}\label{eq:muxmcd}
\limsup_{z\downarrow0}\,z^{\kappa} \log \mu^\bq_{z,\myp
k(z)}(A_{z,\myp k(z)})<0,
\end{equation}
for any $k=k(z)$ such that $z k(z)\to T\in(0,\infty]$ as
${z\downarrow0}$.

\smallskip
\textup{(a)} If\/ $T<\infty$ then for any sequence $(k_n)$
satisfying
\begin{equation}\label{eq:kn1}
\lim_{n\to\infty} \frac{k_n}{\sqrt{n}}=\frac{T}{\vartheta_q(T)},
\end{equation}
there exists a sequence $(z_{n})$ such that
\begin{equation}\label{eq:xnT}
\lim_{n\to\infty} z_{n}\sqrt{n}=\vartheta_q(T)
\end{equation}
and
\begin{equation}\label{eq:nunmcd}
\limsup_{n\to\infty}\,n^{-\kappa/2}\log \nu^\bq_{n,\myp
k_n}(A_{z_{n},\,k_n})<0.
\end{equation}

\smallskip
\textup{(b)} Let\/ $T=\infty$ and $q=0$, and assume in addition that
$z^{2/(\beta+1)}k(z)\to 0$ as $z\downarrow 0$. Then for any sequence
$(k_n)$ satisfying
\begin{equation}\label{eq:knass}
\lim_{n\to\infty}\frac{k_n }{k(\pi/\sqrt{6n}\mypp)}=1,
\end{equation}
there exists a sequence $(z_{n})$ such that the asymptotic relations
\eqref{eq:xnT} and \eqref{eq:nunmcd} hold true, with the right-hand
side of \eqref{eq:xnT} reducing to $\vartheta_0(\infty)\equiv
\vartheta_0=\pi/\sqrt{6}$ \,\textup{(}see~\eqref{eq:q=0}\textup{)}.
\end{proposition}

These two propositions are instrumental for our method; their proof,
being rather technical, is postponed until
Appendix~\ref{sec:prop12pf}.

\section{Number of parts in a typical MDP}\label{sec:m}

In this section, our ultimate goal is to show that if Assumption
\ref{as:main} holds then, under the measures $\nu^\bq_n$ on the
$\MDP$-space $\varLambda_\bq(n)$, the typical length $K(\lambda)$
(i.e., the number of parts) of a partition $\lambda\in\CP_\bq(n)$ of
large weight $N(\lambda)=n$ is concentrated around $c\myp\sqrt{n}$
(with a suitable constant $c>0$) if $q>0$, or grows slightly faster
than $\sqrt{n}$ if $q=0$. To this end, we will first study the
distribution of $K(\lambda)$ under the measure $\mu_z^\bq$ in the
space $\CP_\bq$.

\subsection{Preparatory lemmas}

For $z>0$, denote
\begin{equation}\label{eq:gamma}
\eta_{k}(z):= \frac{\re^{-z\myp Q_{k}}}{1-\re^{-zk}},\qquad k\in\NN,
\end{equation}
where $Q_k$ is given by \eqref{eq:assq*}. For every $z>0$, the
sequence $(\eta_k(z))_{k\ge1}$ is decreasing, and in particular
$$
\eta_k(z)\le \eta_1(z)=\frac{\re^{-z q_{0}}}{1-\re^{-z}},\qquad
k\in\NN.
$$
Furthermore,
$$
0\le \lim_{k\to\infty} \eta_k(z)\le \re^{-z q_{0}}<1.
$$
Thus, the set $\{k\colon \eta_k(z)\ge 1\}$ is always finite
(possibly empty). Define
\begin{equation}\label{eq:def-k0}
k_*\equiv k_*(z):=\begin{cases}\max\{k\in\NN\colon\eta_{k}(z)\ge
1\}&\text{if}\ \ \eta_1(z)\ge1,\\
1&\text{if}\ \ \eta_1(z)<1. \end{cases}
\end{equation}

\begin{remark}
Note that $\lim_{z\downarrow0}\eta_k(z)=+\infty$ for any fixed
$k\in\NN$, and so $k_*(z)>1$ for all $z>0$ small enough.
\end{remark}

First, let us record a few auxiliary statements that do not depend
on Assumption~\ref{as:main}.
\begin{lemma}\label{lm:mukseq}\mbox{}
\begin{itemize}
\item[\rm (a)] For every $z>0$, we have
$$
\max_{k\in\NN}\mu^\bq_z\{K(\lambda)=k\}=\mu^\bq_z\{K(\lambda)=k_*\},
$$
where $k_*=k_*(z)$ is defined in \eqref{eq:def-k0}. Moreover\/,
\begin{equation}\label{eq:mukseq}
\mu^\bq_z\{K(\lambda)=k_*\}>\mu^\bq_z\{K(\lambda)=k_*+1\}>\mu^\bq_z\{K(\lambda)=k_*+2\}>\cdots,
\end{equation}
and, for $k_*\ge 2$,
\begin{equation}\label{eq:mukseq1}
\mu^\bq_z\{K(\lambda)=k_*\}\ge
\mu^\bq_z\{K(\lambda)=k_*-1\}>\dots>\mu^\bq_z\{K(\lambda)=1\}.
\end{equation}

\item[\rm (b)] The function $z\mapsto k_*(z)$ is non-increasing and
has no jumps larger than~$1$. Moreover, $k_*(z)\to +\infty$ as
$z\downarrow0$.
\end{itemize}
\end{lemma}

\proof (a) Using \eqref{eq:gfpnm}, \eqref{eq:s} and \eqref{eq:muxk},
we can rewrite \eqref{eq:gamma} (for $k\ge2$) as
\begin{equation}\label{eq:gamma2}
\eta_{k}(z) =\frac{F_\bq(z, k)}{F_\bq(z,k-1)}
=\frac{\mu^\bq_z\{K(\lambda)=k\}}{\mu^\bq_z\{K(\lambda)=k-1\}}.
\end{equation}
The definition of $k_*=k_*(z)$ (see~\eqref{eq:def-k0})
implies that $\eta_{k}(z)<1$ for $k>k_*$, and \eqref{eq:mukseq}
follows. Similarly, assuming that $k_*\!\ge2$, we have
$\eta_{k_*}\mynn(z)\ge 1$ and $\eta_{k}(z)>1$ for $k<k_*$\myp, which
is the same as~\eqref{eq:mukseq1}.

\smallskip
(b) For $k\in\NN$, let $z=\zeta_k$ be the (unique) solution of the
equation
\begin{equation}\label{eq:gamma=1}
\eta_{k}(z)=1.
\end{equation}
From the formulas \eqref{eq:gamma} and \eqref{eq:gamma=1}, it is
clear that the sequence $(\zeta_k)_{k\ge1}$ is decreasing and,
moreover, $\zeta_k\downarrow0$ as $k\to\infty$. If $z=\zeta_k$
\,($k\ge2$) then
$\mu^\bq_{\zeta_k}\{K(\lambda)=k\}=\mu^\bq_{\zeta_k}\{K(\lambda)=k-1\}$
are the two maxima of the sequence
\strut{}$(\mu^\bq_{z}\{K(\lambda)=j\})_{j\ge 1}$, whereas for
$z\in(\zeta_{k+1},\zeta_{k})$ the unique maximum of this sequence is
attained exactly at $j=k$. Hence, $k_*(z)\equiv k$ for
$z\in(\zeta_{k+1},\zeta_k]$, that is, $z\mapsto k_*(z)$ is a
non-increasing (left-continuous) step function with unit downward
jumps at points $\zeta_k$ ($k\ge2$). Since
$\lim_{k\to\infty}\zeta_k=0$, it also follows that
$\lim_{z\downarrow0}k_*(z)=+\infty$.
\endproof

\begin{remark}
Willing to use a ``one-sided'' version of the notation
$f(z)=O(g(z))$ ($z\downarrow 0$), in what follows we write $f(z)\le
O(g(z))$ ($z\downarrow0$) \,if \,$\limsup_{z\downarrow0} f(z)/
g(z)<+\infty$.
\end{remark}

\begin{lemma}\label{lm:new1}
Uniformly in $k\in\NN$, as $z\downarrow0$,
\begin{align}
\label{eq:mld:t2} \log \mu^\bq_z\{K(\lambda)=k\} &\le
z^{-1}\bigl(\Li_2(\re^{-zk_*}\myn)-\Li_2(\re^{-zk})\bigr)
+z\myp(s_{k_*}\!-s_k)+O\bigl(\log\tfrac1z\bigr)\\[.3pc]
\label{eq:mld:t3} &\le (k_*-k)\log\myn(1-\re^{-zk_{*}}\myn)
+z\myp(s_{k_*}\!-s_k)+O\bigl(\log\tfrac1z\bigr).
\end{align}
\end{lemma}

\proof Recalling \eqref{eq:muxk}, for each $k\in\NN$ we can write
(see~\eqref{eq:Fk})
\begin{equation}\label{eq:new2}
\log \mu^\bq_z\{K(\lambda)=k\} =\log \frac{F_\bq(z,k)}{F_\bq(z)}\le
\log F_\bq(z,k)-\log F_\bq(z,k_*),
\end{equation}
where, according to \eqref{eq:gfpnm},
\begin{equation}\label{eq:logFk}
\log F_\bq(z,k)=-z s_k-\sum_{j=1}^k\log\myn(1-\re^{-z j})\qquad
(k\in\NN).
\end{equation}
By the well-known Euler--Maclaurin sum formula \cite[23.1.36,
p.\,806]{AS} applied to the function $x\mapsto
\log\myn(1-\re^{-zx})$, we get, uniformly in $k\in\NN$ as
$z\downarrow0$,
\begin{align}
\notag \sum_{j=1}^k\log\myn(1-\re^{-zj})&=\int_1^k
\log\myn(1-\re^{-zx})\,\rd{x} +O(1)\log\myn(1-\re^{-z})
+O(1)\int_1^k
\frac{z\mypp\re^{-zx}}{1-\re^{-zx}}\,\rd{x}\\
\label{eq:sum-log-Psi}
&=z^{-1}\bigl(\Li_2(\re^{-zk})-\Li_2(\re^{-z})\bigr)+O\bigl(\log
\textstyle{\frac1z}\bigr),
\end{align}
where $\Li_2(\cdot)$ is the dilogarithm function
(see~\eqref{eq:Dilog}). Thus, substituting \eqref{eq:sum-log-Psi}
into \eqref{eq:logFk} and returning to \eqref{eq:new2}, we
obtain~\eqref{eq:mld:t2}.

Furthermore, since the derivative
$(\Li_2(\re^{-t}))'=\log\mynn(1-\re^{-t}) $ is increasing in
$t\in(0,\infty)$, the function $t\mapsto\Li_2(\re^{-t})$ is convex,
hence
\begin{equation*}
\Li_2(\re^{-zk_*}\myn)-\Li_2(\re^{-zk})\le z\mypp
(k_*-k)\log\myn(1-\re^{-zk_*}\myn),\qquad k\in\NN.
\end{equation*}
Combining this bound with \eqref{eq:mld:t2}
yields~\eqref{eq:mld:t3}.
\endproof

\begin{lemma}\label{lm:mukseqasmain} Suppose that Assumption \textup{\ref{as:main}} is in force, that is,
the sequence $\bq=(q_i)$ satisfies \eqref{eq:assq*} with $q\ge0$ and
$0\le\beta<1$.
\begin{itemize}
\item[\rm (a)] If $q>0$ then
\begin{equation} \label{eq:k*}
k_*(z)= z^{-1}\myp T_q+O(z^{-\beta})\qquad (z\downarrow 0),
\end{equation}
where $T_q$ is defined in~\eqref{eq:deftd}.

\item[\rm (b)]  If $q=0$ then
\begin{equation}\label{eq:k0sublin}
1-\beta\le \liminf_{z\downarrow0}\frac{z\myp k_*(z)}{\log\frac1z}\le
\limsup_{z\downarrow0}\frac{z\myp k_*(z)}{\log\frac1z}\le 1.
\end{equation}
\end{itemize}
In particular, for all $q\ge0$,
\begin{equation} \label{eq:k0}
\lim_{z\downarrow 0} z k_*(z)= T_q.
\end{equation}

\end{lemma}

\proof (a) Like in the proof of Lemma~\ref{lm:mukseq}(b), denote by
$\zeta_k$ ($k\in\NN$) the solution of the equation
\eqref{eq:gamma=1}. Using the definition \eqref{eq:gamma}, equation
\eqref{eq:gamma=1} is expressed at $z=\zeta_k$ as
\begin{equation}\label{eq:eqn-k*}
k^{-1}Q_{k}=-(k\zeta_k)^{-1} \log\myn(1-\re^{-k\zeta_k}\myn).
\end{equation}
Comparing this with equation \eqref{eq:Tq-intro},
%\eqref{eq:deftd},
observe that $k\myp\zeta_k=T_{\tilde{q}_k}$, where
$\tilde{q}_k:=k^{-1}Q_{k}\to q>0$ as $k\to\infty$, due to the
limit~\eqref{eq:q*}, and therefore
$\lim_{k\to\infty}T_{\tilde{q}_k}=T_q$, thanks to continuity of the
mapping $q\mapsto T_q$.
% mentioned after the definition~\eqref{eq:deftd}.

To see why this implies \eqref{eq:k*}, recall from the proof of
Lemma~\ref{lm:mukseq}(b) that $k_*(z)\equiv k$ for
$z\in(\zeta_{k+1},\zeta_k]$ \,($k\in\NN$) and the limit
$z\downarrow0$ is equivalent to $k\to\infty$. Hence,
$$
k_*z=k\myp\zeta_k-k\myp(\zeta_k-z)\to T_q\qquad (z\downarrow0),
$$
because $k\myp\zeta_k\to T_q$ and
$$
0\le k\myp(\zeta_k-z)\le  k\myp\zeta_k-k\myp \zeta_{k+1}\to0\qquad
(k\to\infty).
$$

Furthermore, by a standard perturbation analysis it is easy to
estimate the corresponding remainder term in the
limit~\eqref{eq:k*}. Indeed, setting $\delta_k:=k\myp\zeta_k-T_q\to
0$ and using the asymptotic relation \eqref{eq:assq*}, we can
rewrite \eqref{eq:eqn-k*} in the form
$$
(T_q+\delta_k)\left(q+O(k^{\beta-1})\right)
=-\log\myn\bigl(1-\re^{-T_q}\myn\bigr)-\frac{\re^{-T_q}}{1-\re^{-T_q}}\,\delta_k+O(\delta_k^2),
$$
which yields, in view of the identity \eqref{eq:deftd}, that
$\delta_k=O(k^{\beta-1})$.

In turn, for $\zeta_{k+1}<z\le \zeta_k$ we get
\begin{align*}
k_* z-T_q&=(k\myp
\zeta_k-T_q)-k\myp(\zeta_k-z)\\
&=\delta_k+O(1)\mypp(\delta_k+\delta_{k+1})\\
&=O(k^{\beta-1})=O(z^{1-\beta})\qquad (z\downarrow0),
\end{align*}
and the estimate~\eqref{eq:k*} follows.

\smallskip
(b) Fix $\varepsilon\in(0,1-\beta)$. For $z>0$ small enough,
$\eta_{k_*}(z)=\re^{-z\myp Q_{k_*}}(1-\re^{-z k_*})^{-1}\ge 1$ by
the definition~\eqref{eq:def-k0}, so
\[
\re^{-zk_*}\ge 1-\re^{-z\myp Q_{k_*}} \ge 1-\re^{-z}\ge
z^{1+\varepsilon},
\]
because $Q_{k_*}\ge q_0\ge 1$. Thus,
\begin{equation}\label{eq:zk0le}
zk_*(z)\le (1+\varepsilon)\log\tfrac{1}{z},
\end{equation}
which implies the last inequality in \eqref{eq:k0sublin}, since
$\varepsilon>0$ can be taken arbitrarily close to~$0$.

On the other hand, from \eqref{eq:deftd} we also have
$\eta_{k_*+1}(z)<1$, that is,
\begin{equation}\label{eq:liminf_k*}
z\myp k_*(z)>\log \tfrac{1}{z}-z-\log Q_{k_*+1}.
\end{equation}
Furthermore, using the asymptotic bound \eqref{eq:assq*} for $k=k_*$
(with $q=0$) and the estimate \eqref{eq:zk0le}, we obtain
$$
\log
Q_{k_*+1}=O(1)+\beta\log\tfrac{1}{z}+\beta\log\log\tfrac{1}{z}\qquad
(z\downarrow0).
$$
Substituting this into \eqref{eq:liminf_k*}, it is easy to see that
$$
\liminf_{z\downarrow0}\frac{z\myp k_*(z)}{\log \frac{1}{z}}\ge
1-\lim_{z\downarrow0}\frac{z}{\log \frac{1}{z}}-
\lim_{z\downarrow0}\frac{\log Q_{k_*+1}}{\log \frac{1}{z}}=1-\beta,
$$
and the first inequality in~\eqref{eq:k0sublin} is proved.
\endproof

\begin{remark}\label{rem:k0q=0}
In the case $q=0$, the asymptotic bounds in \eqref{eq:k0sublin} are
optimal in the following sense: under Assumption~\ref{as:main+}
(i.e., when $Q_k\sim \tilde{q}\myp k^{\beta}$ as $k\to\infty$), one
can show that $\lim_{z\downarrow 0}z\mypp(\log\frac1z)^{-1}\myp
k_*(z)=1-\beta>0$.
\end{remark}

\subsection{Asymptotics of $K(\lambda)$ in the space $\varLambda_\bq$: case $q>0$}
We can now give exponential estimates on the asymptotic behaviour of
the random variable $K=K(\lambda)$ (see~\eqref{eq:NK}) under the
measure~$\mu^\bq_{z}$. We start with the case $q>0$.
\begin{theorem}\label{th:mld}
Let the sequence $\bq=(q_i)$ satisfy Assumption
\textup{\ref{as:main}} with $q>0$ and $0\le\beta<1$. Then, for every
$\gamma\in\bigl(0,\frac12(1-\beta)\bigr)$ and any constant $c>0$, we
have
\begin{equation}\label{eq:muxdevk}
\limsup_{z\downarrow0} z^{1-2\gamma}\log
\mu^\bq_z\bigl\{\lambda\in\CP_\bq\colon |K(\lambda)-k_*|>c\myp
z^{\gamma-1}\bigr\}\le -\tfrac12\myp q\myp c^2<0,
\end{equation}
where $k_*=k_*(z)$ is defined in \eqref{eq:def-k0}.
\end{theorem}

\proof From \eqref{eq:muxk} we have
\begin{equation}\label{eq:I_z}
\mu^\bq_z\bigl\{|K(\lambda)-k_*|>c\myp z^{\gamma-1}\bigr\}
=\frac{1}{F_\bq(z)}\sum_{k\in\mathcal{I}_z} F_\bq(z,k),
\end{equation}
where $\mathcal{I}_z:=\{k\in\NN\colon |k-k_*|>c\myp z^{\gamma-1}\}$.
Recalling \eqref{eq:gamma} and \eqref{eq:gamma2}, observe that for
$k>2k_*$
\begin{equation}\label{eq:F/F}
\frac{F_\bq(z,k)}{F_\bq(z,k-1)}=\eta_{k}(z) \le \eta_{2k_*}(z) =
\frac{\re^{-z\myp Q_{2k_*}}}{1-\re^{-2zk_*}}.
\end{equation}
By the asymptotic formulas \eqref{eq:assq*} and \eqref{eq:k0}, this
gives
\begin{align} \notag \limsup_{z\downarrow0}\log
\frac{F_\bq(z,k)}{F_\bq(z,k-1)} &\le
-\lim_{z\downarrow0} z\mypp Q_{2k_*}-\lim_{z\downarrow0}\log\myn(1-\re^{-2zk_*}\myn)\\
\notag
&= -2\myp q\myp T_q-\log\myn(1-\re^{-2\myp T_q})\\[.2pc]
\label{eq:<0} &< -2\myp q\myp T_q-\log\myn(1-\re^{-T_q})=-q\myp
T_q<0,
\end{align}
where the last equality in \eqref{eq:<0} is due to
equation~\eqref{eq:deftd}. Hence, the part of the sum \eqref{eq:I_z}
with $k>2k_*$ is asymptotically dominated by a geometric series with
ratio $\re^{-qT_q}<1$, so that
\begin{equation}\label{eq:sum>2k0}
\frac{1}{F_\bq(z)}\sum_{k>2k_*} F_\bq(z,k) \le
\frac{F_\bq(z,2k_*)}{F_\bq(z)}\cdot \frac{\re^{-q\myp
T_q}}{1-\re^{-q\myp T_q}}.
\end{equation}
Furthermore, with the help of the asymptotic relations
\eqref{eq:sksim} and \eqref{eq:k*} and in view of the
equation~\eqref{eq:deftd}, the estimate \eqref{eq:mld:t3}
specializes as follows
\begin{align}
\notag \log \frac{F_\bq(z,2k_*)}{F_\bq(z)} &\le -
z^{-1}\bigl\{T_q\log\myn(1-\re^{-T_q})+\tfrac{3}{2}\myp
q\myp T_q^{2}+O(z^{1-\beta})\bigr\}+O\bigl(\log\tfrac1z\bigr)\\
\label{eq:frac>2k0}&=-\tfrac{1}{2}\myp q\myp  T_q^{2}
z^{-1}+O(z^{-\beta}) \qquad (z\downarrow0).
\end{align}

Let us now turn to the case $k\le 2k_*$. Denote
$k_-:=\floor{k_*-cz^{\gamma-1}}$,
\,$k_+:=\ceiling{k_*+cz^{\gamma-1}}$ (here and in what follows,
$\floor{\cdot}$ and $\ceiling{\cdot}$ denote the floor and ceiling
functions, respectively). Observe from
Lemma~\ref{lm:mukseqasmain}(a) that $k_\pm\!\to\infty$ as
$z\downarrow 0$ and $k_-\!<k_*\!<k_+\!<2k_*$. Hence, the
monotonicity properties \eqref{eq:mukseq} and \eqref{eq:mukseq1}
yield
\begin{align}
\notag \frac{1}{F_\bq(z)}\sum_{k\in\mathcal{I}_z,\,k\le 2k_*}\!
F_\bq(z,k) &\le \frac{2\myp
k_*}{F_\bq(z)}\max_{k\in\mathcal{I}_z,\,k\le 2k_*}F_\bq(z,k)\\
\label{eq:sum<2k0}
&\le 2\myp k_*\,
\frac{F_\bq(z,k_-)+F_\bq(z,k_+)}{F_\bq(z)}.
\end{align}
Similarly to \eqref{eq:frac>2k0}, from \eqref{eq:mld:t3} we obtain
\begin{align}
\notag \log\frac{F_\bq(z,k_\pm)}{F_\bq(z)}&\le \mp\myp c\myp
z^{\gamma-1}\bigl\{\log\myn(1-\re^{-T_q})+ q\myp
T_q\bigr\}-\tfrac12\myp q\mypp c^2 z^{2\gamma-1}+O(z^{\gamma-\beta})+O\bigl(\log\tfrac1z\bigr)\\
\label{eq:frac<2k0} &=-\tfrac12\myp q\mypp c^2 z^{2\gamma-1}
+O(z^{\gamma-\beta})+O\bigl(\log\tfrac1z\bigr)\qquad (z\downarrow0),
\end{align}
again by making use of the equation~\eqref{eq:deftd}.

Finally, returning to the expansion \eqref{eq:I_z} and combining the
estimates \eqref{eq:sum>2k0}, \eqref{eq:frac>2k0},
\eqref{eq:sum<2k0} and \eqref{eq:frac<2k0}, with the help of the
elementary inequality
\begin{equation}\label{eq:loga+b}
\log\myn(x+y)\le \log 2+\max\myn\{\log x,\log y\},\qquad x,y>0,
\end{equation}
we obtain \eqref{eq:muxdevk}, which completes the proof.
\endproof

Theorem \ref{th:mld} combined with the asymptotic
formula~\eqref{eq:k*} implies the following law of large numbers for
the number of parts $K(\lambda)$ under the measure $\mu_z^\bq$.
\begin{corollary}\label{cor:LLN_K}
Let Assumption \textup{\ref{as:main}} hold with $q>0$. Then, for any $\varepsilon>0$,
\begin{equation*}
\lim_{z\downarrow0} \mu^\bq_z\bigl\{\lambda\in\CP_\bq\colon |z\myp
K(\lambda)-T_q|>\varepsilon\bigr\}=0.
\end{equation*}
\end{corollary}

\subsection{Asymptotics of $K(\lambda)$ in the space $\varLambda_\bq$: case $q=0$}
When Assumption \ref{as:main} holds with $q=0$, the asymptotics for
$k_*(z)$ as $z\downarrow0$ cannot be obtained, as was mentioned in
Remark~\ref{rem:k0q=0}. So there is no hope to find exponential
bounds for $K(\lambda)$ to fit into an interval of order smaller
than $z^{-1}$, as in \eqref{eq:muxdevk}. Nevertheless we can still
find an interval such that $K(\lambda)$ does not hit it with an
exponentially small $\mu^\bq_z$-probability, as $z\downarrow 0$. To
this end, we need some additional notation.

Fix $\gamma\in(0,1)$ and define the function
\begin{equation}\label{eq:defkappa}
z\mapsto k_\gamma\equiv k_\gamma(z):=\inf\{k\in\NN\colon s_k\ge
z^{-2\myp(1-\gamma)}\},\qquad z\in(0,1).
\end{equation}
Recalling that $s_k\ge k$ (see after formula \eqref{eq:s}), from the
definition \eqref{eq:defkappa} it follows that
\begin{equation}\label{eq:kappa_gamma<}
k_\gamma(z)\le \lceil z^{-2\myp(1-\gamma)}\rceil.
\end{equation}
On the other hand, it is clear that $k_\gamma(z)\to\infty$ as
$z\downarrow 0$. Actually we can tell more.

\begin{lemma}\label{lm:k4}
Let\/ Assumption~\textup{\ref{as:main}} hold with $q=0$ and some
$\beta\in[0,1)$. Then, for any
$\gamma\in(0,1)$,
\begin{gather}\label{eq:lm3.5-1}
\lim_{z\downarrow 0}z^{2\myp(1-\gamma)}
s_{k_\gamma(z)}=1,\\
\label{eq:lm3.5-2} \liminf_{z\downarrow
0}z^{2\myp(1-\gamma)/(\beta+1)}
k_\gamma(z)>0.
\end{gather}
Moreover, if\/ $0<\gamma<\frac12$ then
for any $t>0$
\begin{equation}\label{eq:lm3.5-3} \limsup_{z\downarrow
0}z^{1-2\gamma}\myp Q_{k_\gamma(z)-\floor{t/z}}\le t^{-1}.
\end{equation}
\end{lemma}

\proof The definition \eqref{eq:defkappa} implies that
$s_{k_\gamma-1}<z^{-2\myp(1-\gamma)}\le s_{k_\gamma}$. Hence,
recalling notation \eqref{eq:s} and combining the asymptotics
\eqref{eq:assq*} (with $q=0$) and the bound~\eqref{eq:lm3.5-3}, we
have
\begin{align}
\notag z^{-2\myp(1-\gamma)}\le
s_{k_\gamma}=s_{k_\gamma-1}+Q_{k_\gamma}&<z^{-2\myp(1-\gamma)}+Q_{k_\gamma}
\\
&=z^{-2\myp(1-\gamma)}+O(z^{-2\beta\myp(1-\gamma)})\sim
z^{-2\myp(1-\gamma)},
\label{eq:skappa}
\end{align}
since $\beta<1$ and $1-\gamma>0$. Now, the limit \eqref{eq:lm3.5-1}
follows from the two-sided estimate~\eqref{eq:skappa}. Similarly,
using \eqref{eq:sksim} (with $q=0$), we obtain the asymptotic bound
$$
z^{-2\myp(1-\gamma)}\le s_{k_\gamma}=O(k_\gamma^{\beta+1})\qquad
(z\downarrow0),
$$ which implies~\eqref{eq:lm3.5-2}. Finally, since the sequence
$(Q_k)$ is non-decreasing (see~\eqref{eq:assq*}), for $t>0$ we can
write
\[
s_{k_\gamma}\ge \sum_{k=k_\gamma-\floor{t/z}}^{k_\gamma} \! Q_k \ge
\floor{t/z}\cdot Q_{k_\gamma-\floor{t/z}},
\]
and the claim \eqref{eq:lm3.5-3} readily follows in view
of~\eqref{eq:lm3.5-1}.
\endproof

The next result is a counterpart of Theorem~\ref{th:mld} for the
case $q=0$.

\begin{theorem}\label{th:mldsublin} Let Assumption~\textup{\ref{as:main}} hold with $q=0$. Then, for any
$\gamma\in\bigl(0,\frac12(1-\beta)\bigr)$,
\begin{align}\label{eq:muxdevkq=0}
\limsup_{z\downarrow0}&\, z^{1-2\gamma}\log
\mu^\bq_z\bigl\{\lambda\in\CP_\bq\colon K(\lambda)<
z^{-1}\log\log\tfrac1z\bigr\}=-\infty,\\
\label{eq:muxdevkq=0_1} \limsup_{z\downarrow0}&\, z^{1-2\gamma}\log
\mu^\bq_z\{\lambda\in\CP_\bq\colon K(\lambda)> k_\gamma(z)\}\le -1.
\end{align}
\end{theorem}

\proof Put $k^\dag\equiv k^\dag(z):=\floor{z^{-1}\log\log\tfrac1z}$.
In view of the lower bound in \eqref{eq:k0sublin}, it is clear that
$k^\dag(z)/k_*(z)\to0$ as $z\downarrow0$, and hence $k^\dag(z)<
k_*(z)$ for all $z>0$ small enough. Then, using \eqref{eq:muxk} and
\eqref{eq:mukseq1}, we can write
\begin{equation}\label{eq:S3}
\mu^\bq_z\{ K(\lambda)< k^\dag\} =
\sum_{k<k^\dag}\frac{F_\bq(z,k)}{F_\bq(z)} \le
k^\dag\myp\frac{F_\bq(z,k^\dag)}{F_\bq(z)}.
\end{equation}
Furthermore, for any $\varepsilon\in(0,1-\beta)$ and all $z>0$ small
enough, according to \eqref{eq:k0sublin} we have
$$
(1-\beta-\varepsilon)\mypp z^{-1}\log\tfrac{1}{z}\le k_*(z)\le
(1+\varepsilon)\mypp z^{-1}\log\tfrac{1}{z},
$$
which also gives
$s_{k_*}=O\bigl(z^{-\beta-1}(\log\frac{1}{z})^{\beta+1}\bigr)$ by
\eqref{eq:sksim}. Then from \eqref{eq:mld:t2} we get
\begin{align}\notag
\log\frac{F_\bq(z,k^\dag)}{F_\bq(z)}&\le
z^{-1}\bigl\{\Li_2(z^{1-\beta-\varepsilon})-\Li_2(\re^{-zk^\dag})\bigr\}
+O\bigl(z^{-\beta}(\log\tfrac{1}{z})^{\beta+1}\bigr)\\[-.2pc] \notag
&=-z^{-1}\Li_2\bigl((\log\tfrac{1}{z})^{-1}\bigr)+O(z^{-\beta-\varepsilon})\\[.2pc]
\label{eq:Fk3} &\sim -z^{-1}
\bigl(\log\tfrac{1}{z}\bigr)^{-1}\qquad(z\downarrow0),
\end{align}
and \eqref{eq:muxdevkq=0} follows by combining \eqref{eq:S3}
and~\eqref{eq:Fk3}.

Next, to estimate the probability
\begin{equation}\label{eq:S4}
\mu^\bq_z\{ K(\lambda)> k_\gamma\} =
\frac{1}{{F_\bq(z)}}\sum_{k>k_\gamma}F_\bq(z,k),
\end{equation}
observe (cf.~\eqref{eq:F/F}) that, for $k>k_\gamma$ and all $z>0$
small enough, we have
\[
\frac{F_\bq(z,k)}{F_\bq(z,k-1)}=\eta_k(z) =\frac{\re^{-z\myp
Q_k}}{1-\re^{-zk}} < \frac{\re^{-z}}{1-\re^{-zk_\gamma}} \le
1-\tfrac12\myp z.
\]
Indeed, if $2\gamma<1-\beta$ then the asymptotic bound
\eqref{eq:lm3.5-2} implies
$\lim_{z\downarrow0}z^{-1}\myp\re^{-zk_\gamma}=0$, and therefore
$$
\frac{1}{z}\left(\frac{\re^{-z}}{1-\re^{-zk_\gamma}}-1\right)=\frac{\re^{-z}-1}{z\mypp(1-
\re^{-zk_\gamma})}+\frac{\re^{-zk_\gamma}}{z\mypp(1-
\re^{-zk_\gamma})}\to-1\qquad (z\downarrow0).
$$
Thus, we can estimate the right-hand side of \eqref{eq:S4} by the
sum of a geometric progression with ratio $1-\frac12 z<1$, that is,
\begin{equation}\label{eq:S3'}
\mu^\bq_z\{ K(\lambda)> k_\gamma\}\le 2\myp
z^{-1}\mypp\frac{F_\bq(z,k_\gamma)}{F_\bq(z)}.
\end{equation}

Next, using again the estimate \eqref{eq:mld:t2} and also the
asymptotics \eqref{eq:lm3.5-1}, we obtain (cf.~\eqref{eq:Fk3})
\begin{align}
\notag \log\frac{F_\bq(z,k_\gamma)}{F_\bq(z)}& \le
z^{-1}\Li_2(z^{1-\beta-\varepsilon})-zs_{k_\gamma}
+O\bigl(z^{-\beta}(\log\tfrac{1}{z})^{\beta+1}\bigr)\\[-.2pc] \notag
&=-z^{-1+2\gamma}\bigl(1+o(1)\bigr)+O(z^{-\beta-\varepsilon})\\[.2pc]
\label{eq:Fk4} &\sim -z^{-1+2\gamma},
\end{align}
where the asymptotic equivalence in \eqref{eq:Fk4} holds provided
that $0<\varepsilon<1-\beta-2\gamma$. Now, the desired result
\eqref{eq:muxdevkq=0_1} follows by combining \eqref{eq:S3'}
and~\eqref{eq:Fk4}.
\endproof

In the case $q=0$, under the refined Assumption \ref{as:main+} with
$\tilde{q}>0$ (see~\eqref{eq:assq*+}) one can prove the following
analogue of the exponential bound \eqref{eq:muxdevk}: \emph{for any
$c>0$ and $\gamma\in\bigl(0,\frac12\beta\bigr)$,}
\begin{equation}\label{eq:muxdevkq=0refined}
\limsup_{z\downarrow0}\,z^{\beta-2\gamma}\log
\mu^\bq_z\bigl\{\lambda\in\CP_\bq:| K(\lambda)-k_*|>c
z^{\gamma-1}\bigr\} \le -2^{\beta-2}\tilde{q}\mypp \beta\myp c^2<0.
\end{equation}
Here $k_*=k_*(z)$ is again defined by \eqref{eq:def-k0} but now has
the refined asymptotics (cf.~\eqref{eq:k0sublin})
\begin{equation}\label{eq:k*|q=0}
k_*(z)=z^{-1}\left((1-\beta)\log \tfrac{1}{z}-\beta\log \log
\tfrac{1}{z} - \beta\log\myn(1-\beta)- \log \tilde{q} +o(1)\right).
\end{equation}
The exponential bound \eqref{eq:muxdevkq=0refined} together with the
asymptotic formula \eqref{eq:k*|q=0} immediately imply the law of
large numbers for the number of parts (cf.\ Corollary
\ref{cor:LLN_K}): \emph{for any $\varepsilon>0$, }
\[
\lim_{z\downarrow 0}
\mu^\bq_z\bigl\{\lambda\in\CP_\bq\colon\bigl|z\bigl(\log\tfrac{1}{z}\bigr)^{-1}
K(\lambda)-(1-\beta)\bigr|>\varepsilon \bigr\}=0.
\]

Formally, these results do not cover the utterly degenerate case
$\tilde{q}=0$, $\tilde{\beta}=0$ in the asymptotic formula
\eqref{eq:assq*+} of Assumption~\ref{as:main+}; however, as
explained in Remark~\ref{rem:1.4}, it is equivalent to the classical
case of unrestricted partitions, where the asymptotic behaviour of
$K(\lambda)$ (under the measure $\mu_z$ on $\CP$) is described by
the limit theorem~\cite{Fristedt}
\begin{equation}\label{eq:Kexpexp}
\lim_{z\downarrow 0} \mu_z\bigl\{\lambda\in\CP\colon z
K(\lambda)-\log\tfrac1z \le t \bigr\}=\exp(-\re^{-t}),\quad t\in\RR.
\end{equation}
The asymmetry of the limiting distribution \eqref{eq:Kexpexp} (i.e.,
exponential tail on the right and super-exponential tail on the
left) explains the appearance of
the two claims in Theorem~\ref{th:mldsublin}.

\subsection{Asymptotics of $K(\lambda)$ in the space $\varLambda_\bq(n)$}

It is now easy to derive the analogues of Theorems \ref{th:mld}
and~\ref{th:mldsublin} under the measures~$\nu^\bq_n$.

\begin{theorem}\label{th:nu-m}
Suppose that Assumption \textup{\ref{as:main}} holds, with $q\ge0$
and $0\le\beta<1$, and let
$\gamma\in\bigl(0,\frac12(1-\beta)\bigr)$.
\begin{itemize}
\item[\rm (a)] If $q>0$ then there exists a sequence
$(k_n)$ satisfying the asymptotic relation
\begin{equation}\label{eq:knasympt}
k_n\sim\frac{T_q\,\sqrt{n}}{\vartheta_q}\qquad (n\to\infty),
\end{equation}
such that, for any $a>0$,
\begin{equation}\label{eq:log(nu)<0}
\limsup_{n\to\infty}n^{\gamma-1/2}\log\nu^\bq_n
\bigl\{\lambda\in\CP_\bq(n)\colon |K(\lambda)-k_n|>a\mypp
n^{(1-\gamma)/2} \bigr\}<0.
\end{equation}

\item[\rm (b)] If $q=0$ then
\begin{equation}
\label{eq:Ksublin} \limsup_{n\to\infty}n^{\gamma-1/2}\log\nu^\bq_n
\bigl\{\lambda\in\CP_\bq(n)\colon
K(\lambda)<\tfrac12\sqrt{n}\log\log n\text{\ \ or\ \
}K(\lambda)>n^{1-\gamma} \bigr\}<0.
\end{equation}
\end{itemize}
\end{theorem}

\proof (a) Applying Theorem~\ref{th:mld} to the set
$A_z=\{|K(\lambda)-k_*|>c\myp z^{\gamma-1}\bigr\}\subset \CP_\bq$
with $c:=\frac12\myp a\myp\vartheta_q^{1-\gamma}$, we see that $A_z$
satisfies the condition \eqref{eq:limsupmuz} of Proposition
\ref{pr:munu} with $\kappa=1-2\gamma>0$. Hence, setting
$k_n:=k_*(z_n)$ and using \eqref{eq:log(nu)<0a} together with the
property \eqref{eq:xn}, we obtain \eqref{eq:log(nu)<0}, as claimed.
Finally, relation \eqref{eq:knasympt} easily follows from
\eqref{eq:xn} and~\eqref{eq:k0}.

\smallskip
(b) Consider the set
$A_z=\{K(\lambda)<z^{-1}\log\log\tfrac{1}{z}\text{ or }\linebreak[1]
K(\lambda)>k_{\gamma}(z)\}$. By Theorem~\ref{th:mldsublin}, the set
$A_z$ satisfies the condition \eqref{eq:limsupmuz} of
Proposition~\ref{pr:munu}. Moreover, if the asymptotic relation
\eqref{eq:xn} with $q=0$ holds for a sequence $z_n$, then the set
referred to in \eqref{eq:Ksublin} is a subset of $A_{z_n}$, at least
for $n$ large enough, because
\begin{align*}
&z_n^{-1}\log\log\tfrac{1}{z_n}>\tfrac{1}{2}\sqrt{n}\log\log n,\\
&k_{\gamma}(z_n)\le \lceil z_n^{-2(1-\gamma)}\rceil \sim
\left(\frac{6n}{\pi^2}\right)^{1-\gamma}<n^{1-\gamma}\qquad(n\to\infty).
\end{align*}
Thus, the required relation \eqref{eq:Ksublin} follows from
\eqref{eq:log(nu)<0a}.
\endproof

Similarly as before, Theorem \ref{th:nu-m} with $q>0$ implies the
law of large numbers for $K(\lambda)$ under the measure $\nu_n^\bq$,
analogous to Corollary~\ref{cor:LLN_K}.
\begin{corollary}\label{cor:LLN-Kn}
Let Assumption \textup{\ref{as:main}} hold with $q>0$. Then, for any
$\varepsilon>0$,
\begin{equation*}
\lim_{n\to\infty} \nu^\bq_n\biggl\{\lambda\in\CP_\bq(n)\colon
\biggl|\frac{K(\lambda)}{\sqrt{n}}-\frac{T_q}{\vartheta_q}\biggr|>\varepsilon\biggr\}=0.
\end{equation*}
\end{corollary}

If $q=0$ then, under Assumption~\ref{as:main+} with $\tilde{q}>0$
(see~\eqref{eq:assq*+}), one can deduce in a similar fashion the law
of large numbers for $K(\lambda)$: \emph{for any $\varepsilon>0$,}
\begin{equation}\label{eq:KnLLN-q=0}
\lim_{n\to\infty} \nu^\bq_n\biggl\{\lambda\in\CP_\bq(n)\colon
\biggl|\frac{K(\lambda)}{\sqrt{n}\log
n}-\frac{\sqrt{6}\mypp(1-\beta)}{2\pi}\biggr|>\varepsilon\biggr\}=0.
\end{equation}
In fact, an exponential bound for large deviations of $K(\lambda)$
can be obtained by combining \eqref{eq:muxdevkq=0refined} with
Theorem~\ref{th:nu-m}(b), but we omit technical details.

Finally, if $\tilde{q}=0$ and $\tilde{\beta}=0$ in
\eqref{eq:assq*+}, then the classical limit theorem (under the
uniform measure $\nu_n$ on $\CP(n)$)  states
that~\cite{ErdosLehner,Fristedt}
\begin{equation}\label{eq:K-distr2}
\lim_{n\to\infty} \nu_n\biggl\{\lambda\in\CP(n)\colon\frac{\pi
K(\lambda)}{\sqrt{6n}}-\log\frac{\sqrt{6n}}{\pi}\le
t\biggr\}=\exp\myn(-\re^{-t}),\qquad t\in\RR.
\end{equation}
Of course, this result implies the law of large numbers,
\[
\lim_{n\to\infty} \nu_n\biggl\{\lambda\in\CP(n)\colon
\biggl|\frac{K(\lambda)}{\sqrt{n}\log
n}-\frac{\sqrt{6}}{2\pi}\biggr|>\varepsilon\biggr\}=0,
\]
which can be formally considered as the limiting case of
\eqref{eq:KnLLN-q=0} as $\beta\downarrow0$.
\begin{remark}
To be more precise, the results by Erd\H{o}s \& Lehner
\cite{ErdosLehner} and Fristedt \cite{Fristedt}, quoted above as
formulas \eqref{eq:Kexpexp} and \eqref{eq:K-distr2}, are technically
about the maximal part $\lambda_1$, but due to the invariance of the
measures $\mu_z$ and $\nu_n$ under \emph{conjugation} of Young
diagrams (whereby columns become rows and vice versa; see also
Section~\ref{sec:lsalt}), the random variable $\lambda_1$ has the
same distribution as the number of parts $K(\lambda)$.
\end{remark}

\section{Limit shape of the minimal difference partitions}\label{sec:pdls}

\subsection{The parametric family of limit shapes}\label{sec:duality}
Mutual independence of the random variables $(D_j(\lambda))_{j=1}^k$
with respect to the measure $\mu^\bq_{z,\myp k}$ (see
Lemma~\ref{lm:D_j}) provides an easy way to find the limit shape for
MDPs as $z\downarrow0$. It is natural to allow the maximal part $k$
to grow to infinity as $z$ approaches $0$, where the correct growth
rate, as suggested by Theorem~\ref{th:mld}, is of order $z^{-1}$
when $q>0$ and possibly faster, by a logarithmic factor, when $q=0$.
It turns out that if the condition \eqref{eq:assq*} holds and
$\lim_{z\downarrow 0}zk=T<\infty$ then $\mu^\bq_{z,k}$-typical
partitions $\lambda\in \CP_\bq(\DOT,k)$ concentrate around the limit
shape determined by the function
\begin{equation}\label{eq:LS}
\varphi_T(t;q):=\Biggl\{\begin{array}{ll} \displaystyle
q\myp(T-t)
+\log\frac{1-\re^{-T}}{1-\re^{-t}}\mypp,\quad &0<t\le T,\\
0,&t\ge T.
\end{array}
\end{equation}
If $q=0$ then the expression \eqref{eq:LS} is reduced to
\begin{equation}\label{eq:LS0}
\varphi_T(t;0)=\Biggl\{\begin{array}{ll} \displaystyle
\log\frac{1-\re^{-T}}{1-\re^{-t}}\mypp,\quad &0<t\le T,\\
0,&t\ge T,
\end{array}
\end{equation}
which coincides, as one could expect, with the limit shape of plain
(unrestricted) partitions subject to the condition $zk\to T$
(see~\cite{VY-MMJ}).

If $q=0$, one can also allow $zk$ to grow slowly to infinity as
$z\downarrow 0$ (which is actually a typical behaviour), whereby the
limit shape is given by the formula
\[
\varphi_\infty(t;0)=-\log\myn(1-\re^{-t})
\]
(which is formally consistent with \eqref{eq:LS0} if we set
$T=\infty$).

\begin{figure}[h]
\begin{picture}(180,250)
\put(-20,8){
\put(127,0){\includegraphics[width=8cm]{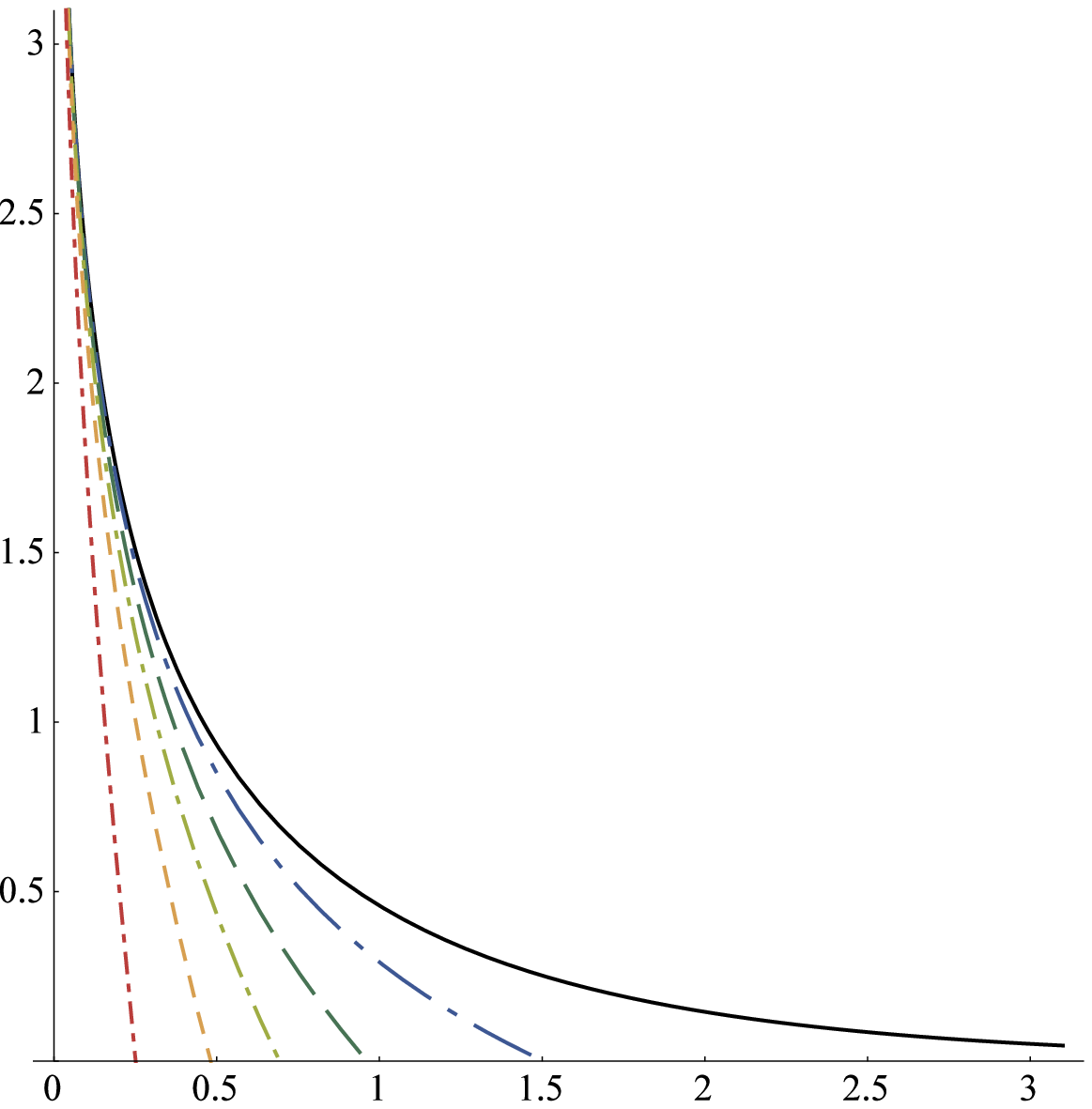}}
\put(132,234){\mbox{\footnotesize$y$}}
\put(354,2.3){\mbox{\footnotesize$x$}} \put(290,75){
\put(0,0){\includegraphics[width=7cm]{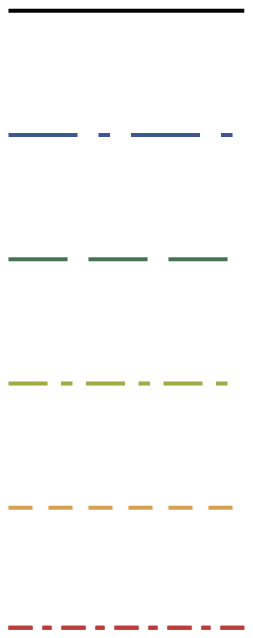}}
\put(50,107){\mbox{\footnotesize $q=0,\;T_q=\infty$}}
\put(50,85){\mbox{\footnotesize$q=\frac16,\:T_q\doteq 1.505482$}}
\put(50,63){\mbox{\footnotesize$q=\frac12,\:T_q\doteq0.962424$}}
\put(50,42){\mbox{\footnotesize$q=1,\:T_q\doteq0.693147$}}
\put(50,21){\mbox{\footnotesize$q=2,\:T_q\doteq0.481212$}}
\put(50,0){\mbox{\footnotesize$q=6,\:T_q\doteq 0.250914$}} } }
\end{picture}
\caption{The parametric family of the limit shapes \eqref{eq:LS-Tq}
plotted in the Cartesian coordinates $x=t$ and
$y=-tq-\log\myn(1-\re^{-t})$ (see equation~\eqref{eq:ls-xy}).}
\label{fig:ls}
\end{figure}

Another simplification of formula \eqref{eq:LS} worth mentioning
occurs for $q>0$ and $T=T_q$ (see~\eqref{eq:deftd}), which
determines the typical behaviour of the number of parts in this case
(see Theorem~\ref{th:mld} and the asymptotic formula~\eqref{eq:k0});
here, the limit shape \eqref{eq:LS} is reduced to
\begin{equation}\label{eq:LS-Tq}
\varphi_{T_q}(t;q)=\left\{\!\begin{array}{ll}
-tq-\log\myn(1-\re^{-t}),\quad&0<t\le T_q\myp,\\[.3pc]
\hphantom{-} 0,&t\ge T_q\myp.
\end{array}\right.
\end{equation}
This coincides with the limit shape found by Comtet et al.\
\cite[Eq.~(19)]{Comtet2}, \cite[Eq.~(11)]{Comtet3}. The limit shape
\eqref{eq:LS-Tq} is illustrated in Fig.~\ref{fig:ls} for various
values of parameter~$q\ge 0$ using Cartesian coordinates $x=t$,
$y=-tq-\log\myn(1-\re^{-t})$, whereby \eqref{eq:LS-Tq} takes the
form
\begin{equation}\label{eq:ls-xy}
\re^{-y}=\re^{qx}\myp(1-\re^{-x}),
\end{equation}
which was already mentioned in Section \ref{sec:1.3}
(see~\eqref{eq:ls-intro-xy}).

\subsection{The limit shape in the spaces \mbox{$\CP_{\bq}(\protect\DOT,k)$} and \mbox{$\CP_{\bq}(n,k_n)$}}
The exact statement is as follows. Recall that the notation
$k_\gamma(z)$ is defined in~\eqref{eq:defkappa}.

\begin{theorem}\label{th:muxmls}
Let Assumption \textup{\ref{as:main}} hold, with $q\ge0$. Then for
every $t>0$ and any\/ $\varepsilon>0$, uniformly in $k=k(z)\in\NN$
such that\/ $\lim_{z\downarrow0}zk(z) =T\in(0,\infty)$,
\begin{equation}\label{eq:pdp:expbound}
\limsup_{z\downarrow0} z\log \mu^\bq_{z,\myp
k}\bigl\{\lambda\in\CP_\bq(\DOT,k)\colon |z\mypp
Y_\lambda(t/z)-\varphi_T(t;q)| > \varepsilon \bigr\} <0.
\end{equation}
Furthermore, if $q=0$ and\/ $\lim_{z\downarrow0}zk(z)=\infty$ but
$k(z)\le k_\gamma(z)$, with some
$\gamma\in\bigl(0,\frac12(1-\beta)\bigr)$, then
\eqref{eq:pdp:expbound} holds with $\varphi_\infty(t;0)$ in place of
$\varphi_T(t;q)$.
\end{theorem}

\proof First, let us show that the curve $t\mapsto\varphi_T(t;q)$ is
the limit of the $\mu^\bq_{z,\myp k}$-mean of the scaled Young
diagrams, that is, for every $t>0$
\begin{equation}\label{eq:pdp:exphi}
\lim_{z\downarrow0}z\mypp \BE^\bq_{z,\myp k}[\myp Y_\lambda(t/z)]=
\varphi_T(t;q).
\end{equation}
To this end, using the definition \eqref{eq:Phi} and the formula
\eqref{eq:E(D)}, we can write, for $0<t<T$,
\begin{equation}\label{eq:pdp:exPhi}
\BE^\bq_{z,\myp k}[\myp Y_\lambda(t/z)]= \sum_{t/z<j\le
k}\BE^\bq_{z,\myp k} [D_j(\lambda)] =\sum_{t/z<j\le k}q_{k-j}
+\sum_{t/z<j\le k}\frac{\re^{-jz}}{1-\re^{-jz}}.
\end{equation}
According to \eqref{eq:assq*}, for $T<\infty$ and $q\ge 0$ the first
sum in \eqref{eq:pdp:exPhi} is asymptotically evaluated as follows
\begin{equation}\label{eq:sumq}
\sum_{t/z<j\le k}q_{k-j}=Q_{k-\floor{t/z}} =
q\left(k-\floor{t/z}\right)+O\bigl((k-t/z)^\beta\bigr)=qz^{-1}(T-t)+o(z^{-1})
\end{equation}
since $zk\to T$ as $z\downarrow0$. If $T=\infty$ and $q=0$, then for
$k\le k_\gamma(z)$ one has $Q_{k-\floor{t/z}}\le
Q_{k_\gamma(z)-\floor{t/z}}=O(z^{-1+2\gamma})=o(z^{-1})$ by
Lemma~\ref{lm:k4} (see~\eqref{eq:lm3.5-3}).

For the second sum in \eqref{eq:pdp:exPhi}, we get (e.g., via the
Euler--Maclaurin sum formula) that
\begin{align}\notag
\sum_{t/z<j\le k}\frac{\re^{-jz}}{1-\re^{-jz}}&\sim \int_{t/z}^{k}
\frac{\re^{-xz}}{1-\re^{-xz}}\,\rd{x}\\
\notag &=
z^{-1}\log\mynn (1-\re^{-zx})\bigr|_{t/z}^{k}\\
\notag
&=z^{-1}\log\frac{1-\re^{-zk}}{1-\re^{-t}}\\
&\sim z^{-1}\log\frac{1-\re^{-T}}{1-\re^{-t}}\qquad (z\downarrow0).
\label{eq:sumchi}
\end{align}
The same calculation is valid when $zk\to\infty$, with the change of
$\re^{-T}$ to $0$. Thus, on substituting the estimates
\eqref{eq:sumq} and \eqref{eq:sumchi} into \eqref{eq:pdp:exPhi} we
get \eqref{eq:pdp:exphi}.

To obtain the exponential bound \eqref{eq:pdp:expbound}, we use a
standard technique often applied in similar problems (see, e.g.,
\cite{Corteeletal}). Suppose that $z k\to T\in(0,\infty]$, and fix
$t\in(0,T)$ and $\varepsilon>0$. In what follows, we always assume
that $z$ is small enough so that $zk>t$ and
\begin{equation}\label{eq:(ii)}
\bigl|z\mypp\BE_{z,\myp k}[\myp
Y_\lambda(t/z)]-\varphi_T(t;q)\bigr|<\tfrac12\myp \varepsilon.
\end{equation}
Then for any $u\in(0,t)$
\begin{align}\notag
\mu^\bq_{z,\myp k}\bigl\{z\myp
Y_\lambda(t/z)-\varphi_T(t;q)>\varepsilon\bigr\}&\le\mu^\bq_{z,\myp
k}\bigl\{Y_\lambda(t/z)\ge
\BE^\bq_{z,\myp k}[\myp Y_\lambda(t/z)]+\tfrac12\myp z^{-1}\varepsilon\bigr\}\\[.4pc]
\notag &\le\exp\bigl(-u\mypp\BE^\bq_{z,\myp k}[\myp
Y_\lambda(t/z)]-\tfrac12\myp uz^{-1}\varepsilon\bigr)
\,\BE^\bq_{z,\myp k}\bigl[\exp(u Y_\lambda(t/z))\bigr]\\[.2pc]
\label{eq:pdp:muxbnd} &=\exp\bigl(-\tfrac12\myp
uz^{-1}\varepsilon\bigr)\prod_{t/z<j\le k} \BE^\bq_{z,\myp
k}\bigl[\exp\bigl(uD_j-u\mypp\BE^\bq_{z,\myp k}(D_j)\bigr)\bigr],
\end{align}
where the first inequality is a consequence of
assumption~\eqref{eq:(ii)}, the second is the exponential Markov
inequality, and the last line follows from the additive structure of
$Y_\lambda(t)$ and independence of~$(D_j)_{j=1}^k$.

Suppose that, for some $w\in(0,1)$ that will be specified later,
\begin{equation}\label{eq:v}
0<u\le \log\!\left(1+\frac{w}{h(t)}\right)=:v(w),
\end{equation}
where we put for short
\begin{equation}\label{eq:Psi}
h(t):=\frac{\re^{-t}}{1-\re^{-t}},\qquad t\in(0,\infty).
\end{equation}
Then for $j\ge t/z$ we have
\[
(\re^u-1)\,h(zj)\le(\re^u-1)\,h(t)\le w.
\]
Applying the elementary inequalities
\begin{alignat*}{2}
-\log\myn(1-x)&{}\le -x w^{-1}\log\myn(1-w)&\qquad&(0<x\le w),\\[.2pc]
\re^u-1&{}\le u\myp v^{-1}(\re^v-1)&\qquad&(0<u\le v),
\end{alignat*}
with $x:=(\re^u-1)\,h(zj)$ and $v:=v(w)$ (see~\eqref{eq:v}), we
obtain
\[
-\log\bigl(1-(\re^u-1)\,h(zj)\bigr)\le u\mypp y(w)\,h(zj),\qquad
y(w):=\frac{-\log\myn(1-w)}{h(t)\,v(w)}.
\]
Hence, for $u\le \min\{v(w),t\}\le jz$
\begin{align}
\notag \log\!\left( \BE^\bq_{z,\myp
k}\bigl[\exp(uD_j-u\mypp\BE^\bq_{z,\myp k}D_j)\bigl]\right)
&=\log\frac{1-\re^{-zj}}{1-\re^{u-zj}}-u\mypp h(zj)\\[.1pc]
\notag
&=-\log\bigl[1-(\re^u-1)\mypp h(zj)\bigr]-u\mypp h(zj)\\[.2pc]
\label{eq:logBE} &\le u\mypp \bigl(y(w)-1\bigr)\mypp h(zj).
\end{align}
Substituting \eqref{eq:logBE} into \eqref{eq:pdp:muxbnd} and
recalling \eqref{eq:sumchi}, we obtain
\begin{align}
\notag z\log \mu^\bq_{z,\myp
k}\bigl\{\lambda\in\CP_\bq(\DOT,k)\colon
z\myp Y_\lambda(t/z)&-\varphi_T(t; q)>\varepsilon\bigr\}\\[.2pc]
\notag
&\le \frac{u}{2}\mypp\bigl(-\varepsilon+(y(w)-1)\,\varphi_T(t;0)\bigr)\\
&\le
\frac{v(w)}{2}\mypp\bigl(-\varepsilon+(y(w)-1)\,\varphi_T(t;0)\bigr)
\label{eq:logmu}.
\end{align}
Since $y(w)\to 1$ as $w\downarrow0$, we can choose $w$ small enough
to make the right-hand side of \eqref{eq:logmu} negative. This
yields the desired bound for the probability of positive deviations
in~\eqref{eq:pdp:expbound}. The probability of negative deviations
is estimated in the same fashion.
\endproof

We are now in a position to state and prove our first main result.

\begin{theorem}\label{th:nunmls}
Suppose that Assumption \textup{\ref{as:main}} is satisfied with
some $q\ge0$, and let $k_n\to\infty$ so that $k_n/\sqrt{n}\to \tau$
as $n\to\infty$, for some $\tau\in(0,\sqrt{2/q}\mypp)$, with the
right bound understood as $+\infty$ when $q=0$. Then, for every
$t_0>0$ and any $\varepsilon>0$,
\begin{equation}\label{eq:limsupnunmsup}
\limsup_{n\to\infty}\,\frac{1}{\sqrt{n}}\log \nu^\bq_{n,\myp
k_n}\!\left\{\lambda\in\CP_\bq(n,k_n)\colon \sup_{t\ge
t_0}\left|z_nY_\lambda(t/z_n) -
\varphi_{T_*}(t; q)\right|
> \varepsilon \right\}<0,
\end{equation}
where
$T_*=T_*(\tau;q)>0$
is the \textup{(}unique\textup{)} solution of the equation
\begin{equation}\label{eq:TofB}
\tau\myp
\vartheta_q(T_*)=T_*
\end{equation}
and
\begin{equation}\label{eq:zn}
z_n:=\frac{T_*}{k_n}\sim
\frac{T_*}{\tau\sqrt{n}}.
\end{equation}
Furthermore, if\/  $q=0$ then the result \eqref{eq:limsupnunmsup} is
also valid in the case $k_n/\sqrt{n}\to\infty$ under the additional
condition $\limsup_{n\to\infty}k_n^{\beta+1}/n^{1-\delta}<\infty$
for some $\delta\in(0,1)$, with $T=\infty$ and\/
$\vartheta_0(\infty)=\pi/\sqrt{6}$.
\end{theorem}

\begin{remark}
The assumption $\tau^2< 2/q$ in Theorem~\ref{th:nunmls} arises
naturally, because if $\lambda\in\CP_\bq(n,k)$ then, due to the MDP
condition~\eqref{eq:restr}, we must have $n\ge s_k=\frac12\myp
qk^2+O(k^{1+\beta})$, which yields $\tau^2\le 2/q$. The boundary
case $\tau^2=2/q$ can in principle be realized, but both the
formulation and analysis should be more accurate, so we do not
consider it with the exception of the important special case $q=0$
when additional difficulties can be treated without much effort.
\end{remark}

\proof[Proof of Theorem~\textup{\ref{th:nunmls}}] Note that  the
equation \eqref{eq:TofB} can be rewritten as
\[
\frac{\Li_2(1-\re^{-T})}{T^2}=\frac{1}{\tau^2}-\frac{q}{2}
\]
with the left-hand side decreasing from $+\infty$ to $0$ as $T$
grows from $0$ to $+\infty$, so its positive solution $T=T_*$ always
exists (and is unique) for any $\tau\in(0,\sqrt{2/q}\mypp)$.

For $0<t_0\le t<T\le\infty$ and $\varepsilon>0$, denote
\begin{align*}
A_{z,\myp k}(t,\varepsilon):&{}=\left\{ \lambda\in\CP(\DOT,k)\colon
\left|z\myp Y_\lambda(t/z)- \varphi_T(t;q) \right|
> \varepsilon \right\},\\[.2pc]
\widehat{A}_{z,\myp k}(t_0,\varepsilon):&{}=\left\{
\lambda\in\CP(\DOT,k)\colon \sup\nolimits_{t\ge t_0}\left|z\myp
Y_\lambda(t/z)- \varphi_T(t;q) \right|
> \varepsilon \right\}.
\end{align*}
Given $t_0>0$ and $\varepsilon>0$,  define $t_i$ recursively by
$\varphi_T(t_i;q)=\varphi_T(t_{i-1},q)-\varepsilon/2$ until
$\varphi_T(t_{s-1})-\varepsilon/2$ becomes negative for some $s$. By
construction,
\begin{equation}\label{eq:hatCA}
{\textstyle\bigcup\limits_{i=0}^{s-1}} A_{z,\myp
k}(t_i,\varepsilon/2) \supset \widehat{A}_{z,k}(t_0,\varepsilon),
\end{equation}
because both $Y_\lambda(t)$ and $\varphi_T(t,q)$ decrease as
functions of $t$.

Now, we aim to apply Theorem \ref{th:muxmls} and
Proposition~\ref{pr:mumnum}. To this end, in the case $T<\infty$
take $k(z)$ to be any integer-valued function such that $zk(z)\to
T$; in the case $T=\infty$ (arising for $q=0$, $\tau=\infty$) let
$k(z):=k_n$ for $z\in(\pi/\sqrt{6(n+1)},\pi/\sqrt{6n}]$, where the
sequence $(k_n)$ is referred to in the theorem. In the latter case
one can write $k(z)=k_{\floor{\pi^2/6z^2}}$, and the additional
requirement $\limsup_{n\to\infty}k_n^{\beta+1}/n^{1-\delta}<\infty$
 combined with \eqref{eq:sksim} implies $s_{k_n}=O(n^{1-\delta})$ which can be rewritten
as $s_{k(z)}=O(z^{-2+2\delta})$. Thus, for $\gamma\in(0,\delta)$ and
$z$ small enough one has $s_{k(z)}<z^{-2+2\gamma}$, and thus
$k(z)<k_\gamma(z)$ (see \eqref{eq:defkappa}).

Hence, Theorem \ref{th:muxmls} implies that for any $t\in(0,T)$ and
$\varepsilon>0$
\begin{equation}\label{eq:limsupCA}
\limsup_{z\downarrow 0} z\myp \log\mu^\bq_{z,\myp
k(z)}\bigl(A_{z,\myp k(z)}(t,\varepsilon)\bigr)<0.
\end{equation}
It follows from the asymptotic bound \eqref{eq:limsupCA} (applied
with $\varepsilon/2$ instead of~$\varepsilon$) and the inclusion
\eqref{eq:hatCA} that for any $t_0>0$
\begin{equation*}
\limsup_{z\downarrow 0}z\myp \log\mu^\bq_{z,\myp k(z)}
\bigl(\widehat{A}_{z,\myp k(z)}(t_0,\varepsilon)\bigr)<0.
\end{equation*}
Furthermore, $k_n/\sqrt{n}\to \tau =T/\vartheta_q(T)$ as
$n\to\infty$; if $q=0$ and $\tau=\infty$ then
$k_n/k(\pi/\sqrt{6n})=1$ by construction and $z^{2/(\beta+1)}k(z)\to
0$ by the assumption
$\limsup_{n\to\infty}k_n^{\beta+1}/n^{1-\delta}<\infty$ As a result,
by Proposition~\ref{pr:mumnum} there exists a sequence $(\tilde
z_n)$ such that for any $t_0>0$
\[
\limsup_{n\to\infty}\frac{1}{\sqrt{n}}\log \nu^\bq_{n,\myp k_n}
\bigl(\widehat{A}_{\tilde z_n,\myp k_n}(t_0,\varepsilon)\bigr)<0.
\]
Finally, it is easy to see that the
sequence $(\tilde z_n)$ of Proposition \ref{pr:mumnum} and the
sequence $(z_n)$ defined by \eqref{eq:zn} are asymptotically
equivalent so can be interchanged in
\eqref{eq:limsupnunmsup}.
\endproof

\subsection{The limit shape in the spaces $\CP_{\bq}$ and $\CP_{\bq}(n)$}
Recall that the function $\varphi_{T_q}(t;q)$ is given by
\eqref{eq:LS-Tq}, where $T_q$ is defined as the unique solution of
the equation \eqref{eq:deftd}.

\begin{theorem}\label{th:muxls}
Let Assumption \textup{\ref{as:main}} hold, with $q\ge0$. Then for
every $t>0$ and any $\varepsilon>0$, $\delta>0$
\begin{equation}\label{eq:muxls}
\limsup_{z\downarrow0}z^{1-\delta}\log
\mu^\bq_z\bigl\{\lambda\in\CP_\bq\colon|z\myp
Y_\lambda(t/z)-\varphi_{T_q}(t;q)|>\varepsilon
\bigr\}<0.
\end{equation}
\end{theorem}

\proof Let $A_z\subset\CP_\bq$ be the set on the left-hand side
of~\eqref{eq:muxls}. Then
\begin{equation*}
\mu^\bq_z(A_z)=\sum_{k=0}^\infty \mu^\bq_z(A_z\cap\CP_\bq(\DOT,k))
=\sum_{k=0}^\infty \mu^\bq_{z,\myp
k}(A_z)\,\mu^\bq_z\{K(\lambda)=k\}\,.
\end{equation*}

Suppose that $q>0$. Take $\gamma\in(0,\min\{\delta/2,(1-\beta)/2\})$
and set $\mathcal{I}_z:=\{k\in\NN\colon |k-k_*|> z^{\gamma-1}\}$,
where $k_* = k_*(z)$ is defined in \eqref{eq:def-k0}. Then
\begin{align}
\notag \mu^\bq_z(A_z)&\le \Biggl(\,\sum_{k\in\mathcal{I}_z}
+\sum_{k\notin\mathcal{I}_z}\Biggr)\,\mu^\bq_{z,\myp
k}(A_z)\,\mu^\bq_z\{K(\lambda)=k\}\\
\label{eq:mu(A)} &\le \mu^\bq_z\{K(\lambda)\in\mathcal{I}_z\}+
\max_{k\notin\mathcal{I}_z}\mu^\bq_{z,\myp k}(A_z).
\end{align}
Using the elementary inequality \eqref{eq:loga+b}, we get from
\eqref{eq:mu(A)}
\[
\log \mu^\bq_z(A_z)\le \log
2+\max\Bigl\{\log\mu^\bq_z\{K(\lambda)\in\mathcal{I}_z\}, \, \log
\max_{k\notin\mathcal{I}_z}\mu^\bq_{z,\myp k}(A_z)\Bigr\}.
\]
Multiplying this by $z^{1-\delta}$ and applying Theorems
\ref{th:mld} and \ref{th:muxmls}, we obtain~\eqref{eq:muxls}.

If $q=0$ then we set $\mathcal{I}_z:= \{k\in\NN\colon
k<z^{-1}\log\log\tfrac1z\}\cup\{k\in\NN\colon k>k_\gamma(z)\}$ and
repeat the above argumentation with a reference to Theorem
\ref{th:mldsublin} instead of Theorem~\ref{th:mld}.
\endproof

Our second main result describes the limit shape under the measure
$\nu^\bq_n$, that is, without any restriction on the number of
parts.

\begin{theorem}\label{th:nu'nls}
Let Assumption \textup{\ref{as:main}} be satisfied, with $q\ge0$.
Then for every $t_0>0$ and any $\varepsilon>0$ and $\delta>0$, we
have
\[
\limsup_{n\to\infty}\, n^{\delta-1/2}\log
\nu^\bq_n\biggl\{\lambda\in\CP_\bq(n)\colon \sup_{t\ge
t_0}\myp\bigl|z_nY_\lambda(t/z_n)
-\varphi_{T_{q}}(t;q)\bigr|>\varepsilon \biggr\}<0,
\]
where $z_n=\vartheta_q/\sqrt{n}$, with $\vartheta_q$ given
by~\eqref{eq:theta-intro-alt}.
\end{theorem}

\proof The claim follows from Theorem \ref{th:muxls} and Proposition
\ref{pr:munu} by the same argumentation as that used to derive
Theorem \ref{th:nunmls} from Theorem \ref{th:muxmls} and
Proposition~\ref{pr:mumnum}.
\endproof

\subsection{Ground state}\label{sec:4.4}

\begin{figure}[h]
\begin{center}
\begin{picture}(300,230)
\put(0,0){\includegraphics[width=0.33\textwidth]{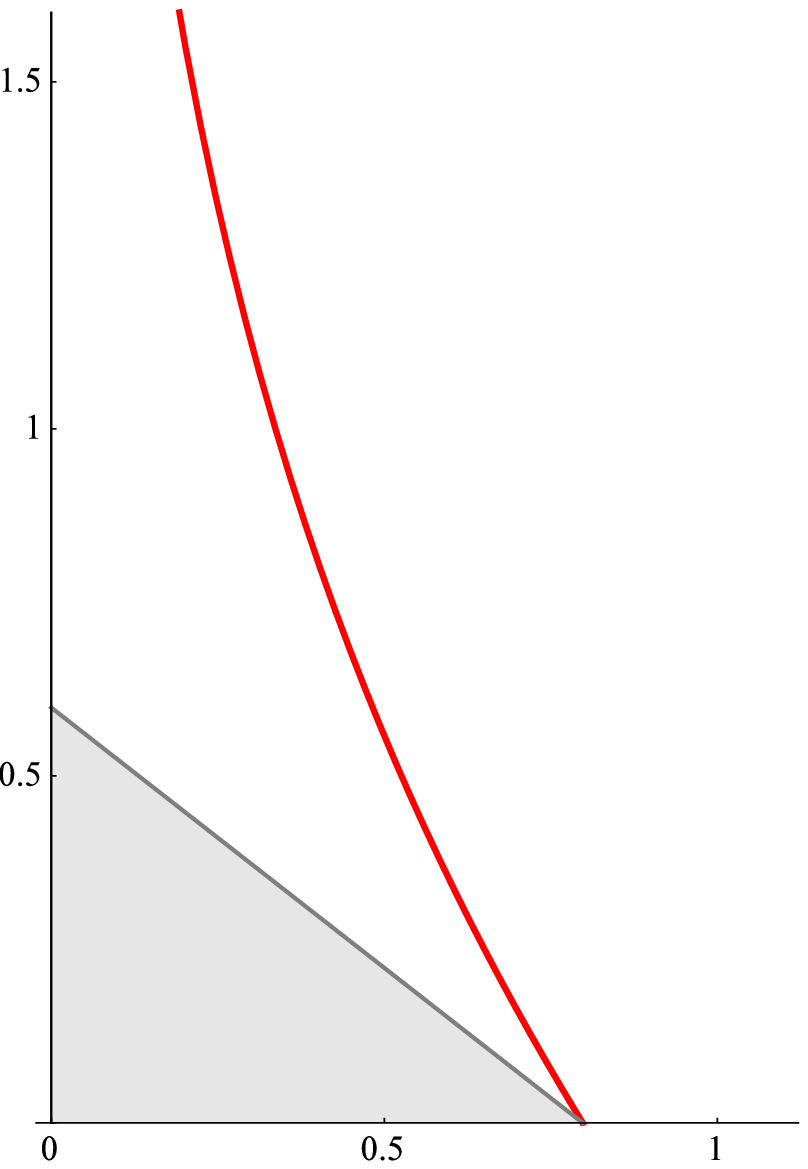}}
\put(200,0){\includegraphics[width=0.33\textwidth]{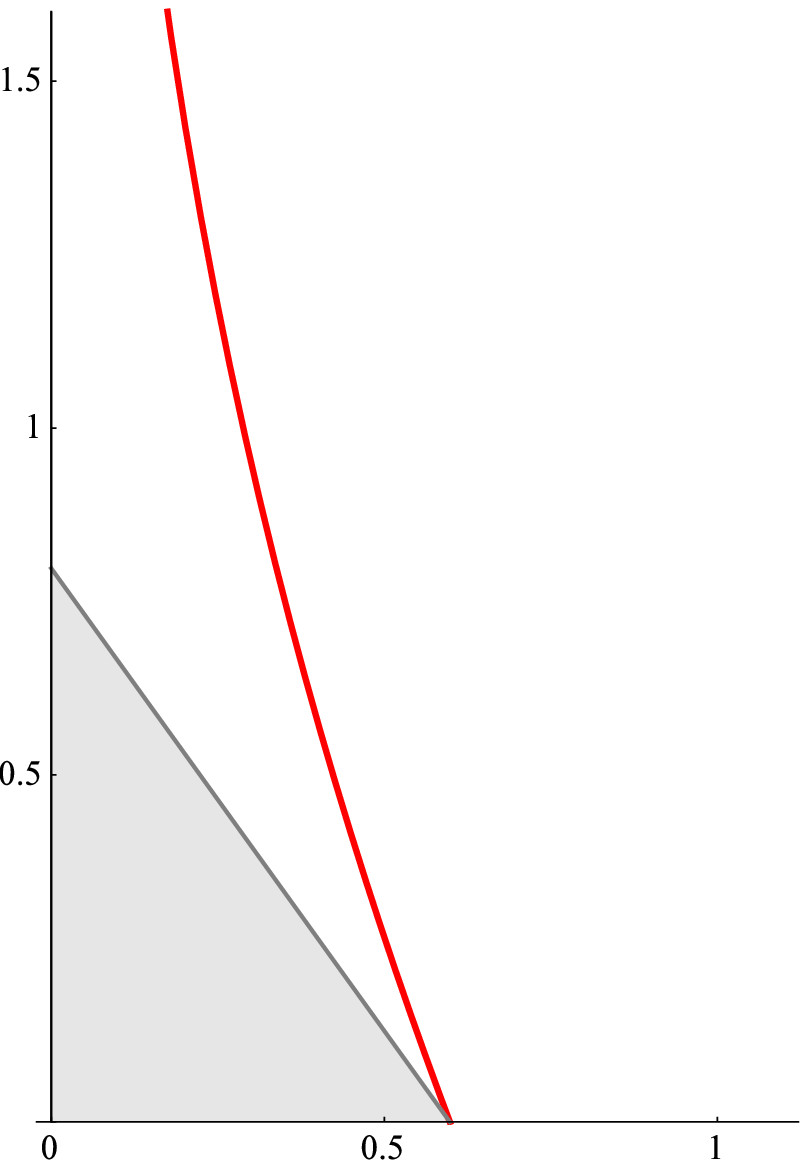}}
\put(110,150){\mbox{\small $q\mapsto 1/q$}}
\put(103,-2){\mbox{\footnotesize $T_q$}}
\put(-8,83){\mbox{\footnotesize $q\myp T_q$}}
\put(173,108){\mbox{\footnotesize $q^{-1}\myp T_{1/q}$}}
\put(280,-2){\mbox{\footnotesize $T_{1/q}$}}
\put(34,27){\mbox{\small $\Delta_q$}}%
\put(227,32){\mbox{\small $\Delta_{1/q}$}}
\end{picture}
\end{center}
\caption{The duality under the transformation $q\mapsto 1/q$
illustrated for $q=\frac43$\myp, where $T_{4/3}\doteq 0.598382$ and
$T_{3/4}\doteq 0.797842$. The ground state triangles $\Delta_q$ and
$\Delta_{1/q}$ (shaded in grey) are obtained from one another by
reflection about the main coordinate diagonal. Thus, in line with
Lemma~\ref{lm:duality_identity}, $T_{1/q}=q\myp T_{q}$ and,
equivalently, $T_{q}=q^{-1}\myp T_{1/q}$; in particular,
$T_{3/4}=\frac43\myp T_{4/3}$. The solid curves (red in the online
version) show the limit shape graphs. According to formula
\eqref{eq:conservation-law}, the areas under the limit shapes sum up
to $\zeta(2)=\frac16\myp\pi^2$.}\label{fig:dual}
\end{figure}

\noindent
Observe that, for $q>0$, the area beneath the limit shape
$t\mapsto \varphi_{T_{q}}(t;q)$
featured in Theorems \ref{th:muxls}
and~\ref{th:nu'nls} contains a
right-angled triangle $\Delta_q$
(shaded in Fig.~\ref{fig:dual}) obtained in the limit from the
(rescaled) partitions in $\CP_\bq(n)$ satisfying the \emph{hard
version} of the \strut{}MDP restrictions \eqref{eq:restr}, that is,
when all inequalities are replaced by equalities. Thus,
we can say that the triangle
$\Delta_q$ represents the \emph{ground state} of the $\MDP(\bq)$
system, while the remaining part of the limit shape corresponds to
additional degrees of freedom in a $\nu_n^\bq$-typical partition.
Note that, according to the $\nu^\bq_n$-typical asymptotic behaviour
of $K(\lambda)$ described in Corollary \ref{cor:LLN-Kn}, under the
scaling of Theorem~\ref{th:nu'nls} the horizontal leg of the
triangle $\Delta_q$ is identified as $T_q$. On the other hand, by
the condition \eqref{eq:q*} the slope of the
hypotenuse of the triangle is given
by~$q$, therefore the vertical leg of $\Delta_q$ is found to be
$qT_q$; in particular, the area of $\Delta_q$ is $\frac12\myp
qT_q^2$. Since the total area of the limit shape is $\vartheta^2_q$
(see \eqref{eq:theta-intro-alt}), the area of the ``free'' part is
given by
\begin{equation}\label{eq:free1}
\vartheta^2_q-\tfrac12\myp qT_q^2=\Li_2(1-\re^{-T_q}).
\end{equation}

This remark helps to clarify the duality identity
\eqref{eq:Tq-identity} of Lemma \ref{lm:duality_identity}. To this
end, consider the triangle $\tilde{\Delta}_q$ obtained from
$\Delta_q$ by reflection about the principal coordinate diagonal,
that is, with legs $qT_q$ (horizontal) and $T_q$ (vertical). This
triangle may serve as the ground state of a suitable
$\MDP(\tilde{\bq})$ ensemble. The slope of the hypotenuse of
$\tilde{\Delta}_q$ is $1/q$, which therefore gives the limiting gap
of the space $\MDP(\tilde{\bq})$. But according to the previous
considerations, the legs of the triangle $\tilde{\Delta}_q$ must
have the lengths $T_{1/q}$ (horizontal) and $(1/q)\mypp T_{1/q}$
(vertical). Comparing these values, we arrive at the
identity~\eqref{eq:Tq-identity} (see Fig.~\ref{fig:dual}).

Finally, despite the limit shape of the ensemble $\MDP(\tilde{\bq})$
contains the triangle $\tilde{\Delta}_q=\Delta_{1/q}$ of the same
area as $\Delta_q$, the ``free'' area changes to
(cf.~\eqref{eq:free1})
$$
\Li_2(1-\re^{-T_{1/q}})=\Li_2(1-\re^{-qT_q})=\Li_2(\re^{-T_q}).
$$
Moreover, according to the identity \eqref{eq:Li2=}, the total area
of the free parts in the limit shapes with $q$ and $1/q$ is given by
\strut$\frac16\myp \pi^2-qT_q^2$, which in turn implies that the
total area of both limit shapes including the ground state triangles
equals \strut$\frac16\myp \pi^2$,
\begin{equation}\label{eq:conservation-law}
\vartheta_q^2+\vartheta_{1/q}^2=\frac{\pi^2}{6},
\end{equation}
which may be interpreted as the (asymptotic) law of conservation of
total energy in \emph{dual systems}, that is, with limiting gaps $q$
and $1/q$. It would be interesting to find a
physical explanation of this identity.

\section{Alternative approach to the limit shape}\label{sec:lsalt}

Iterating the MDP condition \eqref{eq:restr}, for any partition
$\lambda=(\lambda_1,\dots,\lambda_k)\in\CP_\bq(n,k)$
we get the explicit constraints on its parts,
\begin{equation}\label{eq:MDP-iterated}
\lambda_i\ge q_0+\dots+q_{k-i} \qquad (i=1,\dots,k).
\end{equation}
Note that equalities in \eqref{eq:MDP-iterated} correspond to what
was called the ``ground state'' in the discussion in
Section~\ref{sec:4.4}. Now, it is natural to ``subtract'' the ground
state by shifting the parts of
$\lambda\in\CP_\bq(n,k)$ so as to
lift the constraints \eqref{eq:MDP-iterated} (apart from the default
condition that all parts are not smaller than $1$). Specifically,
consider the mapping
\begin{equation}\label{eq:map}
\mathfrak{I}\colon
\lambda=(\lambda_1,\dots,\lambda_k)\to\rho=(\rho_1,\dots,\rho_k)
\end{equation}
defined by
\begin{equation}\label{eq:rho}
\rho_i:=\lambda_i+1-q_0-\dots-q_{k-i}\ge1\qquad (i=1,\dots,k).
\end{equation}

\begin{remark}
The mapping \eqref{eq:rho} is the (shifted) inverse of the
generalized Sylvester transformation mentioned in
Remark~\ref{rm:1.2}.
\end{remark}

Using \eqref{eq:rho} and \eqref{eq:restr}, note that
\begin{equation*}
\rho_i-\rho_{i+1}=\lambda_i-\lambda_{i+1}-q_{k-i}\ge 0,\qquad
\rho_k=\lambda_k+1-q_0\ge 1,
\end{equation*}
and, recalling the notation~\eqref{eq:s},
\begin{equation*}
r:=\sum_{i=1}^k \rho_i=\sum_{i=1}^k
\lambda_i +k-\sum_{i=1}^k i\myp q_{k-i}=n+k-s_k\ge k,
\end{equation*}
where $n\ge s_k$ as long as the set $\CP_\bq(n,k)$ is not empty.
Hence, $\rho=\mathfrak{I}(\lambda)$ is a partition of the same
length $k$ and the new weight $r=n+k-s_k$, but with \emph{no
constraints} on its parts; that is, $\rho\in \CP(r,k)$. Moreover, it
is evident that the mapping \eqref{eq:map} is a bijection of
$\CP_\bq(n,k)$ onto $\CP(r,k)$, for each $k\in\NN$ and any $n\ge
s_k$. In particular, if $\nu^\bq_{n,k}$ is the uniform measure on
$\CP_\bq(n,k)$ then the push-forward
$\mathfrak{I}^*\myn\nu^\bq_{n,k}=\nu^\bq_{n,k}\circ
\mathfrak{I}^{-1}$ is the uniform measure on $\CP(r,k)$.

This observation furnishes a more straightforward way to finding the
limit shape of partitions in the MDP spaces $\CP_\bq(n,k)$ and
$\CP_\bq(n)$. The heuristic idea is as follows. Consider a partition
$\lambda \in\CP_\bq(n,k_n)$, where $k_n\sim \tau\sqrt{n}$ with
$0<\tau<\sqrt{2/q}$ (cf.\ the hypothesis in
Theorem~\ref{th:nunmls}). On account of the asymptotics
\eqref{eq:sksim}, for the weight of $\rho=\mathfrak{I}(\lambda)$
this gives
\begin{equation}\label{eq:r:n}
r=n+k_n-s_{k_n}\sim \bigl(1-\tfrac12\myp q\myp \tau^2\bigr)\myp n=
b^2 n,
\end{equation}
where \begin{equation}\label{eq:bdef} \qquad
b=b(q;\tau):=\textstyle{\sqrt{1-\vphantom{T^|}\smash{\tfrac12}\myp
q\myp \tau^2}>0}.
\end{equation}
In particular, $k_n\sim (\tau/b)\mypp\sqrt{r}$. Suppose now that the
limit shape of $\rho\in\CP(r,k_n)$ exists under the usual
$\sqrt{r}$-scaling, so that for $x>0$ and $r\to\infty$ we have
approximately
$$
\frac{\rho_{x\sqrt{r}}}{\sqrt{r}}\approx \phi(x).
$$
By the relation \eqref{eq:rho} and the asymptotic formulas
\eqref{eq:assq*} and \eqref{eq:r:n}, this implies
\begin{align}
\notag
\frac{\lambda_{x\sqrt{n}}}{\sqrt{n}}=
\frac{\rho_{x\sqrt{n}}}{\sqrt{n}}-\frac{1}{\sqrt{n}}+\frac{Q_{k_n-x\sqrt{n}}}{\sqrt{n}}&\approx
b\,\phi(x/b)+\frac{q\left(k_n-x\sqrt{n}\myp\right)}{\sqrt{n}}\\
\label{eq:LS-heuristic} &\approx b\,\phi(x/b)+q\mypp(\tau-x),
\end{align}
which yields the limit shape for $\lambda\in\CP_\bq(n,k_n)$ as
$n\to\infty$. Note that the last term in \eqref{eq:LS-heuristic}
corresponds to the ground state discussed earlier, whereas the first
term indicates the contribution from the ``free part'' of the
partition $\lambda\in\CP_\bq(n,k_n)$.

Likewise, for partitions $\lambda\in\CP(n)$, assuming that
their length follows the typical behaviour $K(\lambda)\approx
T_q\myp\vartheta_q^{-1}\sqrt{n}$ (see Corollary~\ref{cor:LLN-Kn}),
formula \eqref{eq:LS-heuristic} yields the limit shape
\begin{equation*}
\frac{\lambda_{x\sqrt{n}}}{\sqrt{n}}\approx
b_q\,\phi(x/b_q)+q\left(\frac{T_q}{\vartheta_q}-x\right),
\end{equation*}
where $b_q:=\sqrt{1-\smash{\tfrac12\myp q\myp
T_q^2/\vartheta_q^{2}}\vphantom{T^|}}$ \,(cf.~\eqref{eq:bdef}).

Let us now give a more rigorous argumentation. We confine ourselves
to the case $q>0$ and prove a weaker statement than in the previous
section (i.e., just convergence in probability instead of
exponential bounds on deviations), since known results can be
applied in this case. A similar approach was used by
Romik~\cite{Romik} to find the limit shape of $\MDP(q)$ with $q=2$,
and by DeSalvo \& Pak \cite{Pak2} for any positive integer $q$. The
same technique can be worked out in the case $q=0$, but this
requires a more detailed analysis.

The limit shape for partitions under the uniform measure
$\nu_{r,\myp k}$ on the space $\CP(r,k)$ has been found by Vershik
\& Yakubovich \cite{VY-MMJ} (see also Vershik~\cite{Vershik96}).
Adapted to our notation, this result is formulated as follows.
Recall that a partition $\rho'$ is said to be \emph{conjugate} to
partition $\rho\in \CP(r)$ if their Young diagrams
$\varUpsilon_\rho$ and $\varUpsilon_{\rho'}$ are symmetric to one
another with respect to reflection about the main diagonal of the
coordinate plane. In other words, column blocks of the diagram
$\varUpsilon_\rho$ become row blocks of the diagram
$\varUpsilon_{\rho'}$, and vice versa. Clearly, $\rho'$ has the same
weight as $\rho$, that is, $\rho'\in \CP(r)$. The next result refers
to the conjugate Young diagrams $\varUpsilon_{\rho'}$, but it easily
translates to the original diagrams $\varUpsilon_{\rho}$\myp.

\begin{theorem}[\protect{\cite[Theorem~1]{VY-MMJ}}]
\label{th:5.1}
Let $r,k\to\infty$ so that $k=c\mypp\sqrt{r}+O(1)$
with some $c>0$, then for any $\varepsilon>0$
\begin{equation}\label{eq:lsnucirc}
\nu_{r,\myp
k}\mynn\left\{\rho\in\CP(r,k)\colon\sup\nolimits_{u\ge0}|k^{-1}Y_{\rho'}(ru/k)
-\psi_c(u)|>\varepsilon\right\}\to
0,
\end{equation}
where\footnote{There is a misprint in~\cite[Eq.\:(5),
p.\,459]{VY-MMJ}, where the variable $u$ should be replaced with
$-uc^{-2}\log\myn(1-y_c)$.}
\begin{equation}\label{eq:defphicirc}
\psi_c(u):=\frac{\log\bigl(1-y_c\mypp(1-y_c)^{u/c^{2}}\bigr)}{\log\myn(1-y_c)},\qquad
u\ge0,
\end{equation}
and $y_c\in(0,1)$ is the \textup{(}unique\textup{)} solution of the
equation
\begin{equation}\label{eq:defyc}
c^2\Li_2(y_c)=\log^2(1-y_c).
\end{equation}
\end{theorem}
Equivalently, the statement of Theorem \ref{th:5.1} can be rewritten
as follows: \emph{for any $s_0\in(0,1]$ and $\varepsilon>0$,}
\begin{equation}\label{eq:lsalt:t1}
\nu_{r,\myp k}\bigl\{\rho\in\CP(r,k)\colon \sup\nolimits_{s\in
[s_0,1]}|k\myp r^{-1}\myp Y_\rho(sk)-\phi_c(s)|
>\varepsilon\bigr\}\to 0,
\end{equation}
\emph{where $\phi_c(s)$ is the inverse function,}
\begin{equation}\label{eq:phi}
\phi_c(s):=\psi_c^{-1}(s)=\frac{c^2}{\log\myn(1-y_c)}\,\log\left(\frac{1-(1-y_c)^{s}
}{y_c}\right),\qquad s\in(0,1].
\end{equation}

%Since $\psi_c(\cdot)$ is a decreasing function and its derivative
%is bounded away from zero on any interval $[0,u_0]$, the scaled
%Young diagram of the original partition~$\rho$ is close to the
%inverse function.
Note that the scalings used here along the two axes are both
proportional to $\sqrt{r}$ but different (unless $c=1$).
Unfortunately, the condition $k=c\mypp\sqrt{r}+O(1)$ is too strong
for our purposes. However, tracking the proof given in~\cite{VY-MMJ}
and using the continuity of the expression \eqref{eq:defphicirc}
with respect to $c$, one can verify that the limits
\eqref{eq:lsnucirc} and \eqref{eq:lsalt:t1} hold true provided only
that $k\sim c\mypp\sqrt{r}$.

Returning to the limit shape problem for partitions
$\lambda\in\CP(n,k_n)$, with $k_n\sim\tau\sqrt{n}$, put
\begin{equation}\label{eq:c}
c=\frac{\tau}{b}=\frac{\tau}{\sqrt{1-\vphantom{T^I}
\smash{\tfrac12}\myp q\myp \tau^2}},
\end{equation}
so that $k_n\sim \tau\sqrt{n} \sim c\sqrt{r}$
\,(see~\eqref{eq:r:n}). Let $T_*$ be the solution of the
equation~\eqref{eq:TofB}. Using the definition \eqref{eq:deftheta},
it is straightforward to check that $y_c=1-\re^{-T_*}$ solves the
equation~\eqref{eq:defyc}. Furthermore, expressing $\tau$ from
\eqref{eq:TofB} and using \eqref{eq:deftheta}, formula \eqref{eq:c}
can be rewritten as
$$
c^2=\frac{T_*^2}{\vartheta_q^2(T_*)-\tfrac12\myp q\myp
T_*^2}=\frac{T_*^2}{\Li_2(1-\re^{-T_*})}.
$$
Hence, the expression \eqref{eq:phi} takes the form
\[
\phi_c(s)=\frac{T_*}{\Li_2(1-\re^{-T_*})}\log\frac{1-\re^{-T_*}}{1-\re^{-sT_*}},
\qquad s\in(0,1],
\]
and the asymptotic result \eqref{eq:lsalt:t1}, restated in the new
variable $t=sT^*/\vartheta_q(T_*)$, readily yields
\begin{equation}\label{eq:LS-rho-lambda}
\lim_{n\to\infty}\nu_{n,\myp k_n}^\bq\!\left\{\sup_{t\in[t_0,
T_*/\vartheta_q(T_*)]}
\left|\frac{1}{\sqrt{n}}\,Y_{\mathfrak{I}(\lambda)}(t\sqrt{n}\mypp)
-\frac{1}{\vartheta_q(T_*)}\log\frac{1-\re^{-T_*}}{1-\re^{-t\myp\vartheta_q(T_*)}}\right|>\varepsilon\right\}=0.
\end{equation}
Finally, to see how \eqref{eq:LS-rho-lambda} produces the expression
for the limit shape $\varphi_{T_*}(t;q)$ already obtained in
Theorem~\ref{th:nunmls}, it remains to notice, using
\eqref{eq:assq*} and \eqref{eq:rho}, that
(cf.~\eqref{eq:LS-heuristic})
\[
\frac{Y_{\lambda}(t\sqrt{n}\myp)-Y_{\mathfrak{I}(\lambda)}(t\sqrt{n}\myp)}{\sqrt{n}}
=\frac{Q_{k_n-\lfloor t\sqrt{n}\rfloor}}{\sqrt{n}} \to
q\left(\frac{T_*}{\vartheta_q(T_*)}-t\right),
\]
for all $\lambda\in\CP_\bq(n,k_n)$ and uniformly in
$t\in[0,T_*/\vartheta_q(T_*)]$.

In a similar fashion, one can prove Theorem~\ref{th:nu'nls}. More
specifically, by Corollary \ref{cor:LLN-Kn} $K(\lambda)/k_n\to1$ in
$\nu^\bq_n$-probability, where $k_n=(T_q/\vartheta_q)\myp\sqrt{n}$.
The push-forward $\mathfrak{I}^*\myn\nu^\bq_n=\nu^\bq_n\circ
\mathfrak{I}^{-1}$ under the bijection $\mathfrak{I}$ defined in
\eqref{eq:map} is a measure on partitions $\rho\in\CP$ such that
(random) $r=N(\rho)$ and $k=K(\rho)$ satisfy the relation
$r=n+k-s_k$. Since $K(\lambda)=K(\rho)$, it follows that
$K(\rho)/k_n\to1$ in $(\mathfrak{I}^*\myn\nu^\bq_n)$-probability.
Hence, using \eqref{eq:theta-intro-alt} and \eqref{eq:sksim},  we
obtain, in $(\mathfrak{I}^*\myn\nu^\bq_n)$-probability as
$n\to\infty$,
\[
\frac{r}{k^2}=\frac{n}{k^2}+\frac{1}{k}-\frac{s_k}{k^2}\to\frac{\vartheta_q^2}{T_q^2}-\frac{q}{2}=
\frac{\Li_2(1-\re^{-T_q})}{T_q^2}>0.
\]
Thus, taking $c=T_q/\sqrt{\myn\Li_2(1-\re^{-T_q})}$ it is easy to
see that $y_c=1-\re^{-T_q}$ solves the equation~\eqref{eq:defyc}.
Furthermore, using \eqref{eq:deftd} the expression \eqref{eq:phi} is
reduced to
\[
\varphi(t)=\frac{T_q}{\Li_2(1-\re^{-T_q})}\left(-qT_q-\log\myn(1-\re^{-t\myp
T_q})\right),\qquad t\in(0,1],
\]
and \eqref{eq:lsalt:t1} implies that
\begin{equation}\label{eq:I(lambda)LS}
\lim_{n\to\infty}\nu_n^\bq\left\{\sup_{t\in[t_0, T_q/\vartheta_q]}
\left|\frac{1}{\sqrt{n}}\,Y_{\mathfrak{I}(\lambda)}(t\sqrt{n}\mypp)
-\frac{-qT_q-\log\myn(1-\re^{-t\vartheta_q})}
{\vartheta_q}\right|>\varepsilon\right\}= 0.
\end{equation}
It remains to notice, using condition \eqref{eq:assq*}, that in
$\nu_n^\bq$-probability
\begin{equation*}
\sup_{t\in[t_0,
T_q/\vartheta_q]}\left|\frac{Y_{\lambda}(t\sqrt{n}\myp)-Y_{\mathfrak{I}(\lambda)}(t\sqrt{n}\myp)}{\sqrt{n}}
-q\left(\frac{T_q}{\vartheta_q}-t\right)\right|\to 0,
\end{equation*}
which, together with \eqref{eq:I(lambda)LS}, yields the expression
$\varphi_{T_q}(t;q)$ for the limit shape already obtained in
Theorem~\ref{th:nu'nls}.

\section{Minimal difference partitions with random gaps}\label{sec:MDPRE}

The basic assumption \eqref{eq:assq*}, that the partial sums $Q_k$
of the gap sequence $\bq=(q_i)$ asymptotically grow linearly ($q>0$)
or sub-linearly ($q=0$), may be satisfied not only for fixed
sequences but also for those obtained via some stochastic procedure.
Without attempting to investigate this issue in full generality, we
provide sufficient conditions for the asymptotics \eqref{eq:assq*}
under two simple models for random gaps:
\begin{itemize}
\makeatletter
\renewcommand{\@listI}%
{\leftmargin=\parindent
%\leftmargini
\partopsep=0pc
\parskip=-.5pc
\topsep=-1pc
\itemsep=30pc
\labelwidth=\leftmargini }
\renewcommand{\@listii}{\setlength{\topsep}{-.8pc}
} \setlength{\itemsep}{.2pc} \setlength{\partopsep}{-.5pc}
\setlength{\parskip}{-.0pc}
\makeatother

\vspace{-.4pc}
\item[(i)] $\bq=(q_i)$ is a \emph{sequence of independent
random variables};
\item[(ii)] $\bq=(q_i)$ is
generated using a \emph{random walk in random environment (RWRE)},
that is, a (nearest-neighbour) random  walk with random transition
probabilities.

\vspace{-.4pc}
\end{itemize}

In what follows, abbreviation ``a.s.'' stands for ``almost surely''
with respect to the suitable probability measure (i.e., law of the
sequence $\bq$).

\subsection{Random gaps modelled as an independent sequence}
Suppose that $\bq=(q_i)$ is a sequence of independent (not
necessarily identically distributed) random variables (such that
$q_i\ge0$, $q_0\ge1$), defined on an auxiliary probability space
with probability measure $\PP$\myp; we denote by $\EE$ the
corresponding expectation.

We will need the following standard result about the strong law of
large numbers for independent sequences.
\begin{lemma}[\mbox{\cite[Theorem~6.6, p.\,209]{Petrov}}]\label{lm:Petrov}
Let $(X_{i})_{i\in\NN}$ be a sequence of independent random
variables, and let constants $a_i>0$ be such that
$a_i\uparrow\infty$. If, for some $p\in[0,1]$,
\begin{equation}\label{eq:sumXj}
\sum_{i=1}^\infty \frac{\EE(|X_i|^{p})}{a_i^{p}}<\infty,
\end{equation}
or if\/ $\EE(X_i)=0$ for all $i\in\NN$ and \eqref{eq:sumXj} holds
for some $p\in(1,2]$, then
\begin{equation}\label{eq:SLLN}
\frac{X_1+\dots+X_{k}}{a_k}\to 0\qquad (\text{$\PP$-a.s.}).
\end{equation}
\end{lemma}

\begin{theorem}\label{th:6.2}
Suppose that, for some $p\in(0,1]$ and $\delta\in[0,p)$,
\begin{equation}\label{eq:moments'}
\sum_{i=0}^{k-1}\EE(q_i^{p})=O(k^\delta)\qquad (k\to\infty).
\end{equation}
Then the asymptotic relation \eqref{eq:assq*} holds $\PP$-a.s.\ with
$q=0$ and any $\beta\in(\delta/p,1)$.
\end{theorem}
\proof We wish to apply the first part of Lemma~\ref{lm:Petrov} with
$X_i=q_{i-1}$ and $a_i=i^{\beta}$ ($\beta>\delta/p$). Denoting
$S_k^{(p)}\!:=\sum_{i=1}^{k} \EE(X_i^p)$ \,($S_0^{(p)}:=0$) and
using summation by parts, we obtain
\begin{align}
\notag \sum_{i=1}^{k} \frac{\EE(X_{i}^{p})}{a_i^{p}}
&=\sum_{i=1}^k\frac{S_i^{(p)}-S_{i-1}^{(p)}}{i^{\beta p}}\\
&=\frac{S_k^{(p)}}{k^{\beta p}}+\sum_{i=1}^{k-1}
\left(\frac{1}{i^{\beta p}}-\frac{1}{(i+1)^{\beta
p}}\right)S_i^{(p)}. \label{eq:6.3}
\end{align}
Furthermore, note that
\begin{align*}
\frac{1}{i^{\beta p}}-\frac{1}{(i+1)^{\beta p}}&=\frac{1}{i^{\beta
p}}\left(1-\left(1+\frac{1}{i}\right)^{-\beta
p}\right)\le\frac{\beta p}{i^{\beta p+1}},
\end{align*}
by the elementary inequality $(1+x)^{-\gamma}\ge 1 -\gamma\myp x$
(see \cite[Theorem~41, Eq.\:(2.15.1), p.\,39]{Hardy}) with $x=1/i$
and $\gamma=\beta p$. Hence, returning to \eqref{eq:6.3}, we get
\begin{align}
\sum_{i=1}^{k} \frac{\EE(X_{i}^{p})}{a_i^{p}} \le
\frac{S_k^{(p)}}{k^{\beta p}}+ \beta p
\sum_{i=1}^{k-1}\frac{S_i^{(p)}}{i^{\beta p+1}}. \label{eq:6.4}
\end{align} From the hypothesis \eqref{eq:moments'},
we know that $S_k^{(p)}=O(k^\delta)$, and together with the
assumption $\beta>\delta/p$ this implies that the right-hand side of
\eqref{eq:6.4} stays bounded as $k\to\infty$. Thus, the condition
\eqref{eq:sumXj} is satisfied, and \eqref{eq:assq*} follows due
to~\eqref{eq:SLLN}.
\endproof

Similarly, we can treat the case where the random variables have
finite expected values.
\begin{theorem}\label{th:6.3}
Suppose that the following two conditions are satisfied.
\begin{itemize}
\item[\rm (i)]
For some $q\ge 0$ and $\beta_0\in [0,1)$,
\begin{equation}\label{eq:means}
\sum_{i=0}^{k-1}\EE(q_i)=q\myp k+O(k^{\beta_0})\qquad (k\to\infty).
\end{equation}
\item[\rm (ii)]
For some $p\in(1,2]$ and $\delta\in[0,p)$,
\begin{equation}\label{eq:moments}
\sum_{i=0}^{k-1}\EE\!\left(|q_i-\EE(q_i)|^{p}\right)=O(k^\delta)\qquad
(k\to\infty).
\end{equation}
\end{itemize}
Then the asymptotic relation \eqref{eq:assq*} holds\/ $\PP$-a.s.\
with $q\ge0$ defined in \eqref{eq:means} and $\beta=\beta_0$ if
$\beta_0>\delta/p$, or else with any $\beta\in (\delta/p,1)$.
\end{theorem}
\proof
We can use the second part of Lemma~\ref{lm:Petrov} with
$X_i=q_{i-1}-\EE(q_{i-1})$ and $a_i=i^{\beta}$ ($\beta>\delta/p$).
Indeed, repeating the argumentation in the proof of
Theorem~\ref{th:6.2} and using the assumption \eqref{eq:moments}, we
see that \eqref{eq:SLLN} holds, that is, $\PP$-a.s.
$$
Q_k-\sum_{i=0}^{k-1}\EE(q_i)=o(k^\beta)\qquad (k\to\infty).
$$
Furthermore, on account of the assumption \eqref{eq:means} this
yields
\begin{equation}\label{eq:6.8}
Q_k=qk+O(k^{\beta_0}) + o(k^\beta)\qquad (k\to\infty).
\end{equation}
It remains to notice that if $\beta_0\le \delta/p$ then the combined
error term on the right-hand side of \eqref{eq:6.8} is
$o(k^{\beta})$ (with any $\beta>\delta/p$), while if
$\beta_0>\delta/p$ then this error term is $O(k^{\beta})$ with
$\beta=\beta_0$. This completes the proof of Theorem~\ref{th:6.3}.
\endproof

\begin{example}
To illustrate Theorem \ref{th:6.2}, let $q_i$ have a Bernoulli
distribution,
$$
\PP\bigl(q_i=(i+1)^3\bigr)=(i+1)^{-2},\qquad
\PP(q_i=0)=1-(i+1)^{-2}\qquad (i\in\NN_0). $$ Then for $p\in(0,1]$
we have
$$
\sum_{i=0}^{k-1}\EE(q^p)=\sum_{i=1}^{k} (i+1)^{3p-2}=\begin{cases}
O(k^{3p-1}),&p\in(\frac13,1],\\
O(\log k),& p=\frac13,\\
O(1),& p\in(0,\frac13).
\end{cases}
$$
Thus, the assumption \eqref{eq:moments'} holds with $\delta=3p-1$ if
$p>\frac13$; any $\delta>0$ if $p=\frac13$; and $\delta=0$ if
$p<\frac13$. Hence, the condition $\delta<p$ (as required in
\eqref{eq:moments'}) is satisfied for all $p\in(0,\frac12)$, and
therefore Theorem \ref{th:6.2} is applicable. In contrast, Theorem
\ref{th:6.3} cannot be used, because
$$
\sum_{i=0}^{k-1}\EE(q_i)=\sum_{i=0}^{k-1} (i+1)\sim\tfrac12\myp
k^2\qquad (k\to\infty),
$$
so that the condition \eqref{eq:means} is not fulfilled.
\end{example}

\begin{example}
Consider the particular case where the (independent) random
variables $(q_i)_{i\ge 1}$ are \emph{identically distributed}, and
suppose that, for some $p\in(0,2]$,
\begin{equation}\label{eq:finite}
\EE(q_1^{p})<\infty.
\end{equation}
If $p\le 1$ then the condition \eqref{eq:moments'} is satisfied only
with $\delta\ge1$, unless $\EE(q_1^p)=0$, that is, $q_1=0$
($\PP$-a.s.) when $\delta=0$. Hence, in a non-degenerate case, we
always have $\delta/p\ge1$ and Theorem~\ref{th:6.2} cannot be used.
However, the situation becomes more meaningful if $1<p\le2$. Here,
the conditions \eqref{eq:means} and \eqref{eq:moments} are satisfied
with $q=\EE(q_1)\ge0$, $\beta_0=0$ and $\delta\ge 1$ (assuming that
$\PP(q_1>0)>0$). Hence, by Theorem~\ref{th:6.3}, the asymptotic
relation \eqref{eq:assq*} holds with any $\beta\in(1/p,1)$. Note
that no moment assumption is required on $q_0$, because $q_0/k\to0$
($\PP$-a.s.). If $p=2$ (i.e., $q_1$ has finite variance), then the
law of the iterated logarithm shows that one cannot take
$\beta=\frac12$; in the general case $p\in(1,2)$, the optimality of
the lower bound $\beta>1/p$ follows from \strut{}\cite[\S\myp7.5.16,
p.\,258]{Petrov}.
\end{example}

\subsection{Random gaps modelled via RWRE}
RWRE is a random process $(X_k)_{k\ge0}$ on $\ZZ$ constructed in two
steps: (i) first, the environment $\omega\in\varOmega$ is chosen at
random (under some probability measure $\PP$) and fixed; (ii)
conditional on $\omega$, $(X_k)$ is a time-homogeneous random walk
(Markov chain) with state-dependent transition probabilities
determined by the environment. More precisely, let
$p_j=p_j(\omega)\in(0,1)$ ($j\in\ZZ$) be a family of
\emph{independent and identically distributed} random variables,
defined on a sample space $\varOmega=\{\omega\}$. Denoting by
$P^{\myp\omega}_0$ the \emph{quenched} probability law of the random
walk $(X_k)$ conditioned on the environment $\omega\in\varOmega$
(where the subscript $0$ indicates the starting position of the
walk, $X_0=0$), we have, for all $k\in\NN$ and $j\in\ZZ$,
$$
P^{\myp\omega}_0(X_k=j+1\mypp|\mypp X_{k-1}=j)=p_j(\omega),\qquad
P^{\myp\omega}_0(X_k=j-1\mypp|\mypp X_{k-1}=j)=1-p_j(\omega).
$$
By averaging the quenched measure $P^{\myp\omega}_0$ with respect to
the environment distribution $\PP$, we obtain the \emph{annealed}
measure $\BP_{\myn0}:= \PP\times P^{\myp\omega}_0\equiv \EE\myp
P^{\myp\omega}_0$.
%skew product
For a general review of RWRE, with further details and references,
see, for example, Bogachev~\cite{BogachevRWRE} or
Zeitouni~\cite{Zeitouni}.

Now, given RWRE $(X_k)$, we can generate the gap sequence
$\bq=(q_i)$ as follows:
\begin{equation}\label{eq:Q_RWRE}
q_i=a_i+b\mypp(X_{i+1}-X_{i})\qquad (i\in\NN_0),
\end{equation}
where $b>0$, $a_0\ge b+1$, $a_i\ge b$ for $i\ge1$, and
$q_i\in\NN_0$. Hence, recalling that $X_0=0$, we get
\begin{equation}\label{eq:Q1}
Q_k=\sum_{i=0}^{k-1}q_i=A_k+bX_k,
\end{equation}
where $A_k:=\sum_{i=0}^{k-1} a_i$. To obtain asymptotics
\eqref{eq:assq*} for the sequence \eqref{eq:Q1}, it is natural to
assume that the leading sequence $(a_i)$ itself satisfies a similar
condition,
\begin{equation}\label{eq:Ak}
A_k=ak+O(k^{\beta_0})\qquad (k\to\infty),
\end{equation}
with some $a\ge b$ and $\beta_0\in[0,1)$. In turn, the long-time
behaviour of the RWRE $(X_k)$ is described by the following results
due to Solomon~\cite{Solomon} (for a quick orientation, see also
\cite[Theorems 1 and 2, pp.\,355--356]{BogachevRWRE}).
\begin{lemma}[\mbox{\cite[Theorem~(1.7), p.\,4]{Solomon}}]\label{lm:Solomon1}
Set $\rho_0:=(1-p_0)/p_0$ and $\eta:=\EE\myp(\log \rho_0)$.

\begin{itemize}
\item[\rm (a)] If\/ $\eta<0$ then
$\lim_{k\to\infty}X_k=+\infty$, while if\/ $\eta>0$ then\/
$\lim_{k\to\infty}X_k=-\infty$
\textup{(}$\BP_{\myn0}$-a.s.\textup{)}.

\item[\rm (b)] If\/ $\eta=0$ then $-\infty=\liminf_{k\to\infty} X_k<\limsup_{k\to\infty}
X_k=+\infty$ \textup{(}$\BP_{\myn0}$-a.s.\textup{)}.
\end{itemize}
\end{lemma}

Note that, by Jensen's inequality, $\EE\myp(\log \rho_0)\le \log
\EE\myp(\rho_0)$ and $\{\EE\myp(\rho_0)\}^{-1}\myn\le
\EE\myp(\rho_0^{-1})$, with all inequalities strict unless $\rho_0$
is a deterministic constant.
\begin{lemma}[\mbox{\cite[Theorem~(1.16),
p.\,7]{Solomon}}]\label{lm:Solomon2} The limit
\,$v:=\lim_{k\to\infty} X_k/k$ \,exists $\BP_{\myn0}$-a.s.\ and is
given by
\begin{equation}\label{drift}
 v=  \begin{cases}
  \displaystyle
  \hphantom{-}\frac{1-\EE\myp(\rho_0)}{1+\EE\myp(\rho_0)}&\text{if}\,\
  \EE\myp(\rho_0)<1,\\[.7pc]
  \displaystyle-\frac{1-\EE\myp(\rho_0^{-1})}{1+\EE\myp(\rho_0^{-1})}&\text{if}\,\
  \EE\myp(\rho_0^{-1})<1,\\[.7pc]
  \hphantom{-}\mypp0&\text{if}\,\ \{\EE\myp(\rho_0)\}^{-1}\myn\le 1\le
  \EE\myp(\rho_0^{-1}).
\end{cases}
\end{equation}
\end{lemma}

\smallskip
\begin{remark}\label{rm:|v|}
Formula \eqref{drift} implies that $|v|<1$ in all cases.
\end{remark}

From Lemma \ref{lm:Solomon2} and the condition \eqref{eq:Ak}, we
immediately deduce a strong law of large numbers for the sequence
$Q_k$ (see \eqref{eq:Q1}),
$$
\frac{Q_k}{k}=\frac{A_k}{k}+\frac{bX_k}{k}\to a+b\myp v, \qquad
k\to\infty\qquad (\text{$\BP_{\myn0}$-a.s.}),
$$
so that the limit \eqref{eq:q*} holds $\BP_{\myn0}$-a.s.\ with
$q=a+b\myp v$. \begin{remark} By the inequality $a\ge b$ and
Remark~\ref{rm:|v|}, in the model \eqref{eq:Q_RWRE} we always have
$q>a-b\ge0$.
\end{remark}

To estimate the error term in a way similar to \eqref{eq:assq*}, we
need information about the fluctuations of the RWRE $(X_k)$ as
$k\to\infty$. In the non-critical case (i.e., $\eta\ne0$, see
Lemma~\ref{lm:Solomon1}), this was investigated by Kesten, Kozlov \&
Spitzer \cite{KKS} (see also discussion and commentary in
\cite[pp.\,357--359]{BogachevRWRE}). We will state below (a
corollary from) their results adapted to our purposes. A probability
law on $\RR$ is called \emph{non-arithmetic} if it is not supported
on a set $c\mypp\ZZ$. We write $Y_k=O_p(1)$ if $(Y_k)$ is
\emph{stochastically bounded} (in $\BP_{\myn0}$), that is, if for
any $\varepsilon>0$ there is $M>0$ such that
 $\limsup_{k\to\infty}\BP_{\myn0}(|Y_k|>M)\le \varepsilon$. The
results from \cite{KKS} are transcribed using that if $Y_k$ weakly
converges (to a proper distribution) then $Y_k=O_p(1)$. Recall the
notation $\rho_0=(1-p_0)/p_0$ and $\eta=\EE\myp(\log \rho_0)$.

\begin{lemma}[\mbox{\cite[pp.\mypp146--148]{KKS}}]\label{lm:KKS} Assume that $-\infty\le
\eta<0$ and the distribution of\/ $\log\rho_0$ \textup{(}excluding a
possible atom at $-\infty$\textup{)} is non-arithmetic. Let
$\varkappa\in(0,\infty)$ be such that
$$
\EE(\rho_0^\varkappa)=1\quad\text{and}\quad
\EE(\rho_0^\varkappa\log^+\!\mynn\rho_0 )<\infty,
$$
where $\log^+\!\myn u:= \max\{\log u, 0\}$. Then RWRE $(X_k)$ has
the following asymptotics as $k\to\infty$.
\begin{itemize}
\item[\rm(a)] If\/ $0<\varkappa<1$ then
$$
X_k=O_p(k^\varkappa).
$$
\item[\rm(b)] If\/ $\varkappa=1$ then
$$
X_k=O_p\bigl(k/\log k\bigr).
$$

\item[\rm(c)] If\/ $\varkappa>1$ then
$$
X_k=vk+O_p(k^{\beta_1}),
$$
where $v$ is defined in~\eqref{drift} and
$\beta_1:=\max\{1/2,1/\varkappa\}$.
\end{itemize}
\end{lemma}
Combining Lemma \ref{lm:KKS} and the assumption \eqref{eq:Ak}, we
arrive at the following result. Recall that $\beta_0\in(0,1)$ is
defined in \eqref{eq:Ak}.
\begin{theorem}\label{th:RWRE-Qk}
Under the hypotheses of Lemma~\textup{\ref{lm:KKS}}, the following
asymptotics hold for $Q_k$ as $k\to\infty$.
\begin{itemize}
\item[\rm(a)] If\/ $0<\varkappa<1$ then
$$
Q_k=ak+O_p(k^\beta),
$$
where $\beta=\max\{\beta_0,\varkappa\}$.
\item[\rm(b)] If\/ $\varkappa=1$ then
$$
Q_k=ak+O_p\bigl(k/\log k\bigr).
$$

\item[\rm(c)] If\/ $\varkappa>1$ then
$$
Q_k=(a+v)\myp k+O_p(k^{\beta}),
$$
where $\beta:=\max\{\beta_0,1/2,1/\varkappa\}$.
\end{itemize}
\end{theorem}

Thus, in the RWRE model \eqref{eq:Q_RWRE} the asymptotic formula
\eqref{eq:assq*} is valid in a $\BP_{\myn0}$-stochastic version,
that is, with the error term estimated using $O_p(\cdot)$,
\begin{equation}\label{eq:assq*RWRE}
Q_k=q\myp k+O_p(k^{\beta}) \qquad (k\to\infty),
\end{equation}
where $q=a+v>0$ and $0\le\beta<1$. To be more precise, formula
\eqref{eq:assq*RWRE} with $\beta<1$ holds in all cases except for
$\varkappa=1$, where the error bound becomes logarithmically close
to $k$.

Furthermore, a careful inspection of all the proofs shows that a
stochastic version \eqref{eq:assq*RWRE} of the asymptotics
\eqref{eq:assq*} is sufficient to guarantee convergence (in
$\BP_{\myn0}$-probability) of the scaled Young boundary
$Y_\lambda(t)$ to the limit shape, as described in
Section~\ref{sec:pdls}. As for the special case $\varkappa=1$, it is
natural to expect that the error bound of order $k/\log k$ should be
enough for the limit shape, but verification of the technical
details is tedious, so this is left as a conjecture.

Finally, let us mention the critical case $\eta=0$ not covered by
Lemmas \ref{lm:Solomon1} and~\ref{lm:Solomon2}. Here, RWRE $(X_k)$
is recurrent (see Lemma~\ref{lm:Solomon1}\myp(b)), and its
asymptotic behaviour is characterized by the so-called \emph{Sinai's
localization} \cite{Sinai} (see discussion and commentary in
\cite[pp.\,359--360]{BogachevRWRE}). We state a corollary from this
result adapted to our purposes.
\begin{lemma}[\cite{Sinai}]\label{lm:Sinai}
Suppose that\/ $\PP(\rho_0=1)<1$ and $c_1\le \rho_0\le c_2$
\textup{(}$\PP$-a.s.\textup{)}, with some deterministic constants
$0<c_1<c_2<\infty$. If\/ $\eta=0$ then
$$
X_k=O_p(\log^2k)\qquad (k\to\infty).
$$
\end{lemma}

Combined with \eqref{eq:Ak}, this immediately implies asymptotics of
$Q_k$ (cf.~\eqref{eq:assq*}), which ensures the validity of our
limit shape result.
\begin{theorem}\label{th:RWRE-Qk-cr}
Under the hypotheses of\ Lemma~\textup{\ref{lm:Sinai}},
$$
Q_k=ak+O_p(k^{\beta_0})\qquad (k\to\infty),
$$
where $\beta_0\in(0,1)$ is defined in \eqref{eq:Ak}.
\end{theorem}

\appendix
\section{Appendix: Proof of Propositions~\ref{pr:munu}
and~\ref{pr:mumnum}}\label{sec:prop12pf}

\subsection{Auxiliary lemmas}
According to the representation \eqref{eq:NK} and independence of
$\{D_j\}$ under the measure $\mu^\bq_{z,\myp k}$ (see
Section~\ref{sec:defs}), the weight $N(\lambda)$ of partition
$\lambda\in\CP_\bq$ is the sum of $k\to\infty$ independent random
variables, so one may expect a local limit theorem to hold (cf.\
\cite{Bogachev,FVY,Fristedt,GSE}). For our purposes, it suffices to
obtain an asymptotic lower bound for the probability of the event
$\{N(\lambda)=n\}$. To this end, we need some auxiliary technical
results (for simplicity, we suppress the dependence on $z$ in the
notation of some functions introduced below).

\begin{lemma}\label{lm:log-chi}
Let $\chi_{j}(u):=\BE^{\bq}_{z,\myp k}\bigl[\re^{\myp\ii\myp
ujD_j}\bigr]$ \textup{($u\in\RR$)} be the characteristic function of
the random variable\/ $jD_j$ \textup{(${1\le j\le k}$)}. Then, as
$z\downarrow0$, uniformly in $j\in\NN$ and $u\in\RR$
\begin{equation}\label{eq:chi-lim}
\log\chi_{j}(u)=\ii\myp \bigl(q_{k-j}+h(zj)\bigr)\myp ju
-\tfrac12\mypp h(zj)\bigl(1+h(zj)\bigr)\mypp j^2u^2 +R_{j}(u),
\end{equation}
where $\log\mynn(\cdot)$ denotes the principal branch of the
logarithm, $h(\cdot)$ is given by \eqref{eq:Psi} and
\begin{equation}\label{eq:chi-rem}
R_{j}(u)=\bigl(h(zj)+h(zj)^2+h(zj)^4\log\tfrac{1}{z}\bigr)\,
O(j^3u^3).
\end{equation}
\end{lemma}

\proof An easy computation shows that
\[
\chi_{j}(u)= \sum_{r=0}^\infty \re^{\ii
uj(r+q_{k-j})}\re^{-zjr}(1-\re^{-zj}) =\re^{\ii
ujq_{k-j}}\frac{1-\re^{-zj}}{1-\re^{-zj+\ii uj}}.
\]
Hence
\begin{align*}
\log\chi_{j}(u)&=\ii ujq_{k-j}-\log
\frac{1-\re^{-zj}-\re^{-zj}(\re^{\ii uj}-1)}{1-\re^{-zj}}\\
&=\ii ujq_{k-j}-\log\bigl(1-h(zj)(\re^{\ii uj}-1)\bigr).
\end{align*}
It is easy to check that the function
\begin{equation}\label{eq:gj}
\zeta\mapsto g_j(\zeta):=-\log\bigl\{1-h(zj)\myp(\zeta-1)\bigr\}
\end{equation}
is analytic in the half-plane $\Re\mypp\zeta<1+1/h(zj)$. Hence,
Taylor's formula for complex-analytic functions (see, e.g.,
\cite[\S\mypp5.2, p.\,244]{SaSn}) gives for $|\zeta|<1+1/h(zj)$
\begin{equation}\label{eq:Taylor}
g_j(\zeta)=g_j(1)+g'_j(1)(\zeta-1)+\frac{g''_j(1)}{2}(\zeta-1)^2
+\frac{(\zeta-1)^3}{2\pi\ii}\oint_{\varGamma_j}
\frac{g_j(\xi)}{(\xi-1)^3(\xi-\zeta)}\,\rd{\xi},
\end{equation}
where $\varGamma_j$ is the circle of radius $1+1/(2h(zj))$ about the
origin, positively oriented.

Note from \eqref{eq:gj} that for $\xi\in\varGamma_j$ we have
$$
|\re^{-g_j(\xi)}|=\bigl|1+h(zj)-h(zj)\mypp \xi\bigr|\le \tfrac32
+2\mypp h(zj),\qquad |\myn\arg \re^{-g_j(\xi)}|\le\frac{\pi}{2}.
$$
Using that $|\myn\log\myn (r\myp\re^{\ii\myp \theta})|\le |\myn\log
r|+\pi/2$ \,($r>0$, $|\theta|\le \pi/2$), this yields
\begin{align}
\notag |g_j(\xi)|=\bigl|\log\myn (\re^{-g_j(\xi)})\bigr|
&\le\log\mynn\bigl(\tfrac{3}{2}+2\mypp h(zj)\bigr) +\frac{\pi}2\\
\notag
&\le\log\myn\bigl(1+h(z)\bigr)+\log 2+\frac{\pi}{2}\\
&\le \log \frac{z+1}{z}+\log 2+\frac{\pi}{2}, \label{eq:|log|}
\end{align}
by virtue of monotonicity of $h(\cdot)$ and the elementary bound
$$
1+h(z)=\frac{1}{1-\re^{-z}}\le \frac{z+1}{z}.
$$
Furthermore, for any $\xi\in\varGamma_j$ and $|\zeta|=1$ (in
particular, $\zeta=1$) we have $|\xi-\zeta|^{-1}\le 2\mypp h(zj)$.
Thus, computing the derivatives of $g_j(\cdot)$ at $1$ and
substituting \eqref{eq:|log|} into \eqref{eq:Taylor}, we get
\begin{equation}\label{eq:Taylor1}
g_j(\zeta)=h(zj)(\zeta-1)+\frac{h(zj)^2}{2}\,(\zeta-1)^2
+(\zeta-1)^3\mypp h(zj)^4\mypp O(\log\tfrac{1}{z})\qquad
(z\downarrow 0),
\end{equation}
where the estimate $O(\cdot)$ is uniform in $j\in\NN$ and $\zeta$
such that $|\zeta|=1$.

Now, using the Taylor expansion
$$
\re^{\ii x}=\sum_{\ell=0}^{m-1}\frac{(\ii
x)^\ell}{\ell!}+R_m(x),\qquad |R_m(x)|\le\frac{|x|^m}{m!},
$$
which is valid for all $m\in\NN$ and any real $x$ (see, e.g.,
\cite[\S\myp{}XV.4, Lemma~1, p.\,512]{Feller}), we substitute
$\zeta=\re^{\ii ju}$ into \eqref{eq:Taylor1} to obtain, as
$z\downarrow0$,
\begin{multline*}
g_j(\re^{\ii ju}) =h(zj)\left(\ii
ju-\tfrac12\myp j^2u^2+O(j^3u^3)\right)-\tfrac{1}{2}\mypp h(zj)^2\mynn\left(j^2u^2+O(j^3u^3)\right)\\
+O(j^3u^3) \,h(zj)^4 \log\tfrac{1}{z}\myp ,
\end{multline*}
where all $O$-estimates are uniform in $j\in\NN$, $u\in\RR$. [Note
that it is convenient to use the representation
$(\zeta-1)^2=(\zeta^2-1)-2(\zeta-1)$.] Finally, rearranging the
terms we obtain \eqref{eq:chi-lim} and \eqref{eq:chi-rem}.
\endproof

\begin{lemma}\label{lm:sum_j}
For $r,\ell\in\NN$, denote
\begin{equation}\label{eq:Sigma}
\varSigma_{z,\myp k}(r,\ell):=\sum_{j=1}^k j^r\myp h(zj)^\ell.
\end{equation}
Then, uniformly in $k\ge t_1/z$ \textup{(}for any $t_1>0$\textup{)},
as $z\downarrow 0$,
\begin{gather}
\label{eq:sum-Psi-1} \varSigma_{z,\myp
k}(1,1)=z^{-2}\Li_2(1-\re^{-zk})+O(z^{-1}),\\[.5pc]
\label{eq:sum-Psi-2}
\varSigma_{z,\myp k}(2,2)>\tfrac12\mypp
z^{-3}(1-\re^{-2t_1}),\\[.5pc]
\label{eq:sum-Psi-3} \varSigma_{z,\myp k}(3,\ell)=O(z^{-4})\quad
(\ell=1,2,3),\qquad \varSigma_{z,\myp k}(3,4)=O(z^{-4}\log\tfrac1z).
\end{gather}
\end{lemma}
\proof Using the Euler--Maclaurin sum formula like in the proof of
Lemma \ref{lm:new1}, we obtain
\begin{align*}
\sum_{j=1}^k j\mypp
h(zj)&=\int_1^k\frac{x\mypp\re^{-zx}}{1-\re^{-zx}}\,\rd{x}
+O(1)\frac{\re^{-z}}{1-\re^{-z}} +O(1)\int_1^k
\frac{|(zx-1)\re^{-zx}+\re^{-2zx}|}{(1-\re^{-zx})^2}\,\rd{x}\\
&=z^{-2}
\int_0^{zk}\frac{y\,\re^{-y}}{1-\re^{-y}}\,\rd{y}+O(z^{-1})=z^{-2}
\Li_2(1-\re^{-zk})+O(z^{-1}),
\end{align*}
using the substitution $u=1-\re^{-y}$ and formula \eqref{eq:Dilog}.
Hence, \eqref{eq:sum-Psi-1} is proved.

Similarly, \eqref{eq:sum-Psi-2} follows from the asymptotic estimate
\[
\sum_{j=1}^k j^2\mypp h(zj)^2\sim z^{-3}
\int_0^{zk}\frac{y^2\,\re^{-2y}}{(1-\re^{-y})^2}\,\rd{y}>z^{-3}
\int_0^{t_1}\re^{-2y}\,\rd{y}= \tfrac12\mypp z^{-3}(1-\re^{-2t_1}).
\]
Finally, noting that $y\myp(1-\re^{-y})^{-1}\le \re^{y/2}$ for all
$y>0$, we obtain
\[
\sum_{j=1}^k j^3\mypp h(zj)^\ell\sim z^{-4}
\int_z^{zk}\frac{y^3\,\re^{-\ell
y}}{(1-\re^{-y})^\ell}\,\rd{y}<z^{-4}
\int_z^{\infty}y^{3-\ell}\mypp\re^{-\ell y/2}\,\rd{y},
\]
which is $O(z^{-4})$ for $\ell<4$ and $O(z^{-4}\log\tfrac1z)$ for
$\ell=4$, and \eqref{eq:sum-Psi-3} follows.
\endproof
\begin{remark}
Formula \eqref{eq:sum-Psi-1} may be obtained from
\eqref{eq:sum-log-Psi} by formal differentiation with respect to
$z$, using the dilogarithm identity \eqref{eq:Li2=}.
\end{remark}
\begin{lemma}\label{lm:muxnm}
Let\/ $v>\frac34$ \strut{}and $t_1>0$ be some constants. Then there
exists $\delta>0$ such that, for any $z\in(0,\delta)$ and
\strut{}all $k\ge t_1/z$, the inequality
\[
\mu^\bq_{z,\myp k}\{N(\lambda)=n\}\ge n^{-v}
\]
holds for all $n\in\NN$ satisfying the bound
\begin{equation}\label{eq:nbound}
\left|n-s_k-z^{-2}\Li_2(1-\re^{-zk})\right|\le z^{-4/3}.
\end{equation}
\end{lemma}

\proof Let us start by pointing out that, for $z$ sufficiently
small, the inequality \eqref{eq:nbound} has many integer solutions
$n$. Moreover, since $\Li_2(1-\re^{-t})$ increases in~$t$, it
follows from \eqref{eq:nbound} that for all $z>0$ small enough and
for every $k\ge t_1/z$,
\begin{equation}\label{eq:zsqrtn}
n\ge z^{-2}\Li_2(1-\re^{-zk})+s_k-z^{-4/3} \ge \tfrac12\myp
z^{-2}\Li_2(1-\re^{-t_1}\myn)>0.
\end{equation}

Now, using the decomposition $N(\lambda)=\sum_j j D_j(\lambda)$ and
independence of $D_j(\lambda)$ for different $j$ (see
Lemma~\ref{lm:D_j}), by the Fourier inversion formula we have
\begin{align}
\notag \mu^\bq_{z,\myp k}\{N(\lambda)=n\}
&{}=\frac{1}{2\pi}\int_{-\pi}^\pi \prod_{j=1}^k
\chi_{j}(s)\,\re^{-\ii sn}\,\rd{s} = \frac{1}{\pi}\int_{0}^\pi \Re
\prod_{j=1}^k
\chi_{j}(s)\,\re^{-\ii sn}\,\rd{s}\\
\notag &{}=\frac{1}{\pi}\int_{0}^{z^{7/5}} \Re \prod_{j=1}^k
\chi_{j}(s)\,\re^{-\ii sn}\,\rd{s}
+\frac{1}{\pi}\int_{z^{7/5}}^{\pi} \Re\prod_{j=1}^k
\chi_{j}(s)\,\re^{-\ii sn}\,\rd{s}\\
\label{eq:I1+I2}
&{}=:I_1+I_2.
\end{align}
First, we shall obtain a suitable lower bound for $I_1$ and then
show that $I_2$ is small.

Using Lemma \ref{lm:log-chi} and recalling the notation
\eqref{eq:Sigma}, we have
\begin{align}
\notag I_1&=\frac{1}{\pi}\int_{0}^{z^{7/5}} \Re\exp\left\{-\ii
un+\sum_{j=1}^{k}
\log\chi_{j}(u)\right\}\rd{u}\\
&\begin{aligned} =\frac{1}{\pi}\int_{0}^{z^{7/5}}
\Re\exp\Bigl\{-\ii& u\bigl(n -s_k-\varSigma_{z,\myp k}(1,1)\bigr)
-\tfrac12 u^2\bigl(\varSigma_{z,\myp k}(2,1)+\varSigma_{z,\myp k}(2,2)\bigr)\\
& +O(u^3)\bigl(\varSigma_{z,\myp k}(3,1)+\varSigma_{z,\myp
k}(3,2)+\varSigma_{z,\myp k}(3,4)\log\tfrac1z\bigr)\Bigr\}\,\rd{u}.
\end{aligned}
\label{eq:I1_1}
\end{align}
Due to the estimate \eqref{eq:sum-Psi-1} and the
assumption~\eqref{eq:nbound},
\begin{gather*}
n-s_k-\varSigma_{z,\myp k}(1,1)=O(z^{-4/3})\qquad (z\downarrow 0).
\end{gather*}
Next, using \eqref{eq:sum-Psi-2} we get
\begin{gather*}
\varSigma_{z,\myp k}(2,1)+\varSigma_{z,\myp
k}(2,2)>\varSigma_{z,\myp k}(2,2)\ge \tfrac12\myp
z^{-3}(1-\re^{-2t_1}\myn)\qquad (z\downarrow 0).
\end{gather*}
Finally, by virtue of \eqref{eq:sum-Psi-3}
\[
\varSigma_{z,\myp k}(3,1)+\varSigma_{z,\myp
k}(3,2)+\varSigma_{z,\myp
k}(3,4)\log\tfrac1z=O\bigl(z^{-4}(\log\tfrac1z)^2\bigr)\qquad
(z\downarrow 0).
\]
Substituting these three estimates into \eqref{eq:I1_1} and changing
the variable $u=z^{3/2}v$, we obtain, after some simple
calculations,
\begin{align}
\notag I_1&\ge\frac{z^{3/2}}{\pi}\int_{0}^{z^{-1/10}}
\Re\exp\Bigl\{-\ii v\mypp O(z^{1/6}) -\tfrac14\myp
v^2(1-\re^{-2t_1})
+O(v^3)\bigl(z^{1/2}(\log\tfrac1z)^2\bigr)\Bigr\}\,\rd{v}\\
\label{eq:I1} &\sim \frac{z^{3/2}}{\pi}\int_{0}^{\infty}
\exp\Bigl\{-\tfrac14\myp v^2(1-\re^{-2t_1})
\Bigr\}\,\rd{v}=\frac{z^{3/2}}{\sqrt{\pi(1-\re^{-2t_1})}}\qquad
(z\downarrow0).
\end{align}

Estimation of $I_2$ is based on the inequality
\begin{align*}
|\chi_{j}(s)|^2 &{} =\frac{(1-\re^{-zj})^2}{|1-\re^{-zj+\ii sj}|^2}
=1-\frac{|1-\re^{-zj+\ii sj}|^2-(1-\re^{-zj})^2}{|1-\re^{-zj+\ii
sj}|^2}
\\
&{}= 1-\frac{2\myp\re^{-zj}(1-\cos sj)}{|1-\re^{-zj+\ii sj}|^2} \le
1-\frac{2\myp\re^{-zj}(1-\cos sj)}{(1+\re^{-zj})^2} \le
1-\frac{\re^{-zj}(1-\cos sj)}{2}\,.
\end{align*}
This implies, for $k>k_1:=\floor{t_1/z}$ as in the statement of the
lemma, that
\begin{align}
\notag |I_2|\le \frac{1}{\pi}\int_{z^{7/5}}^{\pi}
\prod_{j=1}^{k}|\chi_{j}(s)|\,\rd{s}
&=\frac{1}{\pi}\int_{z^{7/5}}^{\pi}
\exp\biggl\{\frac{1}{2}\sum_{j=1}^{k}
\log|\chi_{j}(s)|^2\biggr\}\,\rd{s}\\
\notag &{}\le\frac{1}{\pi}\int_{z^{7/5}}^{\pi}
\exp\biggl\{\frac{1}{2}\sum_{j=1}^{k_1}
\log\Bigl(1-\frac{\re^{-zj}}{2}(1-\cos sj)\Bigr)\biggr\}\,\rd{s}\\
\label{eq:llt1}&\le\frac{2}{\pi}\int_{z^{7/5}/2}^{\pi/2}
\exp\biggl\{-\frac{\re^{-t_1}}{4}\sum_{j=0}^{k_1} (1-\cos
2ju)\biggr\}\,\rd{u},
\end{align}
where the substitution $s=2u$ is made in the last line. The last sum
in \eqref{eq:llt1} can be easily estimated: for
$u\in[0,\frac12\myp\pi]$
\begin{equation}\label{eq:sumcos}
\sum_{j=0}^{k_1}(1-\cos
2ju)=\frac{2k_1+1}{2}-\frac{\sin((2k_1+1)u)}{2\sin u}
\ge\min\left\{\frac{k_1^3u^2}{3},\frac{2k_1+1}{4}\right\},
\end{equation}
where for $u\in[0,\pi/(2k_1+1)]$ the inequality \eqref{eq:sumcos}
follows from the elementary inequalities $u-u^3/6\le \sin u\le u$
and $\sin x\le x-x^3/12$, $x\in[0,\pi]$ (applied with
$x=(2k_1+1)\myp u$), while for $u\in[\pi/(2k_1+1),\pi/2]$
\eqref{eq:sumcos} follows from the inequalities $|\sin x|\le 1$ and
$\sin u\ge 2u/\pi\ge 2/(2k_1+1)$. Hence, for $u\in[\frac12\myp
z^{7/5}\myn,\frac12\myp\pi]$ and small $z>0$, the sum
\eqref{eq:sumcos} is bounded below by $t_1^3z^{-1/5}/12$, and this
estimate combined with \eqref{eq:llt1} yields
\begin{equation}\label{eq:I2}
|I_2|\le \exp\{-t_1^3\mypp\re^{-t_1}z^{-1/5}/48\}.
\end{equation}
Plugging \eqref{eq:I1} and \eqref{eq:I2} in to \eqref{eq:I1+I2} and
using \eqref{eq:zsqrtn} to reformulate the obtained estimate in
terms of $n$ yields the result.
\endproof

\subsection{Proof of Proposition~\ref{pr:mumnum}}
Consider case~(a). Substituting \eqref{eq:kn1} into
\eqref{eq:sksim}, we obtain
\begin{equation}\label{eq:skn}
s_{k_n}=\frac{q\mypp T^{2}}{2\,\vartheta_q(T)^2}\,
n+O(n^{\beta+1})\qquad (n\to\infty),
\end{equation}
where the first term disappears for $q=0$. From \eqref{eq:deftheta}
and \eqref{eq:skn} it follows, for $q\ge 0$,
\begin{equation}\label{eq:n-skn}
n-s_{k_n}\sim n\biggl(1-\frac{q\mypp T^2}{2\vartheta_q(T)^2}\biggr)
=n\,\frac{\Li_2(1-\re^{-T})}{\vartheta_q(T)^2}\qquad (n\to\infty).
\end{equation}

Let $z_n>0$ be the unique solution of the equation
\begin{equation}\label{eq:defzn}
(n-s_{k_n})\myp z^2=\Li_2(1-\re^{-k_nz}).
\end{equation}
Using the asymptotic equations \eqref{eq:n-skn} and \eqref{eq:kn1}
one can verify that the limit
\begin{equation}\label{eq:xi}
\xi:=\lim_{n\to\infty}\frac{z_n\sqrt{n}}{\vartheta_q(T)}
\end{equation}
must satisfy the equation
\[
\xi^2\Li_2(1-\re^{-T})=\Li_2(1-\re^{-T\xi}),
\]
which has the unique root $\xi=1$. As a result, the relation
\eqref{eq:xnT} holds for such $z_n$; it also follows that $z_nk_n\to
T$ as $n\to\infty$.

On the other hand, by Lemma~\ref{lm:muxnm} we obtain, for
$v>\frac{3}{4}$ and large enough~$n$,
\begin{equation}\label{eq:muzn}
\mu^\bq_{z_n,\myp k_n}\{N(\lambda)=n\}\ge n^{-v}.
\end{equation}
Let the event $A_{z,\myp k}$ be as given in Proposition
\ref{pr:mumnum}, then
\[
\nu^\bq_{n,\myp k_n}(A_{z_n,\myp k_n}) =\frac{\mu^\bq_{z_n,\myp
k_n}(A_{z_n,\myp k_n}\cap\{N(\lambda)=n\})} {\mu^\bq_{z_n,\myp
k_n}\{N(\lambda)=n\}} \le \frac{\mu^\bq_{z_n,\myp k_n}(A_{z_n,\myp
k_n})} {\mu^\bq_{z_n,\myp k_n}\{N(\lambda)=n\}}
\]
and an application of \eqref{eq:muzn} and \eqref{eq:muxmcd} with
$z=z_n$ and $k(z_n)=k_n$ readily gives~\eqref{eq:nunmcd}.

Case (b) is considered in a similar manner.
 The
assumption $k(z)=o(z^{-2/(\beta+1)})$ (with $\beta<1$) and
\eqref{eq:knass} imply that $k_n\sim \pi k/\sqrt{6n}
=o(n^{1/(\beta+1)})$ as $n\to\infty$. In turn, it follows from
\eqref{eq:sksim} that $s_{k_n}=o(n)$. Hence, if $z_n>0$ is the
solution of \eqref{eq:defzn} then
$\xi:=\lim_{n\to\infty}z_n\sqrt{n}/\vartheta_0$ with
$\vartheta_0\equiv\vartheta_0(\infty)=\pi/\sqrt{6}$
(see~\eqref{eq:q=0}) satisfies
\[
\xi^2\mypp\vartheta_0^2=\Li_2(1)=\frac{\pi}{\sqrt{6}},
\]
which readily implies that $\xi=1$. The rest of the proof is the
same as for case (a) above.

\subsection{Proof of Proposition~\ref{pr:munu}}
For any $z>0$,
\begin{equation}\label{eq:nulemu}
\nu^\bq_n(A_z)=\frac{\mu^\bq_z\bigl(A_z\cap\{N(\lambda)=n\}\bigr)}{\mu^\bq_z\{N(\lambda)=n\}}
\le \frac{\mu^\bq_z(A_z)}{\mu^\bq_z\{N(\lambda)=n\}}\,.
\end{equation}
The upper bound for the numerator on the right-hand side of
\eqref{eq:nulemu} is guaranteed by condition \eqref{eq:limsupmuz},
and the denominator can be bounded below as follows. Recall that the
measure $\mu^\bq_{z,\myp k}$ is the probability measure $\mu^\bq_z$
conditioned on the event $\{K(\lambda)=k\}$; hence, by the total
probability formula we have
\begin{equation}\label{eq:muxn:t1}
\mu^\bq_{z}\{N(\lambda)=n\} = \sum_{k=0}^\infty \mu^\bq_{z,\myp
k}\{N(\lambda)=n\}\cdot \mu^\bq_z\{K(\lambda)=k\},\qquad n\in\NN_0.
\end{equation}
By virtue of Lemma~\ref{lm:mukseq}(a), $k=k_*\equiv k_*(z)$ defined
in \eqref{eq:def-k0} maximizes $\mu^\bq_{z}\{K(\lambda)=k\}$, and if
$q>0$, Theorem \ref{th:mld} applied with $c=\frac14$ guarantees
that, for any $\gamma\in\bigl(0,\frac12(1-\beta)\bigr)$ and for
$z>0$ small enough,
\begin{equation}\label{eq:muKge}
\mu^\bq_z\{K(\lambda)=k_*\} \ge
\frac{1-\mu^\bq_z\{|K(\lambda)-k_*|>cz^{\gamma-1}\}}{1+2c\myp
z^{\gamma-1}} \sim \tfrac12\myp c^{-1}z^{1-\gamma}=2z^{1-\gamma}.
\end{equation}
If $q=0$ we refer to Theorem \ref{th:mldsublin} instead, which
gives, for any $\gamma\in\bigl(0,\frac12(1-\beta)\bigr)$ and $z>0$
small enough,
\begin{equation}\label{eq:muKgeq=0}
\mu^\bq_z\{K(\lambda)=k_*\} \ge
\frac{1-\mu^\bq_z\{K(\lambda)>k_\gamma\}}{k_\gamma} \ge \tfrac12
z^{2(1-\gamma)},
\end{equation}
because $k_\gamma(z)\le \lceil z^{-2\myp(1-\gamma)}\rceil$
(see~\eqref{eq:kappa_gamma<}).

Let $(z_n)$ be a positive sequence satisfying, for large enough
$n\in\NN$, the inequality
\begin{equation}\label{eq:ineqzn}
\left|n-s_{k_*(z)}-z^{-2}\Li_2(1-\re^{-z\myp k_*(z)})\right| \le
z^{-4/3}.
\end{equation}
It is easy to see that $z_n$ must vanish in the limit as
$n\to\infty$. Solutions of \eqref{eq:ineqzn} exist despite the
discontinuities of the function $z\mapsto s_{k_*(z)}$, because
$k_*(z)$ has unit jumps and, consequently, the condition
\eqref{eq:assq*} and the asymptotic formula \eqref{eq:k0} imply that
the jumps of $s_{k_*(z)}$ are bounded by $Q_{k_*(z)}=O(z^{-1})$ for
$q>0$, while for $q=0$ the upper bound in \eqref{eq:k0sublin} gives
$Q_{k_*(z)}=O\bigl(z^{-\beta}(\log\frac1z)^\beta\bigr)$. Thus, the
left-hand side of \eqref{eq:ineqzn} has discontinuities of order
$O(z^{-1})$ as $z\downarrow0$, which is much smaller than the term
$z^{-4/3}$ on the right-hand side. Furthermore, note that $z_n\mypp
k_*(z_n)\to T_{q}$ (see~\eqref{eq:k0}). Hence, in the same fashion
as in the proof of Proposition \ref{pr:mumnum}, we obtain that due
to \eqref{eq:ineqzn} the limit
$\xi:=\lim_{n\to\infty}z_n\sqrt{n}/\vartheta_q$ satisfies the
equation
\[
\xi^2\mypp \vartheta_q^2-\tfrac12\myp q\myp
T_{q}^2=\Li_2(1-\re^{-T_{q}}).
\]
[For $q=0$, use the values $T_0=\infty$, $qT_q^2|_{q=0}=0$ and
$\vartheta_0=\sqrt{\Li_2(1)}=\pi/\sqrt{6}$ (see~\eqref{eq:q=0}).]
Comparing this with equation \eqref{eq:deftheta}, we conclude that
$\xi=1$, and \eqref{eq:xn} readily follows.

With $z=z_n$ and $k=k_*(z_n)$, the conditions of
Lemma~\ref{lm:muxnm} are satisfied, so \eqref{eq:muxn:t1} and
\eqref{eq:muKge} (or \eqref{eq:muKgeq=0} for $q=0$) yield that, for
any $v>\frac34$ and for $n$ large enough,
\begin{equation}\label{eq:polyn}
\mu^\bq_{z_n}\{N(\lambda)=n\} \ge\mu^\bq_{z_n,\myp
k_*(z_n)}\{N(\lambda)=n\}\cdot \mu^\bq_{z_n}\{K(\lambda)=k_*(z_n)\}
\ge \tfrac12\mypp n^{-v}z_n^\sigma,
\end{equation}
where $\sigma={1-\gamma}$ when $q>0$ and $\sigma=2\myp(1-\gamma)$
when $q=0$. But $z_n\sim\const\cdot n^{-1/2}$, so \eqref{eq:polyn}
provides a lower bound which is polynomial in $n\to\infty$. The
claim of the proposition now follows from the estimates
\eqref{eq:polyn} and~\eqref{eq:nulemu}.

\subsection*{Acknowledgements}
%At an early stage, this research was partially supported by DFG
%Grant 436 RUS 113/722.
%Both authors acknowledge hospitality of the Bielefeld Center for
%Interdisciplinary Research (Zentrum f\"ur interdisziplin\"are
%Forschung, ZiF) during the programme ``Stochastic Dynamics:
%Mathematical Theory and Applications'' (2012). The first-named
%author was also supported by a Leverhulme Research Fellowship in
%2010--2011.
The second-named author was supported by a visiting grant from the
School of Mathematics at the University of Leeds in August 2018. We
thank the anonymous reviewers for their thoughtful comments.

\end{document}